\documentclass[manuscript]{informs3} % for review, not blinded

\usepackage{natbib}
 \NatBibNumeric
 \def\bibfont{\small}%
 \bibpunct[, ]{[}{]}{,}{n}{}{,}%

\usepackage[colorlinks=true,breaklinks=true,bookmarks=true,urlcolor=blue,
     citecolor=blue,linkcolor=blue,bookmarksopen=false,draft=false]{hyperref}
\usepackage{multirow}
\def\EMAIL#1{\href{mailto:#1}{#1}}% When hyperref is used, otherwise outcomment 
\usepackage[mathscr]{eucal}
\usepackage{epsfig,epsf,psfrag}
\usepackage{amssymb,amsfonts,latexsym}
\usepackage{amsmath}
\usepackage{graphicx}
\usepackage{epstopdf}
\usepackage{bm,xcolor,url}
\usepackage{fixltx2e}%ordering of single and double column floats
\usepackage{array}%array and tabular environments
\usepackage{verbatim}
\usepackage{algpseudocode, cite}
\usepackage{algorithm}
\usepackage{verbatim}
\usepackage{textcomp}
\usepackage{mathrsfs}
\usepackage[referable]{threeparttablex}
\usepackage{subfig}
\usepackage{enumerate}
\usepackage{footnote}
\usepackage{enumitem}
\usepackage{xr}
\usepackage{mathtools}
\catcode`~=11 \def\UrlSpecials{\do\~{\kern -.15em\lower .7ex\hbox{~}\kern .04em}} \catcode`~=13 

\allowdisplaybreaks[3]

\newcommand{\tnorm}[1]{{\left\vert\kern-0.25ex\left\vert\kern-0.25ex\left\vert #1 
    \right\vert\kern-0.25ex\right\vert\kern-0.25ex\right\vert}}
\newcommand{\tnormt}[1]{{\vert\kern-0.25ex\vert\kern-0.25ex\vert #1 
    \vert\kern-0.25ex\vert\kern-0.25ex\vert}}

\newcommand{\minus}{\scalebox{0.8}{$-$}}

\newcommand{\norm}[1]{\left\Vert#1\right\Vert}
\newcommand{\normt}[1]{\Vert#1\Vert}
\newcommand{\abst}[1]{\vert#1\vert}
\newcommand{\abs}[1]{\left\lvert#1\right\rvert}

\newcommand{\as}{{\rm \;\;a.s.}}
\newcommand{\nn}{\nonumber}

\newcommand{\defeq}{\triangleq}

\newcommand{\nt}{\addtocounter{equation}{1}\tag{\theequation}} 

\newcommand{\dom}{\mathsf{dom}\,}

\newcommand{\inter}{\mathsf{int}\,}
\newcommand{\bdry}{\mathsf{bd}\,}
\newcommand{\diag}{\mathsf{diag}\,}

\newcommand{\calA}{\mathcal{A}}
\newcommand{\calB}{\mathcal{B}}

\newcommand{\calD}{\mathcal{D}}
\newcommand{\calE}{\mathcal{E}}
\newcommand{\calF}{\mathcal{F}}
\newcommand{\calG}{\mathcal{G}}

\newcommand{\calI}{\mathcal{I}}

\newcommand{\calK}{\mathcal{K}}

\newcommand{\calM}{\mathcal{M}}
\newcommand{\calN}{\mathcal{N}}

\newcommand{\calP}{\mathcal{P}}

\newcommand{\calS}{\mathcal{S}}

\newcommand{\calU}{\mathcal{U}}
\newcommand{\calV}{\mathcal{V}}

\newcommand{\calX}{\mathcal{X}}
\newcommand{\calY}{\mathcal{Y}}
\newcommand{\calZ}{\mathcal{Z}}

\newcommand{\tilcalX}{\widetilde{\calX}}

\newcommand{\barcalX}{\bar{\calX}}

\newcommand{\rmc}{\mathrm{c}}

\newcommand{\rmd}{\mathrm{d}}
\newcommand{\rmD}{\mathrm{D}}

\newcommand{\rmE}{\mathrm{E}}

\newcommand{\rmP}{\mathrm{P}}

\newcommand{\bbD}{\mathbb{D}}
\newcommand{\bbE}{\mathbb{E}}

\newcommand{\bbI}{\mathbb{I}}

\newcommand{\bbN}{\mathbb{N}}

\newcommand{\bbP}{\mathbb{P}}

\newcommand{\bbR}{\mathbb{R}}
\newcommand{\bbS}{\mathbb{S}}

\newcommand{\bbU}{\mathbb{U}}
\newcommand{\bbV}{\mathbb{V}}

\newcommand{\bbX}{\mathbb{X}}
\newcommand{\bbY}{\mathbb{Y}}
\newcommand{\bbZ}{\mathbb{Z}}

\newcommand{\barbbR}{\overline{\bbR}}

\DeclareMathAlphabet{\mathbsf}{OT1}{cmss}{bx}{n}
\newcommand{\rvA}{\mathsf{A}}

\newcommand{\tilf}{\widetilde{f}}

\newcommand{\hatG}{\widehat{G}}

\newcommand{\tilG}{\widetilde{G}}

\newcommand{\tilh}{\widetilde{h}}

\newcommand{\tilp}{\widetilde{p}}

\newcommand{\tilq}{\widetilde{q}}

\newcommand{\hatS}{\widehat{S}}

\newcommand{\hatx}{\hat{x}}

\newcommand{\tilx}{\widetilde{x}}

\newcommand{\haty}{\hat{y}}

\newcommand{\barx}{\overline{x}}
\newcommand{\bary}{\overline{y}}
\newcommand{\barz}{\overline{z}}

\newcommand{\barG}{\overline{G}}

\newcommand{\barM}{\overline{M}}

\newcommand{\barQ}{\overline{Q}}

\newcommand{\barS}{\bar{S}}

\newcommand{\uS}{\underline{S}}

\newcommand{\iid}{i.i.d.\ }

\newcommand{\ceil}[1]{\lceil{#1}\rceil}
\newcommand{\floor}[1]{\lfloor{#1}\rfloor}
\newcommand{\lrangle}[2]{\left\langle{#1},{#2}\right\rangle}
\newcommand{\lranglet}[2]{\langle{#1},{#2}\rangle}

\newcommand{\eqa}{\stackrel{\rm(a)}{=}}

\newcommand{\eqc}{\stackrel{\rm(c)}{=}}

\newcommand{\lea}{\stackrel{\rm(a)}{\le}}
\newcommand{\leb}{\stackrel{\rm(b)}{\le}}
\newcommand{\lec}{\stackrel{\rm(c)}{\le}}
\newcommand{\led}{\stackrel{\rm(d)}{\le}}

\newcommand{\gea}{\stackrel{\rm(a)}{\ge}}
\newcommand{\geb}{\stackrel{\rm(b)}{\ge}}

\DeclareMathOperator{\tr}{tr}

\newcommand{\qednew}{\nobreak \ifvmode \relax \else
      \ifdim\lastskip<1.5em \hskip-\lastskip
      \hskip1.5em plus0em minus0.5em \fi \nobreak
      \vrule height0.75em width0.5em depth0.25em\fi}

\TheoremsNumberedBySection  % (Theorem 1.1, Lema 1.1, Theorem 1.2)

\EquationsNumberedBySection % (1.1), (1.2), ...

\begin{document}

% Outcomment only when entries are known. Otherwise leave as is and default values will be used.
%\setcounter{page}{1}
%\DOI{10.1287/xxxx.0000.0000}%

\RUNAUTHOR{Zhao}
%\RUNAUTHOR{Zhao, Haskell, and Tan}

% Title or shortened title suitable for running heads. 
\RUNTITLE{Accelerated Algorithms for Stochastic Convex-Concave SPP}

\graphicspath{{figures/}}
% Full title. 
\TITLE{Accelerated Stochastic Algorithms for Convex-Concave Saddle-Point Problems}

\ARTICLEAUTHORS{
\AUTHOR{Renbo Zhao}
\AFF{Operations Research Center, Massachusetts Institute of Technology, \EMAIL{renboz@mit.edu}}
%\AUTHOR{William B.\ Haskell}
%\AFF{National University of Singapore, \EMAIL{isehwb@nus.edu.sg}}
%\AUTHOR{Vincent Y.\ F.\ Tan}
%\AFF{National University of Singapore, \EMAIL{vtan@nus.edu.sg}}
} 

\ABSTRACT{
We develop stochastic first-order primal-dual algorithms to solve a class of convex-concave saddle-point problems. When the saddle function is strongly convex in the primal variable, we develop the first stochastic restart scheme for this problem. When the gradient noises obey sub-Gaussian distributions, the oracle complexity of our restart scheme is strictly better than any of the existing methods, even in the deterministic case. Furthermore, for each problem parameter of interest, whenever the lower bound exists, the oracle complexity of our restart scheme is either optimal or nearly optimal (up to a log factor). The subroutine used in this scheme is itself a new stochastic algorithm developed for the problem where the saddle function is non-strongly convex in the primal variable. This new algorithm, which is based on the primal-dual hybrid gradient framework, achieves the state-of-the-art oracle complexity and may be of independent interest.  
%is developed .  
% strongly convex regime %and the new , together with this new algowhenever the lower complexity bound exists in the literature, the complexity obtained by our algorithm is optimal in the non-strongly convex regime and nearly optimal (up to log-factor) . %Additionally, in both regimes, we derive the oracle complexity to obtain an $\epsilon$-duality gap  with high probability, under the sub-Gaussian assumption of the gradient noise. 
}%

\MSCCLASS{Primary: 90C47; secondary: 90C15, 90C25}
\ORMSCLASS{Primary: Programming: Nonlinear: Algorithms;
 secondary: Programming: Stochastic}
%  semi-Markov; programming: infinite dimensional}
\HISTORY{}

\maketitle

% Samples of sectioning (and labeling) in MOOR.
% NOTE: (1) all section levels end with a period,
%       (2) capitalization is as shown (sentence style, not title style).
%
\section{Introduction.}\label{sec:intro} 
Let $\bbX$ and $\bbY$ be two finite-dimensional real normed spaces with dual spaces $\bbX^*$ and $\bbY^*$ %and norms $\norm{\cdot}_{\bbX}$ and $\norm{\cdot}_{\bbY}$ 
respectively. 
Consider the following saddle-point problem (SPP)
\begin{equation}
\min_{x\in\calX}\max_{y\in\calY} \big[S(x,y)\defeq f(x) + g(x) + \Phi(x,y) - J(y)\big],\label{eq:main}
\end{equation}
where %$\calX$ and $\calY$ %
$\calX\subseteq\bbX$ and $\calY\subseteq\bbY$ are nonempty closed and convex sets, % in $\bbX$ and $\bbY$ respectively, %$\rvA:\bbX\to\bbY^*$ is a  non-zero (bounded) linear operator (with adjoint $\rvA^*:\bbY\to\bbX^*$), $\lrangle{\cdot}{\cdot}:\bbY^*\times\bbY\to\bbR$ denotes the duality pairing between $\bbY^*$ and $\bbY$, and 
and the functions $f:\bbX\to\barbbR\defeq(-\infty,+\infty]$, $g:\bbX\to\barbbR$ and $J:\bbY\to\barbbR$ are  convex, closed and proper (CCP). In addition, the function $\Phi:\bbX\times\bbY\to [-\infty,+\infty]$ is convex-concave, i.e., $\Phi(\cdot,y)$ is convex for any $y\in\bbY$ and $\Phi(x,\cdot)$ is concave for any $x\in\bbX$. 
%where $\barbbR\defeq(-\infty,+\infty]$. 
We assume that $f$, $g$, $J$ and $\Phi$ satisfy the following regularity conditions:
\begin{itemize}
\item $f$ is $L$-smooth on $\calX$ (where $L\ge 0$). Specifically, $f$ is differentiable on  $\calX'\supseteq\calX$, where $\calX'$ is an open set in $\bbX$, and its gradient $\nabla f:\bbX\to\bbX^*$ is $L$-Lipschitz on $\calX$, i.e., 
\begin{equation}
\norm{\nabla f(x) - \nabla f(x')}_{\bbX^*}\le L\norm{x-x'}_{\bbX}, \quad\forall\,x,x'\in\calX,
\end{equation}
where $\norm{\cdot}_{\bbX^*}$ and $\norm{\cdot}_{\bbX}$ denote the norms on $\bbX^*$ and $\bbX$ respectively.
%We call 
 %(For notational brevity, 
%(In the sequel, we will simply use $\norm{\cdot}$ and $\norm{\cdot}_*$ to denote the primal and dual norms, whenever the underlying normed space is clear.) 
\item $f$ is $\mu$-strongly convex on $\calX$ (where $\mu\ge 0$), i.e., for any $x,x'\in\calX$,
\begin{equation}
f(x)\ge f(x') + \lrangle{\nabla f(x')}{x-x'} + \frac{\mu}{2}\normt{x-x'}^2. 
\end{equation}
In this work, we will consider both cases where $\mu=0$ and $\mu>0$. 
\item $g$ and $J$ %, we assume that they %are ``simple''. For precise meanings,  
admit tractable Bregman proximal projections on $\calX$ and $\calY$, respectively (see Section~\ref{sec:prelim} for details). Define $\dom g\defeq \{x\in\bbX:g(x)<+\infty\}$ and $\dom J\defeq \{y\in\bbY:J(y)<+\infty\}$. To make~\eqref{eq:main} well-posed, we assume that  $\dom g\cap\calX\ne \emptyset$ and  $\dom J\cap\calY\ne \emptyset$. 
\item $\Phi$ is differentiable on $\calX'\times\calY'$, where $\calY'\supseteq\calY$ is an open set in $\bbY$. For any $(x,y)\in\calX\times\calY$, %$y\in\calY$ and $x\in\calX$, 
denote the gradient of $\Phi(\cdot,y)$ and $\Phi(x,\cdot)$ by $x\mapsto\nabla_x \Phi(x,y)$ and $y\mapsto\nabla_y \Phi(x,y)$ respectively. For all $x,x'\in\calX$ and $y,y'\in\calY$, we assume that there exist constants $L_{xx},L_{yx},L_{yy}\ge 0$ such that 
%\supseteq\calY$ %its gradient $\$
\begin{subequations}
%\begin{equation}
\begin{align}
&\norm{\nabla_x \Phi(x,y)-\nabla_x \Phi(x',y)}_{\bbX^*}\le L_{xx}\norm{x-x'}_{\bbX},\label{eq:Phi_sm(a)} \\
&\norm{\nabla_x \Phi(x,y)-\nabla_x \Phi(x,y')}_{\bbX^*}\le L_{yx}\norm{y-y'}_{\bbY},\label{eq:Phi_sm(b)}\\
&\norm{\nabla_y \Phi(x,y)-\nabla_y \Phi(x',y)}_{\bbY^*}\le L_{yx}\norm{x-x'}_{\bbX},\label{eq:Phi_sm(c)}\\
&\norm{\nabla_y \Phi(x,y)-\nabla_y \Phi(x,y')}_{\bbY^*}\le L_{yy}\norm{y-y'}_{\bbY}.\label{eq:Phi_sm(d)}
\end{align}
%\end{equation}
\end{subequations}
Note that %$L_{xy}=L_{yx}$ and %if the norm on $\bbX\times\bbY$ is defined as $\norm{(x,y)}_{\bbX\times\bbY}\defeq (a\normt{x}_\bbX^p+b\normt{y}_{\bbY}^p)^{1/p}$, for any $a,b\ge 0$, and for any $p\ge 1$ or $p\to+\infty$, 
the gradient operator $(x,y)\mapsto[\nabla_x \Phi(x,y), -\nabla_y \Phi(x,y)]$ is %monotone and 
$M$-Lipschitz on $\calX\times\calY$, where $M\defeq L_{xx} + 2L_{yx} + L_{yy}$.
\end{itemize}
%For the nonsmooth functions 
Based on the assumptions above, %in this paper, 
we aim to design optimal (or nearly optimal) first-order algorithms  that find a saddle-point $(x^*,y^*)\in\calX\times\calY$ of Problem~\eqref{eq:main}, i.e., $(x^*,y^*)$ that satisfies 
\begin{equation}
S(x^*,y)\le S(x^*,y^*)\le S(x,y^*), \quad\forall\,(x,y)\in\calX\times\calY.\label{eq:def_saddle_point}
\end{equation}
For well-posedness, we assume that such a saddle-point exists (see Assumption~\ref{assump:bounded} for conditions that guarantee the existence). 
% and 
%\begin{equation}
%g^*\defeq \inf_{x\in\calX}g(x)>-\infty,\quad J^*\defeq \inf_{y\in\calY}J(y)>-\infty.
%\end{equation}

\subsection{Stochastic first-order oracles.}\label{sec:stoc_oracle}

Since we aim to solve~\eqref{eq:main} via first-order information, we need to properly set up the oracle model. For generality, we do not assume that the exact gradients of $f$, $\Phi(\cdot,y)$ and $\Phi(x,\cdot)$ can be obtained. Rather, we only assume that we have access to the {unbiased} estimators of $\nabla f$, $\nabla \Phi(\cdot,y)$ and $\nabla\Phi(x,\cdot)$ (a.k.a., stochastic gradients), which we denote by $\hat{\nabla} f$, $\hat{\nabla} \Phi(\cdot,y)$ and $\hat{\nabla}\Phi(x,\cdot)$, respectively. In addition, we assume that the gradient noise on $\nabla f$, i.e., $\hat\nabla f-{\nabla} f$, has bounded second moment and we denote this bound as $\sigma_{x,f}^2$. %that is bounded 
Similarly, we also assume that gradient noises on  $\nabla \Phi(\cdot,y)$ and $\nabla\Phi(x,\cdot)$ have bounded second-moments and denote the bounds as    
%$\hat\nabla \Phi(\cdot,y)-\nabla \Phi(\cdot,y)$ and $\hat\nabla\Phi(x,\cdot) - {\nabla}\Phi(x,\cdot)$ have {second-moments} that are bounded by $\sigma^2_{x,f}$, 
$\sigma^2_{x,\Phi}$ and $\sigma^2_{y,\Phi}$, respectively. In some situations, we will further assume that the gradient noises have {sub-Gaussian} distributions. For a formal description of the oracle model and a precise statement of the aforementioned assumptions, readers are referred to Section~\ref{sec:algo_cvx} and Assumption~\ref{assump:noise} respectively. 

Indeed, the oracles described above are standard in the literature on stochastic approximation, which dates back to~\citet{Rob_51} and since then, has become a canonical approach to solve stochastic programming (SP) problems. In the standard SP formulation, the smooth functions $f$ and $\Phi$ are typically represented as expectations (see e.g.,~\citet{Nemi_09}), i.e.,  
\begin{equation}
f(x) \defeq \bbE_{\xi\sim P}[\tilf(x,\xi)] \quad\mbox{and}\quad \Phi(x,y) \defeq \bbE_{\zeta\sim Q}[\widetilde{\Phi}(x,y,\zeta)], \quad \forall\,x\in\bbX, \;y\in\bbY,\label{eq:expectation}
\end{equation}
where $\xi$ and $\zeta$ denote the random variables with distributions $P$ (supported on $\calM$) and $Q$ (supported on $\calZ$) respectively, and the functions $\tilf:\bbX\times\calM\to\bbR$ and  $\widetilde{\Phi}:\bbX\times\bbY\times\calZ\to\bbR$ are such that $f$ and $\Phi$ satisfy the convexity and smoothness assumptions above. In particular, if we take $P=n^{-1}\sum_{i=1}^n \delta_{\xi_i}$ and $Q=m^{-1}\sum_{i=1}^m \delta_{\zeta_i}$, where $\{\xi_i\}_{i=1}^n$ and $\{\zeta_i\}_{i=1}^m$ are deterministic points in $\calM$ and $\calZ$, respectively, and $\delta_{\xi_i}$ denotes the delta measure at $\xi_i$ (and the same for $\delta_{\zeta_i}$), then $f$ and $\Phi$ in~\eqref{eq:expectation} assume finite-sum forms, i.e.,
\begin{equation}
f(x) \defeq \frac{1}{n}\sum_{i=1}^n\tilf(x,\xi_i) \quad\mbox{and}\quad \Phi(x,y) \defeq \frac{1}{m}\sum_{i=1}^m\widetilde{\Phi}(x,y,\zeta_i), \quad \forall\,x\in\bbX, \;y\in\bbY.\label{eq:finite-sum}
\end{equation}
In this case, we can construct the stochastic (first-order) oracle for $f$ by first sampling an index set $\calB$ from $[n]\defeq\{1,\ldots,n\}$ uniformly randomly, and then output the gradient of $f_{\calB}\defeq\abs{\calB}^{-1}\sum_{i\in\calB}\tilf(x,\xi_i)$. The stochastic oracle for $\Phi$ can also be constructed in the same way.  
%  according to the uniform distribution over $[n]$ or $[m]$, 
%with that generate stochasticity in the stochastic gradients. % 

%Talk about empirical distribution and SAA here...

From this point on, we will refer to the class of problems in~\eqref{eq:main} as SPP($L$, $L_{xx}$, $L_{yx}$, $L_{yy}$, $\sigma$, $\mu$), where $\sigma \defeq \max\{\sigma_{x,f},\sigma_{x,\Phi},\sigma_{y,\Phi}\}$ represents the collective stochasticity in the (stochastic) gradients of $f$, $\Phi(\cdot,y)$ and $\Phi(x,\cdot)$. If $\sigma=0$, then~\eqref{eq:main} corresponds to a deterministic optimization problem.

{\bf Oracle complexity.} 
To obtain an $\epsilon$-duality gap in expectation or with high probability (cf.\ Section~\ref{sec:conv_res}), we will focus on analyzing the oracle complexity in terms of its dependence on $L$, $L_{xx}$, $L_{yx}$, $L_{yy}$, $\sigma_{x,f}$, $\sigma_{x,\Phi}$, $\sigma_{y,\Phi}$, $\mu$ and $\epsilon$. 
In addition, for the algorithms in this work and almost all the works in the literature, the number of calls to each of the oracle described above (which returns $\nabla f$, $\nabla \Phi(\cdot,y)$ or $\nabla\Phi(x,\cdot)$ or their stochastic versions) is the same. Therefore, in our complexity analysis and  comparison of complexities with other algorithms, we do not distinguish among the these oracles. Instead, the word ``oracle complexity'' refers to the complexity of each of them.

\subsection{Related work.}\label{sec:related_work}
%Before we review the relevant literature, we first briefly introduce the applications of 
The SPP in~\eqref{eq:main} has a wide range of applications across many fields, including statistics, machine learning, operations research and game theory. When $\Phi$ is bilinear, i.e., there exists a (bounded) linear operator $\rvA:\bbX\to\bbY^*$ such that $\Phi(x,y)=\lranglet{\rvA x}{y}$ (where $\lrangle{\cdot}{\cdot}:\bbY^*\times\bbY\to\bbR$ denotes the duality pairing between $\bbY^*$ and $\bbY$), the applications of~\eqref{eq:main} can be found in numerous previous works, e.g.,~\citet{Juditsky_12a,Juditsky_12b},~\citet{Chambolle_16} and~\citet{Zhao_19}. Beyond bilinear $\Phi$, there are also rich applications, including two-player convex-concave zero-sum game~(\citet{Chen_17}), convex optimization with functional constraints~(\citet{Boyd_04}), and kernel matrix learning~(\citet{Lanc_04}). Other  applications can also be found in~\citet{Bala_16}.  
%. For details, we refer readers to~\citet{Bala_16} and~\citet{Hien_17}. 

%\begin{itemize}[labelindent=0pt]
%\item {\em Stochastic programming with expectation constraints.} Let $f$ be defined in~\eqref{eq:expectation} and consider the following convex programming problem: 
%\begin{align}
%\min_{x\in\calX} f(x) %\big\{f(x)=\bbE_{\xi\sim P}[\tilf(x,\xi)]\big\} 
%\quad \st\quad \big\{\psi_i(x)\defeq\bbE_{\zeta\sim Q}[\tpsi_i(x,\zeta)]\big\} \le 0,\quad\forall\,i\in[m],
%\end{align}
%where $\psi_i:\bbR$ is 
%\item {\em Two-player nonlinear zero-sum game.}
%\end{itemize}

%Due to its generality, Problem~\ref{eq:main} subsumes many SPPs studied in the literature as special cases
%Next, we will review the 
Previous works on solving convex-concave SPPs generally fall into two categories, depending on whether the primal-dual coupling term $\Phi$ is bilinear. Our focus will be on the non-bilinear SPPs (i.e., the problem where $\Phi$ is non-bilinear).  Before doing that, we briefly review the works for bilinear SPPs first. %(i.e., the problem where $\Phi$ is bilinear), and then focus on the works for .

\subsubsection{Bilinear SPPs.}\label{sec:bilinear_SPP}
This class of problems is indeed a special case of~\eqref{eq:main}, i.e., when $L_{xx}=L_{yy}=0$. In recent years, both deterministic (i.e., $\sigma=0$) and stochastic (i.e., $\sigma>0$) versions of this problem have been thoroughly studied,  for both  $\mu=0$  and $\mu>0$. 
For the deterministic problems, some well-known algorithms include Nesterov smoothing (a.k.a., excessive gap technique,~\citet{Nest_05,Nest_05b}), primal-dual hybrid gradient (PDHG,~\citet{Chambolle_11,Chambolle_16}), hybrid proximal extragradient-type algorithm (HPE-type,~\citet{He_16}) and primal-dual operator splitting (e.g.,~\citet{Condat_13},~\citet{Vu_13} and~\citet{Davis_15}). In addition, to tackle the stochastic problems, stochastic versions of these algorithms have also been developed, e.g.,~\citet{Chen_14},~\citet{Zhao_18} and~\citet{Zhao_19}. 
%As a remark, all of the abovementioned algorithms crucially leverage the bilinear structure of $\Phi$, and hence cannot be straightforwardly extended to solve non-bilinear SPPs. 

\subsubsection{Non-bilinear SPPs.}
We first consider the case where $\mu=0$ and $\sigma=0$ (i.e., no primal strong convexity and all oracles are deterministic). 
When $\Phi$ is possibly nonsmooth, algorithms based on primal-dual subgradient have been developed in several works, including~\citet{Nedic_09},~\citet{Nest_09} and~\citet{Juditsky_12a}. However, these methods typically incur high oracle complexity, i.e., $O(\epsilon^{-2})$ (where $\epsilon$ denotes the desired accuracy for the duality gap). As a result, they are not competitive when $\Phi$ is smooth (cf.~\eqref{eq:Phi_sm(a)} to~\eqref{eq:Phi_sm(d)}). The smoothness of $\Phi$ has been exploited in many algorithms to achieve better complexity results. These methods include Mirror-Prox (\citet{Nemi_05}), HPE-type algorithm (\citet{Kolo_17}) and PDHG-type algorithm (\citet{Hamed_18}). In particular, the last two algorithms are the extensions of their counterparts for solving bilinear SPPs. When $\sigma>0$, stochastic extensions of Mirror-Prox have been developed in the literature. Some representative works include the stochastic Mirror-Prox (SMP) method (\citet{Juditsky_11}) and the stochastic accelerated Mirror-Prox (SAMP) method (\citet{Chen_17}). % and~\citet{Hamed_18b}.

Unlike the case where $\mu=0$, there exist very few works that have considered the case where $\mu>0$. When $L_{yy}=0$ and $\sigma=0$ (i.e., the function $\Phi(x,\cdot)$ is linear and all oracles are deterministic),~\citet{Juditsky_12b} and~\citet{Hamed_18} have proposed algorithms, which are based on Mirror-Prox and PDHG respectively, that achieve better complexity than their counterparts that are designed for  $\mu=0$. However, no algorithms exist when $L_{yy}=0$ or $\sigma=0$.  Therefore, in this work, we seek to resolve two important questions:
\vspace{1ex} 
\begin{enumerate}[labelindent=0pt,label=(\Roman*)]
\item \underline{\em Can we improve the oracle complexities of these two algorithms when $L_{yy}=0$ and $\sigma=0$?} \label{item:Q1}
\item \underline{\em Can we develop an algorithm that works for all the cases where $\mu>0$, $L_{yy}>0$ and $\sigma>0$?}\label{item:Q2}
\end{enumerate}
\vspace{1ex}
Indeed, we will provide affirmative answers to both questions above, by developing a stochastic restart scheme that is not only able to deal all with the cases listed in~\ref{item:Q2}, but also significantly improves the oracle complexities of the algorithms in~\citep{Juditsky_12b} and~\citep{Hamed_18}, even in the case where $L_{yy}=0$ and $\sigma=0$. 
%enjoys the oracle complexity that nearly matches the existing lower bounds. 

\subsection{Main contributions.}\label{sec:main_contrib}

Our main contributions are summarized below.

\begin{enumerate}[labelindent=0pt,label=(\Roman*)]

\item 
First, when $\mu>0$ (i.e., $f$ is strongly convex), we develop a  (multi-stage) {stochastic} restart scheme (i.e., Algorithm~\ref{algo:restart_stoc}) for solving~\eqref{eq:main}. %by assuming that the gradient noises obey sub-Gaussian distributions. % by using a modified version of Algorithm~\ref{algo:convex_f} (developed for the case where $\mu=0$) as the subroutine. 
Since SPPs have different structures from convex optimization problems (COPs), our stochastic restart scheme is different from those for COPs (e.g.,~\citet{Ghad_13b}). Specifically, we use a distance-based quantity as the restart criterion, instead of the objective error. In addition, instead of focusing on expectation, we analyze the stochasticity via the error probability, which is obtained using techniques from finite-state Markov chains. (For detailed discussions, we refer readers to Section~\ref{sec:restart}.) In principle, the subroutine used in our restart scheme can be any stochastic algorithm developed for the SPP in~\eqref{eq:main} when $\mu=0$ (i.e., $f$ is non-strongly convex). However, different subroutines may result in different oracle complexities. Therefore, to achieve the desired oracle complexity, we develop a new stochastic algorithm (i.e., Algorithm~\ref{algo:convex_f}) for the case where $\mu=0$, based on the PDHG framework. (For details, see contribution~\ref{item:mu_0} below.)%This algorithm is of independent interest. 

%interface with other subroutines (e.g., those based on Mirror-Prox), we believe that our restart scheme is of independent interest for solving stochastic SPPs with primal strong convexity. 

%When the gradient noises have bounded second moments and follow sub-Gaussian distributions, 
When the gradient noises obey sub-Gaussian distributions, to achieve an $\epsilon$-duality gap with probability  (with probability) at least $1-\varsigma$, our scheme has oracle complexity  
\begin{equation}
O\left(\bigg(\sqrt{\frac{L}{\mu}}+\frac{L_{xx}}{\mu}\bigg)\log\left(\frac{1}{\epsilon}\right)+ \frac{L_{yx}}{\sqrt{\mu\epsilon}} + \frac{L_{yy}}{\epsilon}+ \left(\frac{(\sigma_{x,f}+\sigma_{x,\Phi})^2}{\mu\epsilon}+\frac{\sigma_{y,\Phi}^2}{\epsilon^2}\right)\log\left(\frac{\log(1/\epsilon)}{\varsigma}\right)\right).\label{eq:complexity_sc_hp_intro}
\end{equation}
Note that even in the deterministic case (i.e., $\sigma=0$), this complexity is strictly better than any in the previous works (cf.~Table~\ref{table:sc}). Based on the complexity in~\eqref{eq:complexity_sc_hp_intro}, under very mild conditions on the nonsmooth functions $g$ and $J$, % (cf.\ Assumption~\ref{assump:regularity_g_J}), 
our scheme obtains an $\epsilon$-expected duality gap with oracle complexity 
\begin{equation}
O\left(\bigg(\sqrt{\frac{L}{\mu}}+\frac{L_{xx}}{\mu}\bigg)\log\left(\frac{1}{\epsilon}\right)+ \frac{L_{yx}}{\sqrt{\mu\epsilon}} + \frac{L_{yy}}{\epsilon}+ \left(\frac{(\sigma_{x,f}+\sigma_{x,\Phi})^2}{\mu\epsilon}+\frac{\sigma_{y,\Phi}^2}{\epsilon^2}\right)\log\left(\frac{1}{\epsilon}\right)\right).\label{eq:complexity_sc_exp_intro}
\end{equation}
%Compared to the complexity in~\eqref{eq:complexity_cvx_f_intro}, we observe that the complexities of $L$, $L_{xx}$, $L_{yx}$ and $\sigma_{x,f}+\sigma_{x,\Phi}$ have been greatly improved, due to primal strong convexity. (For detailed discussions, we refer readers to Section~\ref{sec:comp_sc}.) 
Regarding the optimality of the complexity in~\eqref{eq:complexity_sc_exp_intro}, the complexities of $L_{xx}$ and $L_{yx}$ match the lower bounds derived in~\citet{Nemi_79} and~\citet{Ouyang_18}, respectively. Additionally, the complexities of $\sigma_{x,f}+\sigma_{x,\Phi}$ and $\sigma_{y,\Phi}$ nearly match (up to $\log(1/\epsilon)$ factor) the lower bounds derived in~\citet{Raginsky_11}. %Similar to the case where $\mu=0$, 
The complexities of $L_{xx}$ and $L_{yy}$ are the best-known, although their lower bounds are not known, to the best of our knowledge. 
%not available in the literature. %, but these complexities are indeed the best-known. 

\begin{table}[t]\centering
%\TABLE
\caption{Comparison of oracle complexities of algorithms to obtain an $\epsilon$-expected duality gap when $\mu>0$.\label{table:sc}}
{\begin{tabular}{ccc}\hline
Algorithm & Problem Class & Oracle Complexity\\\hline
PDHG-type~\citep{Hamed_18} & SPP($L$, $L_{xx}$, $L_{yx}$, 0, 0, $\mu$) & $O\left(\frac{L+L_{xx}+L_{yx}}{\sqrt{\mu\epsilon}}\right)$ \\\hline
Mirror-Prox-B~\citep{Juditsky_12b} & SPP($L$, $L_{xx}$, $L_{yx}$, 0, 0, $\mu$) & $O\left(\frac{L+L_{xx}}{\mu}\log\big(\frac{1}{\epsilon}\big) + \frac{L_{yx}}{\sqrt{\mu\epsilon}}\right)$ \\\hline
Algorithm~\ref{algo:restart_stoc} & SPP($L$, $L_{xx}$, $L_{yx}$, $L_{yy}$, $\sigma$, $\mu$) & \eqref{eq:complexity_sc_exp_intro}\\\hline
\end{tabular}}
\end{table}

\item \label{item:mu_0}
Second, when $\mu=0$ (i.e., $f$ is non-strongly convex), we develop a stochastic algorithm (i.e., Algorithm~\ref{algo:convex_f}) that is used as a subroutine in our restart scheme (see above), and also may be of independent interest by itself. We develop this algorithm by extending the PDHG framework, which was originally developed for bilinear SPPs (cf.\ Section~\ref{sec:bilinear_SPP}), to handle the non-bilinear case. In addition, we incorporate the stochastic acceleration technique (see e.g.,~\citet{Lan_12}) into our algorithm, which indeed enables us to obtain the optimal oracle complexity for the smooth function $f$. By judicious choices of the algorithm parameters,  we are able to obtain an $\epsilon$-expected duality gap with the {state-of-the-art} oracle complexity (cf.~Table~\ref{table:cvx}) %is %of Algorithm~\ref{algo:convex_f} is
\begin{equation}
O\left(\sqrt{\frac{L}{\epsilon}} + \frac{L_{xx} + L_{yx} + L_{yy}}{\epsilon} + \frac{(\sigma_{x,f}+\sigma_{x,\Phi})^2+\sigma_{y,\Phi}^2}{\epsilon^2}\right), \label{eq:complexity_cvx_f_intro}
\end{equation}
when all the gradient noises have bounded second moments (but the distributions are not necessarily sub-Gaussian). 
Previously, the complexity in~\eqref{eq:complexity_cvx_f_intro} has been achieved by the SAMP algorithm introduced in~\citet{Chen_17}. %This complexity bound is the state-of-the-art (cf.\ Table~\ref{table:cvx}), and has been achieved previously by~\citet{Chen_17}. 
However, since the algorithm in~\citep{Chen_17} is based on Mirror-Prox, it significantly differs from our algorithm (which is based on PDHG). Consequently, our algorithm provides a valuable alternative approach to achieve the complexity in~\eqref{eq:complexity_cvx_f_intro}.  
%As a result, our algorithm is significantly different from theirs. 

Regarding the optimality of the complexity in~\eqref{eq:complexity_cvx_f_intro}, we notice that the complexities for $L$ and $L_{yx}$ match the lower bounds derived in~\citet{Ouyang_18}, and the complexities for $\sigma_{x,f}+\sigma_{x,\Phi}$ and $\sigma_{y,\Phi}$ match the lower bounds derived in~\citet{Nemi_79}. Therefore, all of these complexities are optimal. In addition, the complexities of $L_{xx}$ and $L_{yy}$ are also the best-known, although no lower bounds have been derived in the literature. ({Note that in this case, by taking $\Phi(x,y)=\phi^\rmP(x) + \lranglet{\rvA x}{y} - \phi^\rmD(y)$, where the convex functions $\phi^\rmP$ and $\phi^\rmD$ are  $L_{xx}$- and $L_{yy}$-smooth respectively,  we can obtain lower complexity bounds $O(\sqrt{L_{xx}/\epsilon})$ and $O(\sqrt{L_{yy}/\epsilon})$. However, this essentially returns to the bilinear case. Therefore, these lower bounds may not be tight for the non-bilinear case.}) % where $x$ and $y$ are coupled in a } 

For the analysis of Algorithm~\ref{algo:convex_f}, % as a subroutine in our stochastic restart scheme, 
we also derive the following {\em large-deviation-type} convergence result, which complements the convergence result in expectation (cf.~\eqref{eq:complexity_cvx_f_intro}). Specifically, if all the gradient noises obey sub-Gaussian distributions, we can obtain an $\epsilon$-duality gap with probability at least $1-\varsigma$ with oracle complexity %, again, with the state-of-the-art oracle complexity %is %of Algorithm~\ref{algo:convex_f} is
\begin{equation}
O\left(\sqrt{\frac{L}{\epsilon}} + \frac{L_{xx} + L_{yx} + L_{yy}}{\epsilon} + \frac{(\sigma_{x,f}+\sigma_{x,\Phi})^2+\sigma_{y,\Phi}^2}{\epsilon^2}\log\bigg(\frac{1}{\varsigma}\bigg)\right). \label{eq:complexity_cvx_f_hp_intro}
\end{equation}
%This  is crucial to interface with 
%In particular, the $\log(1/\varsigma)$ factor in~\eqref{eq:complexity_cvx_f_hp_intro} indicates that our algorithm achieves the  (see e.g.,~\citet{Nemi_09}). 

\renewcommand{\arraystretch}{1.2}
\setlength{\tabcolsep}{8pt}
\begin{table}[t]\centering
%\TABLE
\caption{Comparison of oracle complexities of algorithms to obtain an $\epsilon$-expected duality gap when $\mu=0$. \label{table:cvx}}
\begin{threeparttable}
{\begin{tabular}{ccc}\hline
Algorithm\tnotex{ftn:HPE} & Problem Class & Oracle Complexity\\\hline
PDHG-type~\citep{Hamed_18} & SPP($L$, $L_{xx}$, $L_{yx}$, $L_{yy}$, 0, 0) & $O\left({\frac{L}{\epsilon}} + \frac{L_{xx} + L_{yx} + L_{yy}}{\epsilon} \right)$\\\hline
MP~\citep{Nemi_05} & SPP($L$, $L_{xx}$, $L_{yx}$, $L_{yy}$, 0, 0) & $O\left({\frac{L}{\epsilon}} + \frac{L_{xx} + L_{yx} + L_{yy}}{\epsilon} \right)$\\\hline
SMP~\citep{Juditsky_11} & SPP($L$, $L_{xx}$, $L_{yx}$, $L_{yy}$, $\sigma$, 0) & $O\left({\frac{L}{\epsilon}} + \frac{L_{xx} + L_{yx} + L_{yy}}{\epsilon} + \frac{(\sigma_{x,f}+\sigma_{x,\Phi})^2+\sigma_{y,\Phi}^2}{\epsilon^2}\right)$\\\hline
SAMP~\citep{Chen_17} & SPP($L$, $L_{xx}$, $L_{yx}$, $L_{yy}$, $\sigma$, 0) & $O\left(\sqrt{\frac{L}{\epsilon}} + \frac{L_{xx} + L_{yx} + L_{yy}}{\epsilon} + \frac{(\sigma_{x,f}+\sigma_{x,\Phi})^2+\sigma_{y,\Phi}^2}{\epsilon^2}\right)$\\\hline
Algorithm~\ref{algo:convex_f} & SPP($L$, $L_{xx}$, $L_{yx}$, $L_{yy}$, $\sigma$, 0) & $O\left(\sqrt{\frac{L}{\epsilon}} + \frac{L_{xx} + L_{yx} + L_{yy}}{\epsilon} + \frac{(\sigma_{x,f}+\sigma_{x,\Phi})^2+\sigma_{y,\Phi}^2}{\epsilon^2}\right)$\\\hline
\end{tabular}}
\begin{tablenotes}\footnotesize
\item[1] \label{ftn:HPE} {Note that we exclude the HPE-type algorithm (\citet{Kolo_17}) from comparison since in their complexity analysis, a different convergence criterion from the duality gap is used. Moreover, only ``inner iteration'' complexity is analyzed in~\citep{Kolo_17}, which is lower than the actual oracle complexity.} 
\end{tablenotes}
\end{threeparttable}%\vspace{-.5cm}
\end{table}

\end{enumerate}

%Finally, the oracle complexities are compared in Table~\ref{table:cvx} and Table~\ref{table:sc}.

%\newpage
\subsection{Notation.}\label{sec:notation}
Denote the set of natural numbers by $\bbN\defeq\{1,2,\ldots\}$ and define $\bbZ_+\defeq \bbN\cup\{0\}$. 
For any finite-dimensional real normed space $\bbU$, we denote its dual space by $\bbU^*$. We denote the norms on $\bbU$ and $\bbU^*$ by $\norm{\cdot}$ and $\norm{\cdot}_*$ respectively. 
In addition, denote the duality pairing between $\bbU^*$ and $\bbU$ by $\lrangle{\cdot}{\cdot}:\bbU^*\times\bbU\to\bbR$. 
For any CCP function $h:\bbU\to\barbbR$, define its domain as $\dom h\defeq\{u\in\bbU:h(u)<+\infty\}$. In addition, for any nonempty set $\calU$ in $\bbU$, denote its interior by $\inter \calU$, its boundary by $\bdry \calU$ and its diameter by $D_\calU\defeq \sup_{u,u'\in\calU}\norm{u-u'}$. We also denote the indicator function of $\calU$ as $\iota_\calU:\bbU\to\barbbR$, i.e., $\iota_\calU(u)\defeq 0$ if $u\in\calU$ and $\iota_\calU(u)\defeq +\infty$ otherwise.  Note that $\iota_\calU$ is closed and convex if and only if $\calU$ is closed and convex. 
%\begin{equation}
%\begin{cases}
%\end{cases}
%\end{equation}

\subsection{Organization.}
The paper is organized as follows. We first introduce some preliminary material in Section~\ref{sec:prelim}. Then, we present Algorithm~\ref{algo:convex_f} developed for the case where $f$ is non-strongly convex (i.e., $\mu=0$), together with convergence results. The detailed proofs of these results are deferred to Section~\ref{sec:analysis}. In Section~\ref{sec:f_sc}, we develop our stochastic restart scheme for the case where $f$ is strongly convex (i.e., $\mu>0$), by using a modified version of Algorithm~\ref{algo:convex_f} as a subroutine. We analyze the oracle complexity for this scheme and  discuss some technical issues. 

\section{Preliminaries.}\label{sec:prelim}
We first introduce the distance generating function and Bregman proximal projection, followed by the primal and dual functions associated with $S(\cdot,\cdot)$ in~\eqref{eq:main}. 

\subsection{Distance generating function and Bregman proximal projection.}\label{sec:DGF_BPP}
Let $\bbU$ and $h$ be given as in Section~\ref{sec:notation}. 
We say that $h$ is {\em essentially smooth} if $h$ is continuously differentiable on $\inter \dom h\ne \emptyset$ and for any %$u\in \dom h\setminus\inter\dom h$ 
$u\in\bdry\dom h$ and any sequence $\{u^k\}_{k\in\bbZ^+}\subseteq\inter\dom h$  such that $u^k\to u$, $\normt{\nabla h(u^k)}_*\to+\infty$. 
Let $\calU$ be any nonempty closed and convex set in $\bbU$. %define $\calU^o\defeq \calU\cap \inter\dom h$. 
 We call $h_\calU$ a {\em distance generating function} (DGF) on $\calU$ if it is continuous on $\calU$,  essentially smooth, %(hence $\calU\subseteq\dom h_\calU$) 
 and% there exists $1>0$ such  that %$h$ to be $a_h$-strongly convex on $$  
\begin{equation}
D_{h_\calU}(u,u')\defeq h_\calU(u) - h_\calU(u') - \lranglet{\nabla h_\calU(u')}{u-u'}\ge (1/2)\norm{u-u'}^2, \,\forall\,u\in\calU,\,\forall\,u'\in\calU^o, \label{eq:quad_lower_bound}
\end{equation}
where $\calU^o\defeq \calU\cap \inter\dom h_\calU$ and $D_{h_\calU}:\calU \times\calU^o\to\bbR$ is called the Bregman distance associated with $h_\calU$. %Note that~\eqref{eq:quad_lower_bound} in particular implies the $1$-strong convexity of $h_\calU$ on $\calU^o$. 
Based on $D_{h_\calU}(\cdot,\cdot)$, we define the {\em Bregman diameter} of $\calU$ under $h_\calU$  as
\begin{equation}
\Omega_{h_\calU}\defeq \sup\nolimits_{u\in\calU,u'\in\calU^o}D_{h_\calU}(u,u'). \label{eq:Bregman_diam}
\end{equation}

In addition, for any $u'\in\calU^o$ and CCP function $\varphi:\bbU\to\barbbR$, define the {\em Bregman proximal projection} (BPP) of $u'$ on $\calU$ under $\varphi$ (associated with the DGF $h_\calU$, $u^*\in\bbU^*$ and $\lambda>0$) as 
\begin{equation}
u'\mapsto u^+\defeq \argmin_{u\in\calU} \big[\varphi_\lambda(u)\defeq\varphi(u) + \lrangle{u^*}{u} + \lambda^{-1}D_{h_\calU}(u,u')\big].\label{eq:Bregman_projection}
\end{equation}
Note that if $\inf_{u\in\calU}\varphi(u)>-\infty$ and $\calU^o\cap\dom \varphi\ne \emptyset$, then the minimization problem in~\eqref{eq:Bregman_projection} always has a unique solution in $\calU^o\cap\dom\varphi$  (see Lemma~\ref{lem:unique_sln} in Appendix~\ref{app:unique}). %the mapping in is single-valued  
We say that the function $\varphi$ has a {\em tractable} BPP on $\calU$ if there exists a DGF $h_\calU$ on $\calU$ such that the minimization problem in~\eqref{eq:Bregman_projection} has a unique {and easily computable} solution in $\calU^o\cap\dom\varphi$, for any $u^*\in\bbU^*$ and $\lambda>0$. 

We now provide some examples of the quadruple $(\bbU,\varphi,\calU,h_\calU)$ that satisfy the assumptions above. Clearly, if $\bbU$ is a Hilbert space with inner product $\lranglet{\cdot}{\cdot}$ and its induced norm $\norm{\cdot}$, and $\calU=\bbU$, we can take $h_\calU=(1/2)\norm{\cdot}^2$, so that the minimization problem in~\eqref{eq:Bregman_projection} becomes the usual proximal minimization problem associated with the function $\varphi$. Due to this, we will state three examples below where $\bbU$ has non-Hilbertian geometry and $\varphi\equiv 0$ (or equivalently, $\varphi=\iota_\calU$ and the minimization problem in~\eqref{eq:Bregman_projection} becomes unconstrained). 
These examples can be found in several references, e.g.,~\citet{Nest_05,Juditsky_12a}. %For convenience, let $\mathsf{dim}\, \bbU=n$. 

\begin{enumerate}[label={\rm (E\arabic*)}]
\item\label{item:ell_1_entropy} We let $\bbU=(\bbR^n,\norm{\cdot}_1)$, where $\norm{u}_1 \defeq \sum_{i=1}^n \abs{u_i}$, $\calU = \Delta_n\defeq \{u\in\bbR^n:u\ge 0, \sum_{i=1}^n u_i=1\}$ and $h_\calU(u)=\sum_{i=1}^n u_i\ln u_i$. It is well-known that $h_\calU$ is 1-strongly convex on $\Delta_n$ with respect to $\norm{\cdot}_1$ (see e.g.,~\citet{Beck_03}), so the inequality in~\eqref{eq:quad_lower_bound} is satisfied. In addition, the minimization problem in~\eqref{eq:Bregman_projection} has closed-form solution, namely, %for any $i\in[n]$, 
\begin{equation}
u^+_i =  \frac{\exp(\xi_i)}{\sum_{j=1}^n \exp(\xi_j)}, \quad\;\forall\,i\in[n], %\quad\mbox{where}\;\;\xi_i = u^*_i - 1-\ln u'_i, . 
\end{equation}
where $\xi_i \defeq 1+\ln u'_i - \lambda u^*_i$, for any $i\in[n]$. 
\item This example is the matrix analogy of~\ref{item:ell_1_entropy}. Let $\bbU = (\bbS^n, \norm{\cdot}_*)$, where $\bbS^n$ denotes the space of real symmetric $n\times n$ matrices, $\norm{U}_* \defeq \norm{\sigma(U)}_1$ and $\sigma(U)\defeq(\sigma_1(U),\ldots,\sigma_n(U))\in\bbR^n$ denotes the vector of eigenvalues of $U$. Note that $\bbU^* = (\bbS^n,\norm{\cdot}_*)$, where $\normt{U^*}_* \defeq \max_{i\in[n]}\abs{\sigma_i(U^*)}$. In addition, let $\calU = \{U\in\bbS^n:U\succeq 0,\tr(U)=1\} = \{U\in\bbS^n:\sigma(U)\ge 0,\textstyle{\sum}_{i=1}^n\sigma_i(U)=1\}$, 
%\begin{equation}
%\calU = \{U\in\bbS^n:U\succeq 0,\tr(U)=1\} = \{U\in\bbS^n:\sigma(U)\ge 0,\textstyle{\sum}_{i=1}^n\sigma_i(U)=1\},\nn
%\end{equation}
 i.e., the set of positive semi-definite matrices with unit trace %. We choose 
 and $h_\calU(U) = \sum_{i=1}^n \sigma_i(U)\ln \sigma_i(U)$. %, whose gradient $\nabla h_\calU(U)$ 
 Consider 
 \begin{align*}
 h_\calU^*(U^*)&\defeq {\sup}_{U\in\calU} \lranglet{U^*}{U} - h_\calU(U) \\
 &= {\sup}_{\sigma(U)\in\Delta_n} \lranglet{\sigma(U^*)}{\sigma(U)} - \textstyle{\sum}_{i=1}^n \sigma_i(U)\ln \sigma_i(U) \\
 &= \ln\left(\textstyle{\sum}_{i=1}^n\exp(\sigma_i(U^*))\right). 
 \end{align*}
 In~\citet[Section~2.2]{Nest_07}, it is shown that $h_\calU^*$ is differentiable and $\nabla h_\calU^*$ is $1$-Lipschitz with respect to $\norm{\cdot}_*$ on $\bbS^n$. Therefore, we conclude that $h_\calU$ is 1-strongly convex with respect to $\norm{\cdot}$ on $\calU$.
Additionally, the solution of the minimization problem in~\eqref{eq:Bregman_projection} admits closed form. Specifically, let $\Lambda\defeq \nabla h_\calU(U') - \lambda U^*\in\bbS^n$ (see~\citet[Section~6.2]{Ben-Tal_05} for the expression of $\nabla h_\calU(U')$), whose eigen-decomposition is denoted by $P\diag(\sigma(\Lambda))P^T$ (where $P$ is an orthogonal matrix). Then $U^+ = P\diag(\sigma(U^+)) P^T$, where 
%$\sigma_i(U^+) = \exp(\sigma_i(\Xi))/\sum_{j=1}^n\exp(\sigma_j(\Xi))$, for any $i\in[n]$. 
\begin{align*}
\sigma_i(U^+) = \frac{\exp(\sigma_i(\Lambda))}{\sum_{j=1}^n\exp(\sigma_j(\Lambda))}, \quad\;\forall\,i\in[n]. 
\end{align*}
\item Finally, consider $\bbU=(\bbR^n,\norm{\cdot}_p)$, where $p\in(1,2]$ and $\norm{u}_p \defeq \left(\sum_{i=1}^n \abs{u_i}^p\right)^{1/p}$. Consequently, $\bbU^*=(\bbR^n,\norm{\cdot}_q)$, where $q \defeq 1/(1-p^{-1})\in[2,+\infty)$. Let $\calU = \bbR_+^n$, where $\bbR_+\defeq [0,+\infty)$ and $h_\calU(u) = (1/2)\norm{u}_p^2$. From~\citet[Section~8]{Ben-Tal_01}, we know that $h_\calU$ is $(p-1)$-strongly convex with respect to $\norm{\cdot}_p$ on $\bbR^n$. In addition, the minimization problem in~\eqref{eq:Bregman_projection} can be reduced to a one-dimensional optimization problem (which is easy to solve). To see this, from the KKT conditions,  if $u^*_i\ge 0$, then $u^+_i=0$. Define $\calI^{\minus}\defeq \{i\in[n]:u^*_i<0\}$. If $\abst{\calI^{\minus}}\le 1$, the claim is trivially true. Otherwise, fix any $i\in\calI^{\minus}$ and for any $j\in\calI^{\minus}\setminus\{i\}$, we have that $u_j = u_i/(u^*_i/u^*_j)^{q-1}$. Thus the minimization problem in~\eqref{eq:Bregman_projection} still only involves one variable $u_i$. % to solve. 
\end{enumerate}
%In many scenarios, depending on the 
For many constraint sets $\calU$, compared with Hilbertian geometry, working under non-Hilbertian geometry offers better dependence on the dimension of $\bbU$ in the oracle complexities of several first-order algorithms (e.g., the subgradient and gradient methods). For detailed discussions, see e.g.,~\citet{Nemi_09} and~\citet{Juditsky_12a}.

\subsection{Primal function, dual function and duality gap.}\label{sec:primal_dual_functions}
For the SPP~\eqref{eq:main}, % a convex-concave saddle function $S:\bbX\times\bbY\to[-\infty,+\infty]$ (cf.~\eqref{eq:main}),
 we define the associated primal and dual problems as 
\begin{align}
(\bbP):\;\, \min_{x\in\calX}\Big[\barS(x)\defeq \sup_{y\in\calY} S(x,y)\Big] \quad\quad\mbox{and}\quad\quad (\bbD): \;\, \max_{y\in\calY} \left[\underline{S}(x)\defeq \inf_{x\in\calX} S(x,y)\right]. \label{eq:primal_dual_function}
\end{align}
From the definition in~\eqref{eq:def_saddle_point}, we can easily prove the following: Given that a saddle-point $(x^*,y^*)$ exists in~\eqref{eq:main}, both $(\bbP)$ and $(\bbD)$ have nonempty solution sets $\calP^*$ and $\calD^*$, respectively. Furthermore, $x^*\in\calP^*$ and $y^*\in\calD^*$ and $\barS(x^*)=S(x^*,y^*)=\uS(y^*)$. Based on the functions $\barS$ and $\uS$, we define the {\em duality gap}
\begin{align}
G(x,y)&\defeq \barS(x)-\uS(y) = {\sup}_{x'\in\calX,y'\in\calY} \Big[\tilG(x,y;x',y')\defeq S(x,y') - S(x',y)\Big]. \label{eq:duality_gap}
\end{align}

%the projection in~\eqref{eq:Bregman_projection} is {\em tractable} if 
%can be solved in closed form. 
%evaluated in finitely many steps (e.g., )
%on $\calU$ if there exists a DGF $h$ such that for any $u^*\in\bbU^*$, $u'\in\calU^o$ and $\lambda>0$, 
% (Indeed, if $\inf_{u\in\calU}\varphi(u)>-\infty$ and $\calU\cap\dom \varphi_\lambda\ne \emptyset$, then Problem~\eqref{eq:Bregman_projection} always 
%%apart from ``closed form'',  this is always true under mild assumptions; 
%see Lemma~\ref{lem:unique_sln} for details.)

\begin{algorithm}[t!]
\caption{An Optimal Stochastic PDHG-Type Algorithm for Convex $f$ ($\mu\ge 0$)} \label{algo:convex_f}
\begin{algorithmic}
\State {\bf Input}: Interpolation sequence $\{\beta_t\}_{t\in\bbN}$, %\subseteq(0,1]$, %\subseteq[1,+\infty)$, 
dual stepsizes $\{\alpha_t\}_{t\in\bbN}$, primal stepsizes $\{\tau_t\}_{t\in\bbN}$, relaxation sequence $\{\theta_t\}_{t\in\bbN}$, DGFs $h_\calY:\bbY\to\barbbR$ and $h_\calX:\bbX\to\barbbR$%, termination accuracy $$  %number of iterations $K$ %that satisfy~\eqref{eq:cond_seq1},~\eqref{eq:cond_seq2} and~\eqref{eq:cond_seq3}
\State {\bf Initialize}: $x^1\in\calX^o$, $y^1\in\calY^o$, $\barx^1=x^1$, $\bary^1=y^1$, $s^1=\hat{\nabla}_y\Phi(x^1,y^1,\zeta^1_y)$, $t=1$
\State {\bf Repeat} (until some convergence criterion is met)
\vspace{-.2cm}
%\begin{equation}
\begin{align}
%&\mbox{Sample $\bxi^k\!\sim\!\nu$ and define $\bv^k\defeq\nabla_\bx F(\bx,\bxi^k)|_{\bx=\tilbx^k}$}\nn\\ %as in~\eqref{eq:def_v_k}}\nn\\  
&y^{t+1} := \argmin_{y\in\calY}J(y) - \lranglet{s^t}{y-y^t} + \alpha_t^{-1}D_{h_\calY}(y,y^t)%, \mbox{ where } \zeta_x^t\sim p_x^t&
\label{eq:dual_upd}\\
&\tilx^{t+1} :=   (1-\beta_t) \barx^t + \beta_t x^t\label{eq:interp_primal}\\
&x^{t+1} := \argmin_{x\in\calX}g(x) + \lranglet{\hat{\nabla}_x\Phi(x^t,y^{t+1},\zeta^t_x)+\hat{\nabla}f(\tilx^{t+1},\xi^t)}{x-x^t} + \tau_t^{-1}D_{h_\calX}(x,x^t)%\nn\\
%&\hspace{6cm}%\mbox{ where } \zeta_y^t\sim p_y^t,\, \xi^t\sim \nu
\label{eq:primal_upd}\\ 
 %\nn\\[-1ex]%&\hspace{.5cm}\prox_{\tau_k g} (\bx^k- \tau_k(\bA^T\by^{k+1} +\bv^k) )\\
%&z^{t+1} := x^{t+1} +\theta_{t+1} (x^{t+1} - x^t)
&s^{t+1} := (1+\theta_{t+1})\hat{\nabla}_y\Phi(x^{t+1},y^{t+1},\zeta^{t+1}_y) - \theta_{t+1}\hat{\nabla}_y\Phi(x^{t},y^{t},\zeta^{t}_y)\label{eq:extrapolation}\\ 
& \barx^{t+1} :=  (1-\beta_t)\barx^t + \beta_t x^{t+1}\label{eq:ave_primal}\\
& \bary^{t+1} := (1-\beta_t)\bary^t + \beta_t y^{t+1}\label{eq:ave_dual}\\
& t := t+1\label{eq:counter}
\end{align}
\vspace{-1cm}
\State {\bf Output}: $(\barx^t,\bary^t)$ %$\barbx^K$ and  $\barby^K$
\end{algorithmic}%\vspace{-.2cm}
\end{algorithm}

\section{Convex $f$: algorithm and convergence results.}\label{sec:algo_cvx}
We first consider the case where $\mu=0$. We begin with describing our algorithm, followed by the assumptions needed to analyze it, and finally its convergence results. The detailed convergence analysis is deferred to Section~\ref{sec:analysis}. % and the assumptions 

\subsection{Algorithm description.}

The pseudo-code of our algorithm is shown in Algorithm~\ref{algo:convex_f}. 
For input, we require two 
%In steps~\eqref{eq:dual_upd} and~\eqref{eq:primal_upd}, the 
CCP functions $h_\calY:\bbY\to\barbbR$ and $h_\calX:\bbX\to\barbbR$ which are DGFs on $\calY$ and $\calX$ respectively, i.e., they are essentially smooth on their respective domains and %there exist constants $1,1>0$ such that
\begin{align}
D_{h_\calY}(y,y')\ge (1/2)\normt{y-y'}^2, &\;\forall\,y\in\calY, \, \forall\,y'\in\calY^o,\label{eq:sc_calY}\\
D_{h_\calX}(x,x')\ge (1/2)\normt{x-x'}^2, &\;\forall\,x\in\calX, \, \forall\,x'\in\calX^o,\label{eq:sc_calX}
\end{align} 
where $\calY^o\defeq \calY\cap\inter\dom h_\calY$ and $\calX^o\defeq \calX\cap\inter\dom h_\calX$. 
In addition, $h_\calY$ and $h_\calX$ are chosen such that the minimization problems in~\eqref{eq:dual_upd} and~\eqref{eq:primal_upd}  have (unique) and easily computable solutions in $\calY^o\cap\dom J$ and $\calX^o\cap\dom g$ respectively (cf.\ Section~\ref{sec:prelim}). In addition, according to~\eqref{eq:Bregman_diam}, we define the Bregman diameters of $\calX$ and $\calY$ under $h_\calX$ and $h_\calY$ as $\Omega_{h_\calX}$ and $\Omega_{h_\calY}$ respectively,  i.e.,
\begin{equation}
\Omega_{h_\calX}\defeq \sup\nolimits_{x\in\calX,x'\in\calX^o}D_{h_\calX}(x,x'),\quad
\Omega_{h_\calY}\defeq \sup\nolimits_{y\in\calY,y'\in\calY^o}D_{h_\calY}(y,y').  \label{eq:Bregman_diam_XY}
\end{equation}

We next formally describe the stochastic first-order oracles. For each $t\in\bbN$, the oracles return $\hat{\nabla}_y\Phi(x^t,y^{t},\zeta^t_y)$, $\hat{\nabla}_x\Phi(x^t,y^{t+1},\zeta^t_x)$ and $\hat{\nabla}f(\tilx^{t+1},\xi^t)$, which are the unbiased estimators of the gradients $\nabla_y \Phi(x^t,y^t)$, $\nabla_x \Phi(x^t,y^{t+1})$  and $\nabla f(\tilx^t)$, respectively, conditioned on the past information. Here $\zeta^t_y$, $\zeta^t_x$ and $\xi^t$ are the underlying random variables that generate the stochasticity. For analysis purposes, let us define a filtration $\{\calF_{t}\}_{t\in\bbZ_+}$ based on the stochastic process $\{(\zeta^t_y,\zeta^t_x,\xi^t)\}_{t\in\bbN}$.  Specifically, %To analyze the convergence of our algorithms, it is important to state the 
we first define the nested sequence of sets of random variables $\{\Xi_t\}_{t\in\bbZ_+}$
%$\Xi_{-1}\subseteq\Xi_{-1/2}\subseteq\Xi_0\subseteq\cdots$ %of sets of random variables $\{\Xi_{t/2}\}_{t\ge -2}$ 
such that $\Xi_{0}\defeq\{0\}$ and %$\Xi_{-1/2}\defeq \{\zeta_y^0\}$ and 
for any $t\in\bbN$,  $\Xi_t\defeq\{(\zeta_x^i,\zeta_y^i,\xi^i)\}_{i=1}^{t}$. % and $\Xi_{t+1/2}\defeq\Xi_t\cup\{\zeta_y^{t+1}\}$. 
Then, for any $t\in\bbZ_+$, define $\calF_{t}$ to be %a filtration $\{\calF_{\tau/2}\}_{\tau\in\bbZ^+}$ such that for any $\tau\in\bbZ^+$, $\calF_{\tau/2}$ 
the $\sigma$-algebra generated by $\Xi_{t}$, i.e., the minimal $\sigma$-algebra with respect to which $\Xi_{t}$ is measurable. 
%\defeq\sigma(\Xi_{\tau/2})$
%$\calF_0\defeq\emptyset$, $\calF_{1/2}\defeq\sigma(\zeta_y^0)$ and for any $t\in\bbN$,  $\calF_t\defeq\sigma(\{\zeta_x^i,\zeta_y^i,\xi^i\}_{i=0}^{t-1})$ and $\calF_{t+1/2}\defeq\sigma(\{\zeta_x^i,\zeta_y^i,\xi^i\}_{i=0}^{t-1}\cup\{\zeta_y^{t}\})$, where $\sigma(\cdot)$ denotes 
%Equipped with $\{\calF_{t/2}\}_{t\ge -2}$,
In addition, for any $t\in\bbN$, we define the stochastic gradient ``noises'' 
%for any $t\in\bbZ^+$, define  
\begin{align}
\delta^{t}_{y,\Phi} &\defeq \hat{\nabla}_y\Phi(x^t,y^{t},\zeta^t_y) - \nabla_y\Phi(x^t,y^{t}),\label{eq:def_delta_y}\\
 \delta_{x,\Phi}^{t} &\defeq \hat{\nabla}_x\Phi(x^t,y^{t+1},\zeta^t_x) - \nabla_x \Phi(x^t,y^{t+1}),
 \label{eq:def_delta_x}\\
  \delta^{t}_{x,f} & \defeq\hat{\nabla}f(\tilx^{t+1},\xi^t) -  \nabla f(\tilx^{t+1}).\label{eq:def_delta_f} 
\end{align}

We briefly describe the structure of Algorithm~\ref{algo:convex_f}, which falls under the PDHG framework (first introduced by~\citet{Chambolle_11}). % for solving bilinear SPPs. 
Specifically, %in Algorithm~\ref{algo:convex_f}, 
we perform the dual ascent, primal descent and relaxation steps in~\eqref{eq:dual_upd},~\eqref{eq:primal_upd} and~\eqref{eq:extrapolation}, respectively. These three steps constitute the main features of PDHG. Besides these steps, we also perform an interpolation step in~\eqref{eq:interp_primal} to accelerate the convergence of the smooth function $f$. Finally, we perform the averaging steps in~\eqref{eq:ave_primal} and~\eqref{eq:ave_dual} to obtain the ergodic iterate sequence $\{(\barx^t,\bary^t)\}_{t\in\bbN}$.

To measure the progress of Algorithm~\ref{algo:convex_f}, we adopt the {duality gap} (defined in~\eqref{eq:duality_gap}) and analyze the convergence rate of the sequence $\{G(\barx^t,\bary^t)\}_{t\in\bbN}$ (in expectation or with high probability). Using the definition in~\eqref{eq:def_saddle_point}, we easily see that if $\bbE[G(\barx^t,\bary^t)]\to 0$ and $\bbE[(\barx^t,\bary^t)]\to(x^\dagger,y^\dagger)$, then $(x^\dagger,y^\dagger)$ must be a saddle-point of Problem~\eqref{eq:main}. 

\subsection{Assumptions.}\label{sec:assumption}

Before presenting our convergence results, we first describe assumptions on the constraint sets $\calX$ and $\calY$, as well as the stochastic gradient noises $\delta^{t}_{y,\Phi}$, $\delta_{x,\Phi}^{t}$ and $\delta^{t}_{x,f}$. % defined in~\eqref{eq:def_delta_y},~\eqref{eq:def_delta_x} and~\eqref{eq:def_delta_f}. 
%We have the following assumptions about $\delta^{t+1}$
% $\delta_{x,*}^{t+1}$ and $\delta^{t+1}_{y,*}$. % $\varsigma$

\begin{assumption}\label{assump:bounded}\quad
\begin{enumerate}[label={\rm (\Alph*)}]
\item The Bregman diameters $\Omega_{h_\calX}$ and $\Omega_{h_\calY}$ in~\eqref{eq:Bregman_diam_XY} are bounded.\label{assump:bounded_diam} 
\item The set $\calX$ is bounded and the Bregman diameter $\Omega_{h_\calY}$ is bounded. \label{assump:bounded_sets}
\end{enumerate}
%constraint sets $\calX$ and $\calY$ are bounded. 
\end{assumption}

\begin{assumption}\label{assump:noise}
Define the conditional expectation $\bbE_{t}[\cdot]\defeq\bbE[\cdot\,|\,\calF_{t}]$, for any $t\in\bbZ_+$. For any $x\in\calX$ and $y\in\calY$ and any $t\in\bbN$, there exist positive constants $\sigma_{y,\Phi}$, $\sigma_{x,\Phi}$ and $\sigma_{x,f}$ such that %define $\bbE_t[\cdot]\defeq\bbE[\cdot]$%we assume that %the following holds almost surely:
\begin{enumerate}[label={\rm (\Alph*)}]%,wide, labelwidth=!, labelindent=0pt, topsep=2pt]\itemsep.3em
\item (Unbiasedness) \quad\quad $\bbE_{t-1}[\delta_{y,\Phi}^{t}] = 0$, $\bbE_{t-1}[\delta_{x,\Phi}^{t}] = 0$, $\bbE_{t-1}[\delta^{t}_{x,f}] = 0$ a.s., \label{assump:unbiased}
\item (Bounded variance) $\bbE_{t-1}[\normt{\delta_{y,\Phi}^{t}}_*^2] \le \sigma_{y,\Phi}^2$, $\bbE_{t-1}[\normt{\delta_{x,\Phi}^{t}}_*^2] \le \sigma_{x,\Phi}^2$, $\bbE_{t-1}[\normt{\delta^t_{x,f}}_*^2] \le \sigma_{x,f}^2$ a.s.,\label{assump:variance_bound}
\item (Sub-Gaussian distributions)\\
 $\bbE_{t-1}\left[\exp\left(\normt{\delta_{y,\Phi}^{t}}_*^2/\sigma_{y,\Phi}^2\right)\right] \le \exp(1)$, $\bbE_{t-1}\left[\exp\left(\normt{\delta_{x,\Phi}^{t}}_*^2/\sigma_{x,\Phi}^2\right)\right] \le \exp(1)$,\\ 
$\bbE_{t-1}\left[\exp\left(\normt{\delta^t_{x,f}}_*^2/\sigma_{x,f}^2\right)\right] \le \exp(1)$ a.s.. \label{assump:sub-Gaussian}
\end{enumerate} 
\end{assumption}

\subsection{Remarks on assumptions.}
We provide remarks on the Assumptions~\ref{assump:bounded} and~\ref{assump:noise}. % in Section~\ref{sec:assumption}. 

\subsubsection{Remarks on Assumption~\ref{assump:bounded}.}\label{sec:rmk_assump1}
%We make several remarks about the assumptions above. 
First, by~\eqref{eq:sc_calX}, we note that Assumption~\ref{assump:bounded}\ref{assump:bounded_diam} implies Assumption~\ref{assump:bounded}\ref{assump:bounded_sets}. These two assumptions will be used in proving different convergence results of the duality gap in the sequel, and appear to be common in the literature (see e.g.,~\citet{Juditsky_12a,Juditsky_12b}). %for deriving the convergence rate   %and cannot be removed in general (for reasons detailed in Section~\ref{sec:rmk_assump}).
In addition, note that in many scenarios, Assumptions~\ref{assump:bounded}\ref{assump:bounded_diam} and~\ref{assump:bounded_sets} are both {\em equivalent} to the boundedness of $\calX$ and $\calY$. For example, if both $\bbX$ and $\bbY$ are finite dimensional real Hilbert spaces (with inner products $\lranglet{\cdot}{\cdot}_{\bbX}$ and $\lranglet{\cdot}{\cdot}_{\bbY}$ and their induced norms $\norm{\cdot}_{\bbX}$ and $\norm{\cdot}_{\bbY}$ respectively), and we take $h_\calX(x)=(1/2)\norm{x}_{\bbX}^2$ and $h_\calY(y)=(1/2)\norm{y}_{\bbY}^2$, then
\begin{align}
D_{h_\calX}(x,x') = (1/2)\normt{x-x'}_{\bbX}^2,\quad D_{h_\calY}(y,y') = (1/2)\normt{y-y'}_{\bbY}^2. 
\end{align}

Next, we detail the purposes of Assumption~\ref{assump:bounded}. Since either part~\ref{assump:bounded_diam} or part~\ref{assump:bounded_sets} implies the boundedness of the constraint sets $\calX$ and $\calY$, combined with other structural assumptions stated in Section~\ref{sec:intro}, we see that Problem~\eqref{eq:main} has at least one saddle-point, hence being well-posed.  In addition, the compactness of $\calX$ and $\calY$ ensures that the sequence $\{\bbE[(\barx^t,\bary^t)]\}_{t\in\bbN}$ in Algorithm~\ref{algo:convex_f} has at least one limit point in $\calX\times\calY$. Hence if $\bbE[G(\barx^t,\bary^t)]\to 0$ (which will be shown in Theorem~\ref{thm:main_cvx}), then any limit point of $\{\bbE[(\barx^t,\bary^t)]\}_{t\in\bbN}$ is a saddle-point of Problem~\eqref{eq:main}. 
% that there exist a subsequence of  that converges to a point $(x^\ddagger,y^)$. 
Moreover, in the convergence analysis of Algorithm~\ref{algo:convex_f} (cf.\ Section~\ref{sec:analysis}), Assumption~\ref{assump:bounded}\ref{assump:bounded_diam} is also needed, since it provides a uniform  upper bound on the sequence $\{D_{h_\calX}(x,x^t)\}_{t\in\bbN}$ and $\{D_{h_\calY}(y,y^t)\}_{t\in\bbN}$, for any $(x,y)\in\calX\times\calY$. 
\subsubsection{Remarks on Assumption~\ref{assump:noise}.}\label{sec:rmk_assump2}
First, we explain the implications of Assumption~\ref{assump:noise}. 
In Assumption~\ref{assump:noise}, part~\ref{assump:unbiased} states that the stochastic noise process $\{(\delta_{y,\Phi}^{t},\delta_{x,\Phi}^{t}, \delta_{x,f}^{t})\}_{t\in\bbZ_+}$ forms a  (vector-valued) martingale difference sequence (MDS) with respect to the filtration $\{\calF_{t}\}_{t\in\bbZ_+}$. % (in the vector sense). 
Part~\ref{assump:variance_bound} states that the (conditional) second-moment of each of $\{\delta_{y,\Phi}^{t}\}_{t\in\bbN}$, $\{\delta_{x,\Phi}^{t}\}_{t\in\bbN}$ and $\{\delta_{x,f}^{t}\}_{t\in\bbN}$ is uniformly bounded. Under this assumption, we are able to prove certain finite-time convergence results of $\{G(\barx^t,\bary^t)\}_{t\in\bbN}$, both in expectation and with high probability (cf.\ Theorem~\ref{thm:main_cvx}\ref{item:conv_exp}). In particular, for any $\varsigma\in(0,1]$, $\Pr\{G(\barx^T,\bary^T) = O(\varsigma^{-1})\}\ge 1-\varsigma$, where $T$ denotes the total number of iterations of Algorithm~\ref{algo:convex_f}. In this work, however, we are also interested in  obtaining large-deviation-type convergence results of the sequence $\{G(\barx^t,\bary^t)\}_{t\in\bbN}$, i.e.,  $\Pr\{G(\barx^T,\bary^T) = O(\log(1/\varsigma))\}\ge 1-\varsigma$. 
 %, but not enough for showing its convergence with high probability. 
To achieve this, we need to assume that  the  (conditional) distributions of these stochastic noises are ``light-tailed''. Specifically, in part~\ref{assump:sub-Gaussian}, we assume that $\delta_{y,\Phi}^{t}$, $\delta_{x,\Phi}^{t}$ and $\delta_{x,f}^{t}$ are (conditional) {\em sub-Gaussian} random vectors with variance proxies $\sigma_{y,\Phi}^2$, $\sigma_{x,\Phi}^2$ and $\sigma_{x,f}^2$, respectively (see e.g.,~\citet{Vershynin_18}). As we will show in Theorem~\ref{thm:main_cvx}\ref{item:conv_whp},  this assumption indeed enables us to prove the aforementioned large-deviation-type results.  %of the sequence $\{G(\barx^t,\bary^t)\}_{t\in\bbN}$. 

Next, we provide some justifications for Assumption~\ref{assump:noise}\ref{assump:sub-Gaussian}. We will focus on the assumption $\bbE_{t-1}\left[\exp\left(\normt{\delta^t_{x,f}}_*^2/\sigma_{x,f}^2\right)\right] \le \exp(1)$ a.s., as the other  two assumptions can be similarly justified. Note that under the compactness assumption of $\calX$, for any $\xi\in\Xi$, if $\hat{\nabla}f(\cdot,\xi)$ is continuous on $\calX$, then $\barM_f(\xi)\defeq\sup_{x\in\calX} \normt{\hat{\nabla}f(x,\xi)}_*<+\infty$. Similarly, by the $L$-smoothness of $f$ on $\calX$ (cf.\ Section~\ref{sec:intro}), we also have $M_f\defeq \sup_{x\in\calX} \normt{{\nabla}f(x)}_*<+\infty$. Therefore, by the definition of $\delta^t_{x,f}$ , we have for any $t\in\bbN$, 
\begin{align*}
\normt{\delta^t_{x,f}}^2_* \le 2(\normt{\hat{\nabla}f(x,\xi^t)}^2_* + \normt{{\nabla}f(x)}^2_*) \le 2( {\sup}_{\xi\in\Xi}\;\barM^2_f(\xi) + M^2_f)  %<+\infty 
\quad \as
\end{align*}
If ${\sup}_{\xi\in\Xi}\;\barM^2_f(\xi)<+\infty$, which happens if $\barM^2_f(\cdot)$ is continuous on $\Xi$ and $\Xi$ is compact, or $\Xi$ simply has finite size (which implies the finite-sum form of $f$ as in~\eqref{eq:finite-sum}), then we can take $\sigma_{x,f}^2 = 2(\sup_{\xi\in\Xi}\;\barM^2_f(\xi) + M^2_f)$. 

Note that the justification above, although being general, leverages the compactness of $\calX$ and $\Xi$. In fact, for some specific examples of $f$, the compactness assumptions are not needed. Let $\ell:\bbR\to\bbR$ be a loss function with bounded derivative on $\bbR$, i.e., $M_\ell \defeq \sup_{s\in\bbR}\abs{\ell'(s)}<+\infty$. Examples of $\ell$ include the Huber loss function or the quadratically smoothed hinge loss function (see~\citet{Zhang_04}). Let $f(x) = \bbE[\ell(\lranglet{\xi}{x})]$, where $\xi\in\bbX^*$ is a sub-Gaussian random vector, i.e., $\bbE[\xi] = 0$ and $\bbE[\exp(\norm{\xi}^2_*/\sigma_\xi^2)]\le \exp(1)$, for some $\sigma_\xi>0$. (For convenience, we write $\xi\sim\mathsf{subG}(\sigma_\xi^2)$.)
Note that the sub-Gaussian assumption on $\xi$ is standard in statistical learning theory (see e.g.,~\citet{Vershynin_18}). By the Leibniz integral rule, we have $\nabla f(x) = \bbE[\ell'(\lranglet{\xi}{x})\xi]$. At each $t\in\bbN$, let $\hat{\nabla} f(x,\xi^t) = \ell'(\lranglet{\xi^t}{x})\xi^t$, where $\xi^t\sim \mathsf{subG}(\sigma_\xi^2)$. Then %we have
\begin{align*}
\bbE[\exp\big(\normt{\delta^t_{x,f}}^2_*/(4M_\ell^2\sigma_\xi^2)\big)] &\le \bbE[\exp\big(2\big(\normt{\ell'(\lranglet{\xi^t}{x})\xi^t}^2_* + \normt{\bbE[\ell'(\lranglet{\xi}{x})\xi]}^2_*\big)/(4M_\ell^2\sigma_\xi^2)\big)]\\
&\lea  \bbE[\exp\big(\normt{\xi^t}^2_*/(2\sigma_\xi^2) + \bbE[\normt{\xi}^2_*]/(2\sigma_\xi^2)\big)]\\
&\leb  \bbE[\exp\big(\normt{\xi^t}^2_*/\sigma_\xi^2\big)^{1/2}] \exp(1/2)\\
&\lec  \bbE[\exp\big(\normt{\xi^t}^2_*/\sigma_\xi^2\big)]^{1/2} \exp(1/2) \led \exp(1),
\end{align*}
where in (a) we use $\abs{\ell'(s)}\le M_\ell$, for any $s\in\bbR$, in (b) we use $\bbE[\normt{\xi}^2_*]\le \sigma_\xi^2$, in (c) we use Jensen's inequality, and in (d) we use $\bbE[\exp(\norm{\xi^t}^2_*/\sigma_\xi^2)]\le \exp(1)$. Therefore, %in this case, 
we can take $\sigma_{x,f}^2 =4M_\ell^2\sigma_\xi^2$. 
% As we will see in Section~\ref{sec:analysis}, %this sub-Gaussianity 
%such an assumption allows us to invoke concentration inequalities (e.g., Asuma-Hoeffing) to

%and~\ref{assump:sub-Gaussian} bound the variance of $\delta$, $\delta_{x,*}$ and $\delta_{y,*}$, but~\ref{assump:variance_LDbound} is stronger than~\ref{assump:variance_bound}. These two  assumptions will be used in proving different convergence results in the sequel. In addition, Assumption~\ref{assump:sub-Gaussian} is slightly stronger than the usual one used in the literature (e.g.,~\citet{Nemi_09}), where the bounds are only required to hold for $\theta=1$. As we will see in Theorem~\ref{thm:main_cvx}\ref{res:cvx_LD}, this stronger assumption leads to a stronger large-deviation result. 
%\end{remark}

\subsection{Convergence results.}\label{sec:conv_res}
%We analyze the convergence of the duality gap in two aspects, i.e., in expectation and with high probability. 
Before presenting our main results, we first introduce an important proposition that will be used in our stochastic restart scheme (cf.~Section~\ref{sec:f_sc}). Its proof is deferred to Section~\ref{sec:analysis}. 

\begin{proposition}\label{prop:recursion_cvx}
Let Assumptions~\ref{assump:bounded}\ref{assump:bounded_diam} and~\ref{assump:noise}\ref{assump:unbiased} hold. %and $\Omeg1^2,\Omeg1^2<+\infty$. 
Define $\theta_0=0$, $\beta_0=2$, $\alpha_0=\tau_0=1$. (Note that these parameters are merely for analysis purpose, and do not appear in Algorithm~\ref{algo:convex_f}, which starts with $t=1$.) If there exists a nonnegative sequence $\{\gamma_t\}_{t\in\bbZ_+}$ that satisfies $\gamma_0= 0$ and  
\begin{align}
&0\le \theta_t\le 1,\quad \theta_{t-1}\le\theta_t,\quad\alpha_t\theta_t\le \alpha_{t-1}, \quad \gamma_t\theta_t= \gamma_{t-1},\label{eq:cond_set1}\\
&\gamma_{t-1}\beta_{t-1}^{-1} = \gamma_{t}(\beta_{t}^{-1}-1), \quad \gamma_{t-1}/\tau_{t-1} \le \gamma_t/\tau_t,\quad \alpha_t\le (2L_{yy})^{-1},\label{eq:cond_set2}\\
&L\beta_t+L_{xx}-({2\tau_t})^{-1}+{4\alpha_tL_{yx}^2}\le 0, \quad (1+\theta_t)L_{yy}-(8\alpha_t)^{-1}\le 0,\label{eq:cond_set3}
\end{align}
for any $t\in\bbN$, then for any $T\ge 3$, %the following holds almost surely: %
we have the following: %with probability one (with probability1) that 
\begin{enumerate}[label={\rm (\Alph*)}]
\item \label{item:conv_exp} If Assumption \ref{assump:noise}\ref{assump:variance_bound} also holds, then 
\begin{align}
\bbE[G(\barx^T,\bary^T)] &\le ({\beta_T^{-1}-1})^{-1}(B_1(T) + B_2(T)),\label{eq:recur_res_exp}\\
\mbox{where}\quad \quad B_1(T)&\defeq \frac{2\theta_T}{\tau_{T-1}}\Omega_{h_\calX}+\frac{4\theta_T}{\alpha_{T-1}}\Omega_{h_\calY},\\
 B_2(T)&\defeq \frac{4(\sigma_{x,\Phi}^2+\sigma_{x,f}^2)}{\gamma_T }\sum_{t=1}^{T-1}\gamma_t\tau_t  + \frac{22\sigma_{y,\Phi}^2}{\gamma_T }\sum_{t=1}^{T-1}\gamma_t\alpha_t.
\end{align}
In addition, for any $\varsigma\in(0,1]$, with probability at least $1-\varsigma$, we have
\begin{align}
G(\barx^T,\bary^T) &\le \varsigma^{-1}({\beta_T^{-1}-1})^{-1}(B_1(T) + B_2(T)). \label{eq:conv_markov}
\end{align}
\item \label{item:conv_whp} Define the stochastic sequence $\{\hatx^t\}_{t\in\bbN}$ such that $\hatx^1=x^1$ and 
\begin{equation}
\hatx^{t+1} \defeq {\argmin}_{x\in\calX} -\lranglet{\delta_{x,\Phi}^{t}+\delta_{x,f}^{t}}{x} + \tau_t^{-1} D_{h_\calX}(x,\hatx^t), \quad \forall\,t\in\bbN. \label{eq:def_hat_x}
\end{equation}
If Assumption \ref{assump:noise}\ref{assump:sub-Gaussian} also holds, then for any $\varsigma\in(0,1/6)$, with probability at least $1-6\varsigma$, we have 
\begin{align}
\tilG(\barx^T,\bary^T;x,y) &\le \frac{1}{\gamma_T(\beta_T^{-1}-1)}\Bigg\{\sum_{t=1}^{T-1}\left(\frac{\gamma_t}{\tau_t}-\frac{\gamma_{t-1}}{\tau_{t-1}}\right)(D_{h_\calX}(x,x^t)+D_{h_\calX}(x,\hatx^t))+\frac{4\gamma_{T-1}}{\alpha_{T-1}}\Omega_{h_\calY}\nn\\
 &+ \gamma_T(1+\log(1/\varsigma))B_2(T)+\gamma_T(\beta_T^{-1}-1)B_3(T)\Big\},\quad\forall\,(x,y)\in\calX\times\calY,\label{eq:conv_whp_pseudo}\\
\mbox{where}\quad\quad B_3(T) &\defeq  \frac{2\sqrt{\log(1/\varsigma)}}{\gamma_T(\beta_T^{-1}-1)}\left(2{\sigma_{y,\Phi}{D_\calY}}+{(\sigma_{x,\Phi}+\sigma_{x,f}){D_\calX}}\right)\left(\textstyle{\sum}_{t=1}^{T-1}\gamma_t^2\right)^{1/2}. 
\end{align}
Furthermore, 
\begin{align}
G(\barx^T,\bary^T) \le &({\beta_T^{-1}-1})^{-1}(B_1(T) + (1+\log(1/\varsigma))B_2(T))+B_3(T). \label{eq:recur_res_hp}
\end{align}
\end{enumerate}
\end{proposition}

%\proof{Proof.}
%See Section~\ref{sec:analysis}.\Halmos
%\endproof

Based on Proposition~\ref{prop:recursion_cvx}, we obtain the main convergence results for Algorithm~\ref{algo:convex_f}. 

\begin{theorem}\label{thm:main_cvx}
Let Assumptions~\ref{assump:bounded}\ref{assump:bounded_diam} and~\ref{assump:noise}\ref{assump:unbiased} hold. 
In Algorithm~\ref{algo:convex_f}, for any $t\in\bbN$, choose %$\theta_0=0$, $\beta_0=2$, $\alpha_0 =\tau_0=1$ and 
 \begin{align}
&\theta_t = \frac{t-1}{t},\quad \beta_t = \frac{2}{t+1},\quad \alpha_t = \frac{1}{16\left(L_{yx}+L_{yy}+\rho\sigma_{y,\Phi}\sqrt{t}\right)},\label{eq:param_choice_1}\\
&\tau_t = \frac{t}{2\left(2L+(L_{xx}+L_{yx})t+\rho'(\sigma_{x,\Phi}+\sigma_{x,f})t^{3/2}\right)},\label{eq:param_choice_2}
\end{align}
where $\rho,\rho'>0$ are constants independent of $(L,L_{xx},L_{yx},L_{yy},\sigma_{x,f},\sigma_{x,\Phi},\sigma_{y,\Phi},t)$.
\begin{enumerate}[label={\rm (\Alph*)}]
\item \label{res:cvx_exp} If Assumption~\ref{assump:noise}\ref{assump:variance_bound} also holds, then for any $T\ge 3$, we have
\begin{align}
\bbE[G(\barx^T,\bary^T)]\le & B_{\rm E}(T)\defeq\frac{16L}{T(T-1)}\Omega_{h_\calX} + \frac{8(L_{xx}+L_{yx})}{ T}\Omega_{h_\calX} + \frac{128(L_{yx}+L_{yy})}{ T}\Omega_{h_\calY}  \nn\\
 &+ \frac{8\sigma_{y,\Phi}}{\sqrt{T}}\left(\frac{1}{\rho}+{16\rho}\Omega_{h_\calY}\right) + \frac{8(\sigma_{x,f}+\sigma_{x,\Phi})}{\sqrt{T}}\left(\frac{1}{\rho'}+{\rho'}\Omega_{h_\calX}\right).
 %=O\left(\frac{L}{T^2}+\frac{B}{T}+\frac{\sigma+\sigma_{x,*}+\sigma_{y,*}}{\sqrt{T}}\right).
 \label{eq:cvx_bound_rho}
\end{align}
In addition, for any $\varsigma\in(0,1]$, with probability at least $1-\varsigma$, we have
\begin{align}
G(\barx^T,\bary^T) &\le \varsigma^{-1}B_{\rm E}(T). \label{eq:cvx_bound_markov_rho}
\end{align}
\item \label{res:cvx_LD} %Let $\varsigma\in(0,1/6]$ and 
If Assumption \ref{assump:noise}\ref{assump:sub-Gaussian} also holds, then for any $T\ge 3$ and $\varsigma\in(0,1/6]$, with probability at least $1-6\varsigma$, we have 
\begin{align}
G(\barx^T,\bary^T)\le B_{\rm E}(T) &+ \frac{8\sigma_{y,\Phi}}{\sqrt{T}}\left(\frac{\log(1/\varsigma)}{\rho}+\sqrt{{\log(1/\varsigma)}}D_\calY\right)\nn\\
 &+ \frac{8(\sigma_{x,\Phi}+\sigma_{x,f})}{\sqrt{T}}\left(\frac{\log(1/\varsigma)}{\rho'}+\sqrt{{\log(1/\varsigma)}}D_\calX\right). \label{eq:cvx_bound_rhoLD}
\end{align}
%holds with probability $1-4\varsigma$. 
\end{enumerate} 
\end{theorem}
%Recall that $M=L_{xx} + 2L_{yx} + L_{yy}$.%and define $\sigma\defeq \sigma_{x,\Phi}+\sigma_{x,f} + \sigma_{y,\Phi}$.   

\proof{Proof.}
We first verify the choices of the input sequences $\{\beta_t\}_{t\in\bbN}$, $\{\alpha_t\}_{t\in\bbN}$, $\{\tau_t\}_{t\in\bbN}$ and $\{\theta_t\}_{t\in\bbN}$ in Theorem~\ref{thm:main_cvx} indeed satisfy the conditions required in Proposition~\ref{prop:recursion_cvx}. Indeed, based on these choices, we can choose $\gamma_t=t$, for any $t\in\bbN$. We only show the steps to verify the conditions in~\eqref{eq:cond_set3}. First, since $\tau_t^{-1}\ge (4L+2(L_{yx}+L_{xx})t)/t$ and $\alpha_t\le 1/(16L_{yx})$, we have
\begin{align}
L\beta_t+L_{xx}-({1/2})\tau_t^{-1}+{4\alpha_tL_{yx}^2}\le 2L/(t+1)+L_{xx}-(2L+(L_{yx}+L_{xx})t)/t+L_{yx}/4\le 0.\nn
\end{align}
% since the verification for the other conditions are trivial. 
Also, since $\alpha_t^{-1}\ge 16L_{yy}$ and $\theta_t\le 1$, we have 
\begin{equation}
(1+\theta_t)L_{yy}-(1/8)\alpha_t^{-1}\le (1+\theta_t)L_{yy}-2L_{yy} \le 0. 
\end{equation}

%Now, we have
Next, we bound the summation terms appearing in~\eqref{eq:recur_res_exp} and~\eqref{eq:recur_res_hp}. Specifically, 
by noting that $\alpha_t\le \left(16\rho\sigma_{y,\Phi}\sqrt{t}\right)^{-1}$ and  $\tau_t\le (2\rho'\left(\sigma_{x,f}+\sigma_{x,\Phi})\sqrt{t}\right)^{-1}$, we have %for any $t\in\bbN$, 
\begin{align}
 \sum_{t=1}^{T-1} \gamma_t\alpha_t &\le \frac{1}{16\rho\sigma_{y,\Phi}}\sum_{t=1}^{T-1} \sqrt{t}\le \frac{1}{16\rho\sigma_{y,\Phi}}\int_{t=0}^T\sqrt{t} \rmd t \le \frac{T^{3/2}}{16\rho\sigma_{y,\Phi}},
\label{eq:cvx_final_bound1}\\
\sum_{t=1}^{T-1} \gamma_t\tau_t&\le \frac{1}{2\rho'(\sigma_{x,f}+\sigma_{x,\Phi})}\sum_{t=1}^{T-1}\sqrt{t}\le  \frac{T^{3/2}}{2\rho'(\sigma_{x,f}+\sigma_{x,\Phi})}, \label{eq:cvx_final_bound2}\\
\sum_{t=1}^{T-1} \gamma_t^2&=\sum_{t=1}^{T-1} t^2 \le \int_{t=0}^T t^2 \rmd t = \frac{1}{3}T^{3}. \label{eq:cvx_final_bound3}
\end{align}
%where in both (a) and (b) we use $\sum_{t=1}^{T-1}\sqrt{t}\le \int_{t=0}^T\sqrt{t} \rmd t=(2/3)T^{3/2}$. 
%and $\{t\tau_t^{-1}\}_{t\in\bbZ^+}$ are increasing sequences, so that both $$
We then substitute~\eqref{eq:param_choice_1},~\eqref{eq:param_choice_2},~\eqref{eq:cvx_final_bound1},~\eqref{eq:cvx_final_bound2} and~\eqref{eq:cvx_final_bound3} into~\eqref{eq:recur_res_exp},~\eqref{eq:conv_markov} and~\eqref{eq:recur_res_hp} to obtain~\eqref{eq:cvx_bound_rho}, \eqref{eq:cvx_bound_markov_rho} and~\eqref{eq:cvx_bound_rhoLD}. 
\Halmos
\begin{remark}\label{rmk:cvx}
Note that the parameter choices in Theorem~\ref{thm:main_cvx} do not involve the Bregman diameters $\Omega_{h_\calX}$ and $\Omega_{h_\calY}$. However, if they are known (or can be estimated), we can choose $\rho={1}/(4\sqrt{\Omega_{h_\calY}})$ and $\rho'=1/\sqrt{\Omega_{h_\calX}}$ to ``optimize'' the bound in~\eqref{eq:cvx_bound_rho},~\eqref{eq:cvx_bound_markov_rho} and~\eqref{eq:cvx_bound_rhoLD}. In Section~\ref{sec:estimation_diameters}, we will provide some methods to estimate $\Omega_{h_\calX}$ and $\Omega_{h_\calY}$. 
\end{remark}

\begin{remark}\label{rmk:complexity_cvx_f}
%Recall that the output of Algorithm~\ref{algo:convex_f} is denoted by $(\barx^T,\bary^T)$. % denotes the 
Theorem~\ref{thm:main_cvx} indicates that %the expected convergence rate of the duality gap is $O(L/T^2+(L_{xx} + L_{yx} + L_{yy})/T+(\sigma_{x,f} + \sigma_{x,\Phi} + \sigma_{y,\Phi})/\sqrt{T})$. In other words, 
to find a (stochastic) primal-dual pair $(x,y)\in\calX\times\calY$ %with an $\epsilon$-expected duality gap (i.e.,
such that  $\bbE[G(x,y)]\le\epsilon$, the oracle complexity of Algorithm~\ref{algo:convex_f} is
\begin{equation}
O\left(\sqrt{\frac{L}{\epsilon}} + \frac{L_{xx} + L_{yx} + L_{yy}}{\epsilon} + \frac{(\sigma_{x,f}+\sigma_{x,\Phi})^2+\sigma_{y,\Phi}^2}{\epsilon^2}\right). \label{eq:complexity_cvx_f}
\end{equation}
In addition, to find $(x,y)\in\calX\times\calY$ %$\epsilon$-duality gap with probability at least $1-\varsigma$ (i.e., 
such that $\Pr\{G(x,y)\le\epsilon\}\ge 1-\varsigma$, the oracle complexity of Algorithm~\ref{algo:convex_f} is
\begin{equation}
O\left(\frac{1}{\varsigma}\left\{\sqrt{\frac{L}{\epsilon}} + \frac{L_{xx} + L_{yx} + L_{yy}}{\epsilon} + \frac{(\sigma_{x,f}+\sigma_{x,\Phi})^2+\sigma_{y,\Phi}^2}{\epsilon^2}\right\}\right). \label{eq:complexity_cvx_f_markov}
\end{equation}
under Assumption~\ref{assump:noise}\ref{assump:variance_bound} and 
\begin{equation}
O\left(\sqrt{\frac{L}{\epsilon}} + \frac{L_{xx} + L_{yx} + L_{yy}}{\epsilon} + \frac{(\sigma_{x,f}+\sigma_{x,\Phi})^2+\sigma_{y,\Phi}^2}{\epsilon^2}\log\bigg(\frac{1}{\varsigma}\bigg)\right) \label{eq:complexity_cvx_f_hp}
\end{equation}
under Assumption~\ref{assump:noise}\ref{assump:sub-Gaussian}.  Comparing~\eqref{eq:complexity_cvx_f_markov} and~\eqref{eq:complexity_cvx_f_hp},  we observe that
the sub-Gaussian assumption on the gradient noises (cf.\ Assumption~\ref{assump:noise}\ref{assump:sub-Gaussian}) enables us to obtain a significantly better complexity result, in terms of its dependence on $\varsigma$. %Such an observation also applies to the oracle complexity results in Section~\ref{sec:f_sc}, where  we assume that $f$ is $\mu$-strongly convex on $\calX$, with $\mu>0$. 
\end{remark}
%the strong convexity parameter of $f$ i.e., $\mu>0$. 

\begin{remark}
Note that as it becomes customary in the literature (see e.g.,~\citet{Lan_12,Ghad_13b}), the optimality of the dependence of $B_\rmE(T)$ (i.e., the bound of $\bbE[G(\barx^T,\bary^T)]$ in~\eqref{eq:cvx_bound_rho}) on the  diameters $\Omega_{h_\calX}$, $\Omega_{h_\calY}$, $D_\calX$ and $D_\calY$ is not discussed. As a result, we do not include these diameters in~\eqref{eq:complexity_cvx_f} and~\eqref{eq:complexity_cvx_f_hp}. With that said, the dependence of $B_\rmE(T)$ on these diameters is indeed the same as those in the  existing methods (cf.\ Table~\ref{table:cvx}). Similar comments also apply to all the oracle complexity results in Section~\ref{sec:f_sc}. 
\end{remark}

\section{Strongly Convex $f$: restart scheme and complexity analysis.}\label{sec:f_sc}
%We consider the special case where $f$ is $\mu$-strongly convex (s.c.) on $\calX$ (where $\mu>0$), i.e., for any $x,x'\in\calX$,
%\begin{equation}
%f(x)\ge f(x') + \lrangle{\nabla f(x')}{x-x'} + (\mu/2)\normt{x-x'}^2. 
%\end{equation}
We consider the case where $\mu>0$. 
%In addition, we assume that $\Phi(x,\cdot)$ is a linear function for any $x\in\calX$, i.e., there exists a function $\phi:\bbX\to\bbY^*$ such that $\Phi(x,\cdot) = \lrangle{\phi(x)}{\cdot}$. We also assume that $\sigma_{y,\Phi}=0$. 
%As a result, $L_{yy}=0$,  and the stochastic first-order oracle that returns $\hat{\nabla}_y\Phi(x,y)$ degenerates to the deterministic zeroth-order oracle that returns $\phi(x)$.   
%Under these assumptions, 
We aim to develop restart schemes based on Algorithm~\ref{algo:convex_f} that significantly improve the oracle complexities in~\eqref{eq:complexity_cvx_f} and~\eqref{eq:complexity_cvx_f_hp}. %by leveraging the restart technique. 

\subsection{Algorithm~\ref{algo:convex_f} with rescaled DGF.}\label{sec:algo_modified}
We first introduce a variant of Algorithm~\ref{algo:convex_f} (developed for $\mu\ge 0$) that will be used as a subroutine in our restart scheme. Fix any $x_\rmc\in\calX^o$ and $R>0$, and define a new DGF on $\barcalX(x_\rmc,R)\defeq R\calX+x_\rmc$: % and $\calY$, respectively, i.e., %for any $x\in R\calX+x_\rmc$ and $y\in\calY$, 
\begin{equation}
\tilh_{\barcalX(x_\rmc,R)}(x)\defeq R^2 h_\calX\left(\frac{x-x_\rmc}{R}\right). \label{eq:new_DGFs}
\end{equation}
Based on the  definition in~\eqref{eq:new_DGFs}, we easily see that the corresponding Bregman distance %$D_{\tilh_{\barcalX(x_\rmc,R)}}(\cdot,\cdot)$ and $D_{\tilh_\calY}(\cdot,\cdot)$  %(induced by $\tilh_\calY$ and $\tilh_{\barcalX(x_\rmc,R)}$, respectively) 
%have the forms 
\begin{align}
D_{\tilh_{\barcalX(x_\rmc,R)}}(x,x') &= R^2\left\{h_\calX\left(\frac{x-x_\rmc}{R}\right) - h_\calX\left(\frac{x'-x_\rmc}{R}\right)-\lrangle{\nabla h_\calX\left(\frac{x'-x_\rmc}{R}\right)}{\frac{x-x'}{R}}\right\}, %\quad\forall\,x\in\calX,\forall\,x'\in\calX^o.
 \label{eq:Bregman_calX}
%D_{\tilh_\calY}(y,y') &= \eta_y D_{h_\calY}(y,y'). \label{eq}
\end{align}
for any $x\in\calX$ and $x'\in\calX^o$. Since from~\eqref{eq:sc_calX},  we know that $D_{h_\calX}(x,x')\ge (1/2)\normt{x-x'}^2$, by  using the expression of $D_{\tilh_{\barcalX(x_\rmc,R)}}(\cdot,\cdot)$ in~\eqref{eq:Bregman_calX}, we also have $D_{\tilh_{\barcalX(x_\rmc,R)}}(x,x')\ge (1/2)\normt{x-x'}^2$. 

Define %the $\normt{\cdot}$-ball 
$\calB(x_\rmc,R)\defeq \{x\in\bbX:\norm{x-x_\rmc}\le R\}$. % and $\barcalX'(x_\rmc,R)\defeq R^{-1}(\calX-x_\rmc)$. % and $y_{\min}\defeq \argmin\nolimits_{y\in\calY}h_\calY(y)$.  
If $\calB(0,1)\subseteq \dom h_\calX$,  
we then have 
\begin{align}
{\sup}_{x\in\calX\bigcap\calB(x_\rmc,R)}D_{\tilh_{\barcalX(x_\rmc,R)}}(x,x_\rmc)&\le \;\;{\sup}_{x\in\calB(x_\rmc,R)}D_{\tilh_{\barcalX(x_\rmc,R)}}(x,x_\rmc) \le R^2 \Omega'_{h_\calX},\nn  \\
&\quad\mbox{where}\;\;\Omega'_{h_\calX}\defeq {\sup}_{z\in \calB(0,1)}\;D_{h_\calX}(z,0) < +\infty. \label{eq:diameter_Rbound}
%\\
%{\sup}_{y\in\calY}D_{\tilh_\calY}(y,y_{\min}) &= \eta_y\big\{{\sup}_{y\in\calY}h_\calY(y) -{\min}_{y\in\calY}h_\calY(y)\big\}\nn\\
%&\defeq \Omega'_{\tilh_\calY}<+\infty,\label{eq:Omega_y'} 
\end{align}
In words, the quantity $\Omega'_{h_\calX}$ is the Bregman diameter of the unit ball $\calB(0,1)$ under $h_\calX$, which can be interpreted as the {\em normalized} Bregman diameter of $\calB(x_\rmc,R)$ under $\tilh_{\barcalX(x_\rmc,R)}$. 
Note that the condition $\calB(0,1)\subseteq \dom h_\calX$ is satisfied %(for any $R>0$ and $x_\rmc\in\calX^o$) 
when $\bbX$ is a Hilbert space and $h_\calX=(1/2)\norm{\cdot}^2$, in which case $\dom h_\calX=\bbX$. For some other examples, we refer readers to~\citet[Section~4]{Nest_05} and~\citet[Section~5]{Nemi_05}. 

%for certain examples of the triple $(\calX,g,h_\calX)$ that defines a BPP with closed-form solution (see ), the condition $\barcalX'(x_\rmc,R)\cap \calB(0,1)\subseteq \dom h_\calX$ is indeed satisfied, for any $R>0$ and $x_\rmc\in\calX^o$. %set $\calX$ indeed has non-empty intersection with $\calB(0,1)$. 
%where~\eqref{eq:Omega_y'} is guaranteed by Assumption~\ref{assump:bounded}\ref{assump:bounded_sets}. 

%Equipped with the definitions above, we are ready to present our algorithm. 

\renewcommand\thealgorithm{1R}
\begin{algorithm}[t!]
\caption{Algorithm~\ref{algo:convex_f} with Rescaled DGF} \label{algo:convex_f_scaled}
\begin{algorithmic}
\State {\bf Input}:  $x^0\in\calX^o$, positive radius $R>0$, primal constraint set $\calX'\subseteq\calX$, % (such that $x^* \in\calX'\subseteq\calX$), %error probability $\nu>0$, 
number of iterations~$T$, interpolation sequence $\{\beta_t\}_{t\in\bbN}$, %\subseteq(0,1]$, %\subseteq[1,+\infty)$, 
dual stepsizes $\{\alpha_t\}_{t\in\bbN}$, primal stepsizes $\{\tau_t\}_{t\in\bbN}$, relaxation sequence $\{\theta_t\}_{t\in\bbN}$, DGFs $h_\calY:\bbY\to\barbbR$ and $h_\calX:\bbX\to\barbbR$
%, termination accuracy $$  %number of iterations $K$ %that satisfy~\eqref{eq:cond_seq1},~\eqref{eq:cond_seq2} and~\eqref{eq:cond_seq3}
\State {\bf Initialize}: $x^1=x^0$, $y^1\in\calY^o$, $\barx^1=x^1$, $\bary^1=y^1$, $s^1=\hat{\nabla}_y\Phi(x^1,y^1,\zeta^1_y)$ %, $t=1$
\State {\bf Define}: $\barcalX({x^1,R})$ and  $\tilh_{\barcalX({x^1,R})}$ using $h_\calX$, $x^1$ and $R$ as in~\eqref{eq:new_DGFs}, with $x_\rmc$ replaced by $x^1$
\State {\bf For} $t = 1,\ldots,T-1$
\vspace{-.2cm}
%\begin{equation}
\begin{align}
&\mbox{Run steps~\eqref{eq:dual_upd} to~\eqref{eq:ave_dual} in Algorithm~\ref{algo:convex_f}, except changing~\eqref{eq:primal_upd} to}\nn\\
&\quad \quad x^{t+1} := \argmin_{x\in \calX'}g(x) + \lranglet{\hat{\nabla}_x\Phi(x^t,y^{t+1},\zeta^t_x)+\hat{\nabla}f(\tilx^{t+1},\xi^t)}{x-x^t} + \tau_t^{-1}D_{\tilh_{\barcalX({x^1,R})}}(x,x^t).\label{eq:primal_upd_scaled}
\end{align}
\vspace{-.8cm}
\State {\bf Output}: $(\barx^T,\bary^T)$ %$\barbx^K$ and  $\barby^K$
\end{algorithmic}%\vspace{-.2cm}
\end{algorithm}

Our modified algorithm is shown in Algorithm~\ref{algo:convex_f_scaled}. %We make two remarks about it. 
%First, in the input, $x^*$ denotes primal part of any saddle-point $(x^*,y^*)$ of~\eqref{eq:main} (cf.~\eqref{eq:def_saddle_point}). %The way to choose the positive radius $R$ will be explained in our restart 
%In addition, it is the unique minimizer of the primal function $\barS$ on $\calX$ (cf.~\eqref{eq:primal_dual_function}). 
%As will be shown in Section~\ref{sec:det_restart}, when Algorithm~\ref{algo:convex_f_scaled} is used as a subroutine in our restart scheme, the positive radius $R$ will always satisfy $R\ge 2\normt{x^0-x^*}$. % will be guaranteed. 
Compared to Algorithm~\ref{algo:convex_f}, it has two differences. First, in Algorithm~\ref{algo:convex_f_scaled}, we fix the total number of steps $T-1$ before the algorithm starts. This is because Algorithm~\ref{algo:convex_f_scaled} will be used as the subroutine in our restart scheme, where at each stage, 
we terminate it once the total number of iterations is reached. Second, we replace step~\eqref{eq:primal_upd} in Algorithm~\ref{algo:convex_f} with step~\eqref{eq:primal_upd_scaled}.  Compared with step~\eqref{eq:primal_upd}, in step~\eqref{eq:primal_upd_scaled}, we change the constraint set from $\calX$ to $\calX'$ and the Bregman distance from $D_{h_\calX}(\cdot,\cdot)$ to $D_{\tilh_{\barcalX({x^1,R})}}(\cdot,\cdot)$. The motivations for these changes will become apparent in the subsequent analysis. 

In addition, we remark that  %in step~\eqref{eq:primal_upd_scaled}, the set $\calX'$ 
if step~\eqref{eq:primal_upd} has an easily computable solution, then it is reasonable to assume that step~\eqref{eq:primal_upd_scaled} does as well. Three cases that guarantee this are:%where~\eqref{eq:primal_upd} and~\eqref{eq:primal_upd_scaled} have the same form 
\begin{itemize}
\item $\bbX$ is a normed space, $g\equiv 0$ and $\calX'=\calX=\bbX$,
\item $\bbX$ is a Hilbert space, $\calX'=\calX$  and $h_\calX=(1/2)\norm{\cdot}^2$,
\item $\bbX$ is a Hilbert space, $g\equiv 0$, $\calX'$ has an easily computable orthogonal projection  and $h_\calX=(1/2)\norm{\cdot}^2$.
\end{itemize}

%we can take . % (cf.~Section~\ref{sec:rmk_assump}). %If in addition, $g$ is positively homogeneous (i.e., for any $\lambda>0$ and $x\in\bbX$, there exist $r\in\bbR$ independnet of $(\lambda,x)$, such that $g(\lambda x)=\lambda^r g(x)$), we can take $\calX'=\lambda\calX$, for any $\lambda>0$.  %In this case, we can take $\calX'=\calX$. If in addition, $g\equiv 0$ can be any set such that~\eqref
%If in addition, $g\equiv 0$, we can take $h_\calX=(1/2)\norm{\cdot}_\bbX^2$ and $\calX'$ to be any set that admits a closed-form orthogonal projection.  %As a result,~\eqref{eq:primal_upd_scaled} simply becomes a Bregman projection. 
%Finally, if $\bbX$ is a Hilbert space and $g\equiv 0$, 
%We make a few remarks about Algorithm~\ref{algo:convex_f_scaled}. First, 

%We now present the convergence results of Algorithm~\ref{algo:convex_f_scaled}. 

 \subsection{Deterministic restart scheme for strongly convex $f$.}\label{sec:det_restart}
For ease of exposition, we first develop our restart scheme in the case where $\sigma=0$, i.e., we can obtain the {deterministic} gradients of $f$, $\Phi(\cdot,y)$ and $\Phi(x,\cdot)$.  (The restart scheme for the stochastic case, where $\sigma_{x,f}, \sigma_{x,\Phi}, \sigma_{y,\Phi}>0$, will be developed in Section~\ref{sec:stoc_restart}.) %To do so, 
We start with analyzing the convergence properties of Algorithm~\ref{algo:convex_f_scaled}. % in the deterministic setting. 

\begin{proposition}\label{prop:conv_scaled}
Assume that $\sigma=0$, $ \calB(0,1)\subseteq \dom h_\calX$ and $\Omega_{h_\calY}<+\infty$. %In addition, assume that  
In Algorithm~\ref{algo:convex_f_scaled}, fix 
\begin{align}
T\ge \Big\lceil\!\max\Big\{3,&\;64\textstyle{\sqrt{{(L/\mu)\Omega'_{h_\calX}}}},\;{1024(L_{xx}/\mu)\Omega'_{h_\calX}},\;4096L_{yx}(\mu R)^{-1}{\sqrt{\Omega'_{h_\calX}\Omega_{h_\calY}}},\;8192L_{yy}(\mu R^2)^{-1}\Omega_{h_\calY}\Big\}\Big\rceil. \label{eq:choice_T}
\end{align}
If $x^*\in\calX'$, $R\ge 2\normt{x^0-x^*}$, and we choose $\{\beta_t\}_{t=1}^{T}$ and $\{\theta_t\}_{t=1}^T$ as in~\eqref{eq:param_choice_1}, and $\alpha_t=\alpha$ and $\tau_t=t\tau$ for any $t\in[T]$, where
\begin{align}
\alpha = 1/\left(16(\eta^{-1}L_{yx} + L_{yy})\right),\quad \tau = 1/\left(4L+2(L_{xx}+\eta L_{yx})T\right),\quad \eta = (4/R)\textstyle{\sqrt{\Omega_{h_\calY}/\Omega'_{h_\calX}}},\label{eq:param_choice_scaled}
\end{align}
then %and $\rho'>0$ is the constant defined in Theorem~\ref{thm:main_cvx},  then %we have
\begin{align}
G(\barx^T,\bary^T) \le  B^{\rm det}_R(T)\defeq &\frac{16LR^2}{T(T-1)}\Omega'_{h_\calX} + \frac{8L_{xx}R^2}{T-1}\Omega'_{h_\calX}+ \frac{64L_{yx}R}{T-1}{\sqrt{\Omega'_{h_\calX}\Omega_{h_\calY}}} + \frac{128L_{yy}}{T}\Omega_{h_\calY}\le \mu R^2/16,\label{eq:cvx_scaled_bound}
\end{align}
and $\normt{\barx^T-x^*}\le \sqrt{(2/\mu)B^{\rm det}_R(T)}\le R/\big(2\sqrt{2}\big)$. 
\end{proposition}

\proof{Proof.}
From the choices of $\{\beta_t\}_{t=1}^T$, $\{\alpha_t\}_{t=1}^T$, $\{\tau_t\}_{t=1}^T$ and $\{\theta_t\}_{t=1}^T$, we can easily verify that the conditions~\eqref{eq:cond_set1} to~\eqref{eq:cond_set3} in Proposition~\ref{prop:recursion_cvx} continue to hold with $\gamma_t=t$, for any $t\in[T]$. % and $L$, $L_{xx}$, $L_{yx}$ and $L_{yy}$ being replaced by $L'$, $L'_{xx}$, $L'_{yx}$ and $L'_{yy}$. %(Note that since $L_{yy}=0$, the condition $\alpha_t\le (2L_{yy})^{-1}$ becomes void, i.e., not required in the step~\eqref{eq:bound_inner_product_T}.) 
In particular, $\gamma_t/\tau_t=\tau^{-1}$, for any $t\in[T]$ and $\gamma_0/\tau_0=0$. Therefore, by substituting the  parameter choices in Proposition~\ref{prop:conv_scaled} and $x^1=\hatx^1=x^0$ into~\eqref{eq:conv_whp_pseudo}, we have % now becomes
%\begin{align}
%\frac{T(T-1)}{2}\tilG(\barx^T,\bary^T;x,y) \le &\frac{2}{\tau}D_{\tilh_{\barcalX(x^1,R)}}(x,x^1) +\frac{4({T-1})\eta_y}{\alpha}\Omega_{h_\calY} + {2}\tau\sum_{t=1}^{T-1}t^2 \normt{\delta_{x,\Phi}^{t}+\delta_{x,f}^{t}}_*^2 \nn\\
%&+\sum_{t=1}^{T-1}t\lranglet{\delta_{x,\Phi}^{t}+\delta_{x,f}^t}{\hatx^t-x^{t}}.\label{eq:final_scaled}
%\end{align}
\begin{align}
\tilG(\barx^T,\bary^T;x,y) \le &\frac{2}{T(T-1)}\left\{\frac{2}{\tau}D_{\tilh_{\barcalX(x^0,R)}}(x,x^0) +\frac{4({T-1})}{\alpha}\Omega_{h_\calY}\right\},\nn\\
\le &  \left(\frac{16L}{T(T-1)}+ \frac{8L_{xx}}{T-1}\right)D_{\tilh_{\barcalX(x^0,R)}}(x,x^0)+ \frac{8L_{yx}}{T-1}\Bigg\{\eta D_{\tilh_{\barcalX(x^0,R)}}(x,x^0)+16\eta^{-1}\Omega_{h_\calY}\Bigg\}\nn\\
&\hspace{6cm}\;\;\; + \frac{128L_{yy}}{T}\Omega_{h_\calY}, \quad\forall \,T\ge 3.\label{eq:final_scaled}
\end{align}
(Note that since in~\eqref{eq:primal_upd_scaled},  the DGF $\tilh_{\barcalX(x^1,R)} = \tilh_{\barcalX(x^0,R)}$ is used, we need to replace the DGF $h_\calX$ in~\eqref{eq:conv_whp_pseudo} with $\tilh_{\barcalX(x^0,R)}$ in order to obtain~\eqref{eq:final_scaled}.)   
Next, recall from~\eqref{eq:primal_dual_function} that $\barS(x)=f(x) + g(x) + \max_{y\in\calY}\Phi(x,y) - J(y)$ and define 
\begin{align}
\underline{\hatS}_{x^0,R} (y) \defeq {\min}_{x\in\calX\bigcap\calB(x^0,R)} \;S(x,y),\quad x^*_{x^0,R}(y)\defeq {\argmin}_{x\in\calX\bigcap\calB(x^0,R)} \;S(x,y).
\end{align}
Note that since $f$ is $\mu$-strongly convex %(with respect to $\norm{\cdot}$) 
on $\calX$, the same holds for $\barS$. Based on $\barS$ and $\underline{\hatS}_{x^0,R}$, we can define the $R$-restricted duality gap 
\begin{equation}
\hatG_R(\barx^T,\bary^T)\defeq \barS(\barx^T) - \underline{\hatS}_{x^0,R} (\bary^T) = \sup\nolimits_{x\in\calX\cap\calB(x^0,R),y\in\calY} \tilG(\barx^T,\bary^T;x,y). \label{eq:def_R_gap}
\end{equation} 
The second equality in~\eqref{eq:def_R_gap} suggests us to take supremum over $x\in\calX\cap\calB(x^0,R)$ and $y\in\calY$ on both sides of~\eqref{eq:final_scaled}, and by using~\eqref{eq:diameter_Rbound} and the value of $\eta$ in~\eqref{eq:param_choice_scaled}, we have
%By the choices of $\alpha$ and $\tau$ in~\eqref{eq:param_choice_scaled}, we  have that for any $T\ge 3$, 
\begin{align}
\hspace{-.4cm}\hatG_R(\barx^T,\bary^T) \le  \frac{16LR^2}{T(T-1)}\Omega'_{h_\calX} + \frac{8L_{xx}R^2}{T-1}\Omega'_{h_\calX}+ \frac{64L_{yx}R}{T-1}{\sqrt{\Omega'_{h_\calX}\Omega_{h_\calY}}}  + \frac{128L_{yy}}{T}\Omega_{h_\calY}. 
\label{eq:final_scaled_after_sup}
\end{align}
%where the last inequality follows from the choice of $T$ in~\eqref{eq:choice_T}. 
%\begin{align}
%\hatG_R(\barx^T,\bary^T) \defeq& \; \barS(\barx^T) - \underline{\hatS}_{x^1,R} (\bary^T) \\
%\le &\;\frac{4\eta_x}{\tau T(T-1)}\Omega'_{h_\calX} +\frac{8\eta_y}{\alpha T}\Omega_{h_\calY} + \frac{4\tau}{T(T-1)}\sum_{t=1}^{T-1}t^2 (\normt{\delta_{x,\Phi}^{t}+\delta_{x,f}^{t}}'_*)^2 \nn\\
%&+\frac{2}{T(T-1)}\sum_{t=1}^{T-1}t\lranglet{\delta_{x,\Phi}^{t}+\delta_{x,f}^t}{\hatx^t-x^{t}}.
%\end{align}
On the other hand, by the $\mu$-strong convexity of $\barS$, $x^*$ is the unique minimizer of $\barS$ on $\calX$, and % we have
\begin{align}
\normt{\barx^T-x^*}^2\le \frac{2}{\mu}(\barS(\barx^T) - \barS(x^*)) \lea \frac{2}{\mu}(\barS(\barx^T) - \underline{\hatS}_{x^0,R} (\bary^T)) = \frac{2}{\mu}\hatG_R(\barx^T,\bary^T), \label{eq:x_T_dist}
\end{align}
where (a) follows from $\barS(x^*)=S(x^*,y^*)\ge S(x^*,\bary^T)\ge {\min}_{x\in\calX\bigcap\calB(x^0,R)} \;S(x,\bary^T)=\underline{\hatS}_{x^0,R} (\bary^T)$ as $x^*\in\calX\cap\calB(x^0,R)$ (note that $\normt{x^0-x^*}\le R/2$ in Proposition~\ref{prop:conv_scaled}).

Next, we aim to show that $\hatG_R(\barx^T,\bary^T)=G(\barx^T,\bary^T)$. %Fix any $y\in\calY$. 
To start, suppose that 
\begin{equation}
\normt{x^*_{x^0,R}(\bary^T)-x^*}\le R/\big(2\sqrt{2}\big).\label{eq:dist_R/2}
\end{equation}
By the input condition $\normt{x^0-x^*}\le R/2$, we have $\normt{x^*_{x^0,R}(\bary^T)-x^0}<R$. In other words, $x^*_{x^0,R}(\bary^T)\in\calX\cap\inter\calB(x^0,R)$. Note that by its definition, %for any $y\in\calY$, 
there exists $d\in \partial_x S(x^*_{x^0,R}(\bary^T),\bary^T)$
\begin{equation}
\lrangle{d}{x-x^*_{x^0,R}(\bary^T)}\ge 0, \quad \forall\,x\in\calX\cap\calB(x^0,R).\label{eq:optimality_1} 
\end{equation}
Since  $x^*_{x^0,R}(\bary^T)\in\inter\calB(x^0,R)$, for any $x\in\calX\setminus\calB(x^0,R)$, there exists $\lambda\in(0,1)$ such that $\barx\defeq \lambda x+(1-\lambda) x^*_{x^0,R}(\bary^T) \in\calB(x^0,R)$. Moreover, $\barx\in\calX$ %by the convexity of $\calX$. 
since $\calX$ is convex. Thus $\barx\in\calX\cap\calB(x^0,R)$ and
\begin{equation}
\lrangle{d}{\barx-x^*_{x^0,R}(\bary^T)}\ge 0. %, \quad \forall\,x\in\calX\cap\calB(x^0,R). 
\end{equation}
On the other hand, we have $\barx - x^*_{x^0,R}(\bary^T) = \lambda(x-x^*_{x^0,R}(\bary^T))$. Consequently, 
\begin{equation}
\lrangle{d}{x-x^*_{x^0,R}(\bary^T)}\ge 0, \quad \forall\,x\in\calX\setminus\calB(x^0,R). \label{eq:optimality_2}
\end{equation}
Combining~\eqref{eq:optimality_1} and~\eqref{eq:optimality_2}, we have
\begin{equation}
\lrangle{d}{x-x^*_{x^0,R}(\bary^T)}\ge 0, \quad \forall\,x\in\calX. \label{eq:optimality_final}
\end{equation}
This indicates that $x^*_{x^0,R}(\bary^T)=\argmin_{x\in\calX} S(x,\bary^T)$ and hence $\hatG_R(\barx^T,\bary^T)=G(\barx^T,\bary^T)$. 

Now, it remains to show that~\eqref{eq:dist_R/2} holds. First, since $S(\cdot,\bary^T)$ is $\mu$-strongly convex on $\calX$, by the definition of $x^*_{x^0,R}(\bary^T)$ and the fact that $x^*\in \calX\cap\calB(x^0,R)$, we have %that for any $d\in\partial_x S(x^*_{x^0,R}(\bary^T),\bary^T)$, 
\begin{align}
S(x^*,\bary^T) - \underline{\hatS}_{x^0,R} (\bary^T)= S(x^*,\bary^T) - S(x^*_{x^0,R}(\bary^T),\bary^T) %\lrangle{d}{x^* - x^*_{x^0,R}(\bary^T)} + \frac{\mu}{2}\normt{x^* - x^*_{x^0,R}(\bary^T)}^2
\ge \frac{\mu}{2}\normt{x^*_{x^0,R}(\bary^T)-x^*}^2. \label{eq:sc_bound1}
\end{align}
On the other hand, $S(x^*,\bary^T)\le \max_{y\in\calY} S(x^*,y) = \barS(x^*) \le \barS(\barx^T)$, since $x^*$ minimizes $\barS$ on $\calX$. Thus
\begin{equation}
\hatG_R(\barx^T,\bary^T) = \barS(\barx^T) - \underline{\hatS}_{x^0,R} (\bary^T) \ge \frac{\mu}{2}\normt{x^*_{x^0,R}(\bary^T)-x^*}^2. \label{eq:sc_bound2}
\end{equation}
However, note that from~\eqref{eq:final_scaled_after_sup} and the choice of $T$ in~\eqref{eq:choice_T}, we have $\hatG_R(\barx^T,\bary^T) \le \mu R^2/16$. We hence complete the proof. \Halmos
\endproof

From Proposition~\ref{prop:conv_scaled}, we observe that $\normt{x^0-x^*}\le R/2$ and $\normt{\barx^T-x^*}\le R/\big(2\sqrt{2}\big)$. This suggests that if Algorithm~\ref{algo:convex_f_scaled} is used as the subroutine in a restart scheme, then at each stage, the radius $R$ %(i.e., the distance from ) 
can be reduced by a factor of $\sqrt{2}$, and accordingly, the bound of the duality gap (i.e., $\mu R^2/16$) can be halved. 
This observation naturally leads us to the restart scheme in Algorithm~\ref{algo:restart_det}, which comprises $K$ stages. At each stage $k$,  using the output primal variable $\barx_{k-1}^{T_{k-1}}$ from the previous stage as the input, we run Algorithm~\ref{algo:convex_f_scaled} for a sufficiently large number of iterations (i.e., $T_k$ iterations), so as to ensure the output primal variable $\barx_k^{T_k}$ in the current stage satisfies that 
\begin{align}
\normt{\barx_k^{T_k}-x^*}\le R_k/(2\sqrt{2}). \label{eq:suff_close}
\end{align}
%$(\barx_k^{T_k},\bary_k^{T_k})$ satisfies 

\renewcommand\thealgorithm{2D}
\begin{algorithm}[t!]
\caption{Deterministic restart scheme for strongly convex $f$} \label{algo:restart_det}
\begin{algorithmic}
\State {\bf Input}: Diameter estimate $U\ge D_\calX$,  desired accuracy $\epsilon>0$, $K=\big\lceil\!\max\big\{0,\log_2\big(\mu U^2/(4\epsilon)\big)\big\}\big\rceil+1$
\State {\bf Initialize}: $R_1=2U$, $x^1\in\calX^o$ %, $y_0\in\calY^o$
\State {\bf For} $k = 1,\ldots,K$
\vspace{-.2cm}
%\begin{enumerate}
\begin{align}
\hspace{-.1cm}1. \;\;T_k:= \Big\lceil\!\max\Big\{3,\;64\textstyle{\sqrt{{(L/\mu)\Omega'_{h_\calX}}}},&\;{1024(L_{xx}/\mu)\Omega'_{h_\calX}},\nn\\
&\hspace{-.5cm}4096L_{yx}(\mu R_k)^{-1}{\sqrt{\Omega'_{h_\calX}\Omega_{h_\calY}}},\;8192L_{yy}(\mu R_k^2)^{-1}\Omega_{h_\calY}\Big\}\Big\rceil. \label{eq:choice_T_k}
\end{align}
%\end{enumerate}
\vspace{-.8cm}
\State \quad\, 2. Run Algorithm~\ref{algo:convex_f_scaled} for $T_k$ iterations with starting primal variable $x_k$, radius $R_k$, constraint set $\calX_k\equiv \calX$ and other input parameters set as in Proposition~\ref{prop:conv_scaled}. Denote the output as $(\barx_k^{T_k},\bary_k^{T_k})$. 
\State \quad\, 3. $R_{k+1}:=R_k/\sqrt{2}$, $x_{k+1}:=\barx_k^{T_k}$.
\State {\bf Output}: $(x_{K+1},y_{K+1})$ %$\barbx^K$ and  $\barby^K$
\end{algorithmic}%\vspace{-.2cm}
\end{algorithm}

\begin{theorem}\label{thm:restart_det}
Let Assumption~\ref{assump:bounded}\ref{assump:bounded_sets} hold.  In Algorithm~\ref{algo:restart_det}, for any desired accuracy $\epsilon\in(0,\mu U^2/4]$, we have $G(x_{K+1},y_{K+1}) \le \epsilon$, and the total number of oracle calls $C^{\rm det}_\epsilon$ satisfies that 
\begin{align}
C^{\rm det}_\epsilon =\sum_{k=1}^K T_k \le &\Big(3+64\textstyle{\sqrt{{(L/\mu)\Omega'_{h_\calX}}}}+{1024(L_{xx}/\mu)\Omega'_{h_\calX}}\Big) \left(\big\lceil\log_2\big(\mu U^2/(4\epsilon)\big)\big\rceil+1\right)\nn\\
&+8192\big(L_{yx}/\textstyle{\sqrt{\mu\epsilon}}\big)\sqrt{\Omega'_{h_\calX}\Omega_{h_\calY}} + 2048\big(L_{yy}/\epsilon\big)\Omega_{h_\calY}.\label{eq:complexity_det_sc} 
\end{align}
\end{theorem}

Before presenting the proof, from Theorem~\ref{thm:restart_det}, we know that the oracle complexity  for Algorithm~\ref{algo:restart_det} to obtain an $\epsilon$-duality gap is
\begin{equation}
O\left(\big(\sqrt{{L}/{\mu}}+{L_{xx}}/{\mu}\big)\log\left({1}/{\epsilon}\right)+ {L_{yx}}/{\sqrt{\mu\epsilon}} + {L_{yy}}/{\epsilon}\right). 
\end{equation}
 
\proof{Proof.}
Note that $R_k=2^{(3-k)/2}U$, for any $k\in[K]$. Therefore, 
\begin{equation}
G(x_{K+1},y_{K+1}) = G(\barx_K^{T_K},\bary_K^{T_K})\lea \mu R_K^2/16 = \mu U^2 2^{-(K+1)}\leb \epsilon, \label{eq:bound_G_K+1}
\end{equation}
where (a) follows from~\eqref{eq:cvx_scaled_bound} and (b) follows from $K\ge \log_2\big(\mu U^2/(4\epsilon)\big)+1$.  %in the input, we have 
In addition, we can also substitute the value of $R_k$ into~\eqref{eq:choice_T_k} and obtain~\eqref{eq:complexity_det_sc}. \Halmos
\endproof

\begin{remark}
By restricting $\epsilon\in(0,\mu U^2/4]$, we see that  $\max\big\{0,\log_2\big(\mu U^2/(4\epsilon)\big)\big\}$ simply becomes $\log_2\big(\mu U^2/(4\epsilon)\big)$. By doing so, we have indeed simplified the bound in~\eqref{eq:complexity_det_sc} (as compared to the bound derived for $\epsilon>0$). 
On the other hand, note that if $\epsilon>\mu U^2/4$, then $K=1$, and Algorithm~\ref{algo:restart_det} effectively becomes Algorithm~\ref{algo:convex_f_scaled}. For the same reason, we will also focus on analyzing the regime $\epsilon\in(0,\mu U^2/4]$ in the stochastic restart scheme (see Theorem~\ref{thm:restart_stoc}). 
\end{remark}

\begin{remark}\label{rmk:abs_cosntant}
We observe that in Algorithm~\ref{algo:restart_det}, the choice of $T_k$ (cf.~\eqref{eq:choice_T_k}) involves some large absolute constants, e.g., 4096 and 8192. These constants are mainly for theoretical purposes, and not meant to be tight. In practice, to avoid running Algorithm~\ref{algo:convex_f_scaled} for an excessively large number of iterations at each stage $k$ of Algorithm~\ref{algo:restart_det}, we can simply replace all the absolute constants in~\eqref{eq:choice_T_k} (except 3) by a single constant $c>0$, and choose $c$ accordingly. The same remark also applies to the stochastic restart scheme (i.e.,~Algorithm~\ref{algo:restart_stoc}) in Section~\ref{sec:algo_stoc}.    %of our restart scheme (i.e., ). we can simply 
\end{remark}

%\begin{remark}

%\end{remark}

%Theorem~\ref{thm:restart_det} indicates that to obtain an $\epsilon$-duality gap, the oracle complexity is 
%\begin{equation}
%O\left({\sqrt{\frac{L}{\mu}}}\log\bigg(\frac{1}{\epsilon}\bigg)+
%\frac{L_{xx}}{\mu}\log\bigg(\frac{1}{\epsilon}\bigg) + \frac{L_{yx}}{\sqrt{\mu\epsilon}} + \frac{L_{yy}}{\epsilon}\right).
%\end{equation}

\subsection{Stochastic restart scheme.}\label{sec:stoc_restart}
We now consider the general case where $\sigma_{x,f},\sigma_{x,\Phi},\sigma_{y,\Phi}>0$. %we only have access to the {\em stochastic} gradients of $f$, $\Phi(\cdot,y)$ and $\Phi(x,\cdot)$. 

\subsubsection{Intuition.} \label{sec:intuition}
At the first attempt, one may try to combine the techniques used in the proofs of Theorem~\ref{thm:main_cvx}\ref{res:cvx_exp} and Proposition~\ref{prop:conv_scaled}, and then analyze the convergence of $\bbE[G(\barx^T,\bary^T)]$ in Algorithm~\ref{algo:convex_f_scaled}. However, a close inspection shows that such a combination does not work. Indeed, with this combination, in the proof of Proposition~\ref{prop:conv_scaled}, although one can ensure that $\bbE[\hatG_R(\barx^T,\bary^T)] \le \mu R^2/16$ and hence $\bbE\big[\normt{x^*_{x^0,R}(\bary^T)-x^*}\big]\le R/\big(2\sqrt{2}\big)$ (cf.~\eqref{eq:sc_bound2}) by choosing $T$ properly, these conditions cannot guarantee that $\bbE[\hatG_R(\barx^T,\bary^T)]=\bbE[G(\barx^T,\bary^T)]$. This is because it is unclear how to  ensure that  $x^*_{x^0,R}(\bary^T)\in\calX\cap\inter\calB(x^0,R)$ holds in the ``expectation'' sense, in contrast to the deterministic setting, where this condition holds a.s. With that said, when Assumption~\ref{assump:noise}\ref{assump:sub-Gaussian} holds (i.e., the gradient noises are sub-Gaussian), it is  possible to show that  $\hatG_R(\barx^T,\bary^T) \le \mu R^2/16$ (and hence all the rest steps in the proof of Proposition~\ref{prop:conv_scaled}, including $G(\barx^T,\bary^T) \le \mu R^2/16$) holds {\em with high probability} (cf.~Theorem~\ref{thm:main_cvx}\ref{res:cvx_LD}). Moreover, if we ``concatenate'' this result in our  restart scheme, we will still end up with a high-probability bound on the duality gap, as long as we keep the number of stages ``reasonably'' small.

\subsubsection{Algorithmic details.} \label{sec:algo_stoc}
%To start with, note that under the geometry in Section~\ref{sec:stoc_restart}, Assumption~\ref{assump:noise}\ref{assump:sub-Gaussian} now becomes
%\begin{align}
%&\bbE_{t-1}\left[\exp\left\{(\normt{\delta_{y,\Phi}^{t}}'_*)^2/(\sigma'_{y,\Phi})^2\right\}\right] \le \exp(1) \;\;{\rm a.s.},\\ 
%&\bbE_{t-1}\left[\exp\left\{(\normt{\delta_{x,\Phi}^{t}}'_*)^2/(\sigma'_{x,\Phi})^2\right\}\right] \le \exp(1) \;\;{\rm a.s.},\\
%&\bbE_{t-1}\left[\exp\left\{(\normt{\delta_{x,f}^{t}}'_*)^2/(\sigma'_{x,f})^2\right\}\right] \le \exp(1) \;\;{\rm a.s.},
%\end{align}
%where $\sigma'_{y,\Phi} \defeq \sigma_{y,\Phi}/\sqrt{\eta_y}$, $\sigma'_{x,\Phi} \defeq (R/\sqrt{\eta_x})\sigma_{x,\Phi}$ and $\sigma'_{x,f} \defeq (R/\sqrt{\eta_x})\sigma_{x,f}$. Based on these definitions, 
Let us first analyze the convergence of Algorithm~\ref{algo:convex_f_scaled} in the stochastic setting. 

\begin{proposition}\label{prop:conv_scaled_stoc}
Assume that $ \calB(0,1)\subseteq \dom h_\calX$, and let Assumptions~\ref{assump:bounded}\ref{assump:bounded_sets},~\ref{assump:noise}\ref{assump:unbiased} and~\ref{assump:noise}\ref{assump:sub-Gaussian} hold. %In addition, assume that  
Fix any $\varsigma\in(0,1/6]$. In Algorithm~\ref{algo:convex_f_scaled}, choose  $\calX'$ such that $x^*\in\calX'$ and  $D_{\calX'}\le R$,  and choose 
\begin{align}
T\ge \Big\lceil\!\max\Big\{&3,\;64\textstyle{\sqrt{{(L/\mu)\Omega'_{h_\calX}}}},\;{2048(L_{xx}/\mu)\Omega'_{h_\calX}},\;4096L_{yx}(\mu R)^{-1}{\sqrt{\Omega'_{h_\calX}\Omega_{h_\calY}}},\;128^2L_{yy}(\mu R^2)^{-1}\Omega_{h_\calY},\nn\\
&\;512^2(\sigma_{x,f}+\sigma_{x,\Phi})^2(\mu R)^{-2}\big(4\textstyle{\sqrt{(1+\log(1/\varsigma))\Omega'_{h_\calX}}}+2\textstyle{\sqrt{\log(1/\varsigma)}}\big)^2,\nn\\
&\;512^2\sigma_{y,\Phi}^2(\mu R^2)^{-2}\big(8\textstyle{\sqrt{2(1+\log(1/\varsigma))\Omega_{h_\calY}}} + {2\textstyle{\sqrt{\log(1/\varsigma)\Omega_{h_\calY}}}}\big)^2\;\Big\}\Big\rceil. \label{eq:choice_T_stoc}
\end{align}
 If we choose $R\ge 2\normt{x^0-x^*}$, $\{\beta_t\}_{t=1}^{T}$ and $\{\theta_t\}_{t=1}^T$ as in~\eqref{eq:param_choice_1}, and $\alpha_t=\alpha$ and $\tau_t=t\tau$ for any $t\in[T]$, where
\begin{equation}
\begin{aligned}
&\alpha = 1/\big(16\big(\eta^{-1}L_{yx} + L_{yy}+\rho \sigma_{y,\Phi}\sqrt{T}\big)\big), & \rho = (4R)^{-1}\sqrt{(1+\log(1/\varsigma))/(2\Omega'_{h_\calX}\Omega_{h_\calY})},\\
&\tau = 1/\big(4L+2(L_{xx}+\eta L_{yx})T+\rho'(\sigma_{x,\Phi} + \sigma_{x,f})T^{3/2}\big), & \rho'=(8R)^{-1}\sqrt{(1+\log(1/\varsigma))/(\Omega'_{h_\calX}\Omega_{h_\calY})},
\end{aligned}%\label{eq:param_choice_scaled_stoc}
\end{equation}
and $\eta = (4/R)\textstyle{\sqrt{\Omega_{h_\calY}/\Omega'_{h_\calX}}}$, then with probability at least $1-6\varsigma$,  %and $\rho'>0$ is the constant defined in Theorem~\ref{thm:main_cvx},  then %we have
\begin{align}
G(\barx^T,\bary^T) \le  B^{\rm det}_R(T) + B^{\rm var}_R(T)\le \mu R^2/16,\label{eq:cvx_scaled_bound_stoc}
\end{align}
where $B^{\rm det}_R(T)$ is defined in~\eqref{eq:cvx_scaled_bound} and 
\begin{align}
B^{\rm var}_R(T)\defeq\frac{4(\sigma_{x,\Phi}+\sigma_{x,f}) R}{\sqrt{T}}&\left\{4\textstyle{\sqrt{(1+\log(1/\varsigma))\Omega'_{h_\calX}}}+2\sqrt{\log(1/\varsigma)}\right\} \nn\\
&+ \frac{4\sigma_{y,\Phi}}{\sqrt{T}}\Big\{8\textstyle{\sqrt{2(1+\log(1/\varsigma))\Omega_{h_\calY}}} + {2\sqrt{\log(1/\varsigma)\Omega_{h_\calY}}}\Big\}.
\end{align}
Furthermore, $\normt{\barx^T-x^*}\le \sqrt{(2/\mu)(B^{\rm det}_R(T) + B^{\rm var}_R(T))}\le R/\big(2\sqrt{2}\big)$ with probability at least $1-6\varsigma$. 
\end{proposition}

\proof{Proof.}
%Following the proof of Theorem~\ref{thm:main_cvx}\ref{res:cvx_LD}, by the assumption that $D_{\calX'}\le R$ and the choices of $\{\beta_t\}_{t\in\bbN}$, $\{\alpha_t\}_{t\in\bbN}$, $\{\tau_t\}_{t\in\bbN}$ and $\{\theta_t\}_{t\in\bbN}$ in Proposition~\ref{prop:conv_scaled_stoc}, we can easily show that 
The proof is adapted from that of Proposition~\ref{prop:conv_scaled}, but with a careful stochastic argument. 
First, we substitute the values of $\{\beta_t\}_{t=1}^T$, $\{\alpha_t\}_{t=1}^T$, $\{\tau_t\}_{t=1}^T$ and $\{\theta_t\}_{t=1}^T$ into~\eqref{eq:conv_whp_pseudo} and then take supremum over $x\in\calX\cap\calB(x^0,R)$ and $y\in\calY$ on both sides of the resulting inequality. By noting that $D_{\calX'}\le R$, we have that with probability at least $1-6\varsigma$,
\begin{align}
\hatG_R(\barx^T,\bary^T) \le  B^{\rm det}_R(T) + B^{\rm var}_R(T).\label{eq:bound_G_R} %\label{eq:cvx_scaled_bound_stoc}
\end{align}
% This immediately lead to the (high-probability) bound on $\normt{\barx^T-x^*}$ from~\eqref{eq:x_T_dist}. 
Denote the event in~\eqref{eq:bound_G_R} as $\calA_T$, on which we condition till the end of the proof. By the choice of $T$ in~\eqref{eq:choice_T_stoc}, we have with probability 1 that $\hatG_R(\barx^T,\bary^T)\le \mu R^2/16$. Consequently, from~\eqref{eq:x_T_dist} and~\eqref{eq:sc_bound2}, we see that $\normt{\barx^T-x^*}\le R/\big(2\sqrt{2}\big)$ and $\normt{x^*_{x^0,R}(\bary^T)-x^*}\le R/\big(2\sqrt{2}\big)$ respectively with probability~1. As a result, we conclude that $x^*_{x^0,R}(\bary^T)\in\calX\cap\inter\calB(x^0,R)$ and hence $\hatG_R(\barx^T,\bary^T) = G(\barx^T,\bary^T)$. \Halmos %We then complete the proof. 
%the first inequality in~\eqref{eq:cvx_scaled_bound_stoc} holds. Moreover, by the choice of $T$ in~\eqref{eq:choice_T_stoc}, we see that the second inequality in~\eqref{eq:cvx_scaled_bound_stoc} also holds. 
\endproof

Based on Proposition~\ref{prop:conv_scaled_stoc}, we present our stochastic restart scheme in Algorithm~\ref{algo:restart_stoc}. Compared to the deterministic restart scheme (i.e., Algorithm~\ref{algo:restart_det}), a notable difference is that at each stage $k$, the constraint set $\calX'$ becomes $\calX\cap\calB(x_k,R_k/2)$ (as opposed to $\calX$ as in Algorithm~\ref{algo:restart_det}). %the intersection of $\calX\cap\calB(x_k,R_k/2)$ and $\calB^{x_k,R_k/2}_{\norm{\cdot}}$. 
As will be shown in Theorem~\ref{thm:restart_stoc}, this step enables us to obtain the nearly optimal %$O(\log\log(1/\epsilon)/\epsilon)$ 
oracle complexity on the primal noise term $\sigma_{x,f}+\sigma_{x,\Phi}$. Note that in general, the minimization problem in~\eqref{eq:primal_upd_scaled} with constraint set $\calX'=\calX\cap\calB(x_k,R_k/2)$ is harder to solve, compared to the problem with constraint set $\calX$. However, in many special cases, the problem with constraint set $\calX'$ can still be easily solved. 
For instance, let $\bbX$ be a Hilbert space and $h_\calX = (1/2)\norm{\cdot}^2$, so that $D_{\tilh_{\barcalX(x_\rmc,R)}}(x,x')=(1/2)\norm{x-x'}^2$.  We can rewrite $\calB(x_k,R_k/2)$ as the functional constraint $\normt{x-x_k}\le R_k/2$, and then consider the Lagrangian form of this problem, i.e., 
\begin{align}
\min_{x\in\calX}\max_{\lambda\ge 0} \; g(x)+\lranglet{\pi^t}{x-x^t} + ({2\tau_t})^{-1}\normt{x-x^t}^2 + \lambda\big(\normt{x-x_k}^2-R_k^2/4\big),\label{eq:Lagrangian_form}
\end{align}
where $\pi^t\defeq\hat{\nabla}_x\Phi(x^t,y^{t+1},\zeta^t_x)+\hat{\nabla}f(\tilx^{t+1},\xi^t)$. We observe that for any fixed $\lambda\ge 0$, the minimization problem in~\eqref{eq:Lagrangian_form} has the same form as that in~\eqref{eq:Bregman_projection}, hence it has an easily computable solution. Furthermore, since $\lambda$ is a scalar, we can optimize it over $[0,+\infty)$ efficiently (e.g., via Newton's method or bisection).
%(However, since we only have one functional constraint $\normt{x-x_k}\le R_k/2$, this problem can still be solved relatively easily.) 
%%, we can consider its Lagrangian form, i.e., 
%\begin{align}
%\min_{x\in\calX}\max_{\lambda\ge 0} \; g(x)+\lranglet{\pi^t}{x-x^t} + \frac{\eta_x}{2\tau_tR_k^2}\normt{x-x^t}^2 + \lambda\big(\normt{x-x_k}^2-R_k^2/4\big),\label{eq:Lagrangian_form}
%\end{align}
%where $\pi^t\defeq\hat{\nabla}_x\Phi(x^t,y^{t+1},\zeta^t_x)+\hat{\nabla}f(\tilx^{t+1},\xi^t)$. We observe that for any fixed $\lambda\ge 0$, the minimization problem in~\eqref{eq:Lagrangian_form} has the same form as that in~\eqref{eq:Bregman_projection}, hence has easily computable solution. Since $\lambda$ is a scalar, we can solve the maximization problem in~\eqref{eq:Lagrangian_form} efficiently (e.g., via bisection). 
%In addition, for some special cases, e.g., $\calX\subseteq\calB(x_k,R_k/2)$ or $\calB(x_k,R_k/2)\subseteq\calX$, we can directly compute the solution of~\eqref{eq:primal_upd_scaled} without any iterative procedure. %has a 
%To discuss the computational issue of the step~\eqref{eq:primal_upd_scaled} in Algorithm~\ref{algo:convex_f_scaled}, suppose $\bbX$ is a Hilbert space and take $h_\calX = (1/2)\norm{\cdot}_{\bbX}^2$. If $g\equiv 0$ and $\calX=\bbX$, then~\eqref{eq:primal_upd_scaled} becomes an orthogonal projection onto  $\calB^{x_k,R_k/2}_{\norm{\cdot}}$, which clearly admits closed-form solution. If not,
%To solve this problem, for any fixed $\lambda\ge 0$, note that 

To analyze the convergence of Algorithm~\ref{algo:restart_stoc}, let us define the event
\begin{align}
\calE_1 = \Omega,\quad \calE_k\defeq \big\{\hatG_{R_{k-1}}(x_{k},y_{k})\le \mu R_{k-1}^2/16\big\}, \quad \forall\,k=2,\ldots,K+1,  \label{eq:E_k}
\end{align}
where $\Omega$ denotes the sample space for the stochastic process $\{(x_k,y_k)\}_{k=1}^{K+1}$. 
For any $k\ge 2$, by conditioning on $\calE_{k}$, from the proof of Proposition~\ref{prop:conv_scaled_stoc}, we see that both $\normt{x_k-x^*}\le R_k/2$ (since $R_k=R_{k-1}/\sqrt{2}$) a.s.\ and 
\begin{equation}
G(x_k,y_k)=\hatG_{R_{k-1}}(x_k,y_k)\le \mu R_{k-1}^2/16 \quad\quad \mbox{a.s.} \label{eq:duality_gap_wp1}
\end{equation} 
Therefore, given $\calE_k$, $\calX_k$ satisfies all the requirements stated in Algorithm~\ref{algo:convex_f_scaled} and Proposition~\ref{prop:conv_scaled_stoc}, i.e., $x^*\in\calX_k\subseteq\calX$ and $D_{\calX_k}\le R_k$. 
(Note that for $k=1$, we have $\calX_1=\calX$, thus these requirements are clearly satisfied.)
As a result, Proposition~\ref{prop:conv_scaled_stoc} can be applied to the $k$-th stage of Algorithm~\ref{algo:restart_stoc}. Based on this observation, we derive the oracle complexity of Algorithm~\ref{algo:restart_stoc} below.

\renewcommand\thealgorithm{2S}
\begin{algorithm}[t!]
\caption{Stochastic restart scheme for strongly convex $f$} \label{algo:restart_stoc}
\begin{algorithmic}
\State {\bf Input}: Diameter estimate $U\ge D_\calX$, starting primal variable $x_0\in\calX^o$, desired accuracy $\epsilon>0$, error probability $\nu\in(0,1]$, $K=\big\lceil\max\big\{0,\log_2\big(\mu U^2/(4\epsilon)\big)\big\}\big\rceil+1$, $\varsigma=\nu/(6K)$
\State {\bf Initialize}: $R_1=2U$, $x_1=x_0$, $y_0\in\calY^o$
\State {\bf For} $k = 1,\ldots,K$
\vspace{-.2cm}
%\begin{enumerate}
\begin{align}
\hspace{-.1cm}1. \;\;&T_k:= \Big\lceil\!\max\Big\{3,\;64\textstyle{\sqrt{{(L/\mu)\Omega'_{h_\calX}}}},\;{2048(L_{xx}/\mu)\Omega'_{h_\calX}},\;4096L_{yx}(\mu R_k)^{-1}{\sqrt{\Omega'_{h_\calX}\Omega_{h_\calY}}},\nn\\
&\;128^2L_{yy}(\mu R_k^2)^{-1}\Omega_{h_\calY},\;512^2(\sigma_{x,f}+\sigma_{x,\Phi})^2(\mu R_k)^{-2}\big(4\textstyle{\sqrt{(1+\log(1/\varsigma))\Omega'_{h_\calX}}}+2\textstyle{\sqrt{\log(1/\varsigma)}}\big)^2,\nn\\
&\;512^2\sigma_{y,\Phi}^2(\mu R_k^2)^{-2}\big(8\textstyle{\sqrt{2(1+\log(1/\varsigma))\Omega_{h_\calY}}} + {2\textstyle{\sqrt{\log(1/\varsigma)\Omega_{h_\calY}}}}\big)^2\;\Big\}\Big\rceil. \label{eq:choice_T_k_stoc}
\end{align}
%\end{enumerate}
\vspace{-.8cm}
\State \quad\, 2. Run Algorithm~\ref{algo:convex_f_scaled} for $T_k$ iterations with starting primal variable $x_k$, radius $R_k$, constraint set $\calX_k= \{x\in\calX:\normt{x-x_k}\le R_k/2\}$ and other input parameters set as in Proposition~\ref{prop:conv_scaled_stoc}. Denote the output as $(\barx_k^{T_k},\bary_k^{T_k})$. 
\State \quad\, 3. $R_{k+1}:=R_k/\sqrt{2}$, $x_{k+1}:=\barx_k^{T_k}$.
\State {\bf Output}: $(x_{K+1},y_{K+1})$ %$\barbx^K$ and  $\barby^K$
\end{algorithmic}%\vspace{-.2cm}
\end{algorithm}

\begin{theorem}\label{thm:restart_stoc}
Assume that $ \calB(0,1)\subseteq \dom h_\calX$. Also, let Assumptions~\ref{assump:bounded}\ref{assump:bounded_sets},~\ref{assump:noise}\ref{assump:unbiased} and~\ref{assump:sub-Gaussian} hold. 
In Algorithm~\ref{algo:restart_stoc}, for any $x_0\in\calX^o$, desired accuracy $\epsilon\in(0,\mu U^2/4]$ and error probability $\nu\in(0,1]$, it holds that $G(x_{K+1},y_{K+1}) \le \epsilon$ with probability at least $1-\nu$. Furthermore, the total number of oracle calls $C^{\rm st}_\epsilon$ in Algorithm~\ref{algo:restart_stoc} satisfies that %i.e.,
\begin{align}
C^{\rm st}_\epsilon  \le& \Big(3+64\textstyle{\sqrt{{(L/\mu)\Omega'_{h_\calX}}}}+{2048(L_{xx}/\mu)\Omega'_{h_\calX}}\Big) \left(\big\lceil\log_2\big(\mu U^2/(4\epsilon)\big)\big\rceil+1\right)\nn\\
&+256^2\big(L_{yx}/\textstyle{\sqrt{\mu\epsilon}}\big)\sqrt{\Omega'_{h_\calX}\Omega_{h_\calY}} + 64^2\big(L_{yy}/\epsilon\big)\Omega_{h_\calY}\nn\\
&+1024^2\left\{(\sigma_{x,f}+\sigma_{x,\Phi})^2/(\epsilon\mu)\right\}\big\{(4\Omega'_{h_\calX}+1)\log\left(6\left[\log_2\big(\mu U^2(4\epsilon)^{-1}\big)+2\right]/\nu\right)+4\Omega'_{h_\calX}\big\}\nn\\
&+1024^2\left(\sigma_{y,\Phi}^2/\epsilon^2\right)\big\{1+\log\left(6\left[\log_2\big(\mu U^2(4\epsilon)^{-1}\big)+2\right]/\nu\right)\big\}\Omega_{h_\calY}.  \label{eq:complexity_stoc_sc}
\end{align}
\end{theorem}

Before presenting the proof, from Theorem~\ref{thm:restart_stoc}, we see that the oracle complexity of Algorithm~\ref{algo:restart_stoc} to obtain an $\epsilon$-duality gap with probability at least $1-\nu$ is 
\begin{equation}
O\left(\bigg(\sqrt{\frac{L}{\mu}}+\frac{L_{xx}}{\mu}\bigg)\log\left(\frac{1}{\epsilon}\right)+ \frac{L_{yx}}{\sqrt{\mu\epsilon}} + \frac{L_{yy}}{\epsilon}+ \left(\frac{(\sigma_{x,f}+\sigma_{x,\Phi})^2}{\mu\epsilon}+\frac{\sigma_{y,\Phi}^2}{\epsilon^2}\right)\log\left(\frac{\log(1/\epsilon)}{\nu}\right)\right).\label{eq:complexity_sc_hp}
\end{equation}

\proof{Proof.}
For any $k=1,\ldots,K+1$, let $\bbI_{\calE_k}$ denote the indicator function of the event $\calE_k$ (cf.~\eqref{eq:E_k}). It is clear that $\{\bbI_{\calE_k}\}_{k=1}^{K+1}$ forms a (finite-horizon) Markov chain, and therefore
\begin{align}
\Pr\big\{\textstyle{\bigcap}_{k=1}^{K+1}\;\calE_k\big\} &= \Pr\big\{\bbI_{\calE_k}=1,\,\forall\,k=2,\ldots,K+1\big\}\nn\\
 &=\bbP\big\{\bbI_{\calE_2}=1\big\}\prod_{k=3}^{K+1}  \Pr\big\{\bbI_{\calE_k}=1\,\big|\,\bbI_{\calE_{k-1}}=1\big\}\nn\\
&\gea (1-6\varsigma)^{K}\geb 1-6K\varsigma \eqc 1-\nu,
\end{align}
where (a) follows from Proposition~\ref{prop:conv_scaled_stoc}, (b) follows from Bernoulli's inequality and (c) follows from the choice of $\varsigma$ in Algorithm~\ref{algo:restart_stoc}. 
By the choice of $\{R_k\}_{k=1}^{K}$ in Algorithm~\ref{algo:restart_stoc}, we have  
\begin{equation}
R_k=2^{(3-k)/2}U, \;\forall\, k\in[K].\label{eq:R_k}
\end{equation}
Therefore, from~\eqref{eq:bound_G_K+1}, we know that $\mu R_K^2/16 \le \epsilon$. As a result,
\begin{align}
\Pr\big\{G(x_{K+1},y_{K+1})\le \epsilon\big\}&\ge \Pr\big\{G(x_{K+1},y_{K+1})\le \mu R_{K}^2/16\big\}\nn\\
&\ge \Pr\big\{G(x_{K+1},y_{K+1})\le \mu R_{K}^2/16\,\big|\,\calE_{K+1}\big\}\Pr\big\{\calE_{K+1}\big\}\nn\\
&\eqa \Pr\big\{\calE_{K+1}\big\} \ge \Pr\big\{\textstyle{\bigcap}_{k=1}^{K+1}\;\calE_k\big\}\ge  1-\nu, 
\end{align}
where (a) follows from~\eqref{eq:duality_gap_wp1}. From the choice of $\varsigma$ and $K$, we have 
\begin{align}
\log(1/\varsigma) = \log(6K/\nu)\le \log\left(6\left[\log_2\big(\mu U^2(4\epsilon)^{-1}\big)+2\right]/\nu\right).\label{eq:log_1/varsigma}
\end{align}
Now, we substitute both~\eqref{eq:R_k} and~\eqref{eq:log_1/varsigma} into the choice of $T_k$ in~\eqref{eq:choice_T_k_stoc}, and then obtain the oracle complexity in~\eqref{eq:complexity_stoc_sc}. \Halmos
\endproof

%From Theorem~\ref{thm:restart_stoc}, we see that in order to obtain an $\epsilon$-duality gap  with probability at least $1-\nu$, the oracle complexity is 
%\begin{equation}
%O\left(\bigg(\sqrt{\frac{L}{\mu}}+\frac{L_{xx}}{\mu}\bigg)\log\left(\frac{1}{\epsilon}\right)+ \frac{L_{yx}}{\sqrt{\mu\epsilon}} + \frac{L_{yy}}{\epsilon}+ \left(\frac{(\sigma_{x,f}+\sigma_{x,\Phi})^2}{\mu\epsilon}+\frac{\sigma_{y,\Phi}^2}{\epsilon^2}\right)\log\left(\frac{\log(1/\epsilon)}{\nu}\right)\right).\label{eq:complexity_sc_hp}
%\end{equation}

\subsubsection{Complexity of convergence in expectation.} \label{sec:conv_expectation}
Based on the complexity results in Theorem~\ref{thm:restart_stoc}, we aim to analyze the oracle complexity to obtain an $\varepsilon$-{\em expected} duality gap. To do so, we need one additional assumption on the nonsmooth functions $g$ and $J$.
%result for $\bbE[G(\barx_{K+1},\bary_{K+1})]$

\begin{assumption}\label{assump:regularity_g_J}
The nonsmooth functions $g$ and $J$ have closed domains. %, and are continuous on $\calX$ and $\calY$, respectively. 
%\begin{enumerate}[label={\rm (\Alph*)}]
%\item The function $g$ has closed domain and continuous on $\calX$. 
%\item The function $g$ has closed domain and continuous on $\calX$. 
%\end{enumerate}
\end{assumption}

\begin{remark}
Note that Assumption~\ref{assump:regularity_g_J} % the assumption of closed domains 
is satisfied by many nonsmooth functions, e.g., an indicator function of a closed convex set or an absolutely homogeneous function (e.g., the norm function). 
%In addition, since $g$ and $J$ are closed and convex, if $\calX$ and $\calY$ are contained in $\dom g$ and $\dom J$, respectively, then the continuity assumption is satisfied. 
\end{remark}

To see the implication of Assumption~\ref{assump:regularity_g_J}, let us write down the explicit form of the primal and dual functions $\barS$ and $\uS$ in~\eqref{eq:primal_dual_function}:
\begin{align}
\barS(x) &= f(x) + g(x) + \big\{{\max}_{y\in\calY\bigcap\dom J} \Phi(x,y) - J(y)\big\},\label{eq:barS}\\
\uS(y) &= \big\{{\min}_{x\in\calX\bigcap\dom g} f(x) + g(x) + \Phi(x,y)\big\} - J(y). \label{eq:uS}
\end{align}
In~\eqref{eq:barS}, by the closedness of $\dom J$ and the compactness of $\calY$ (cf.\ Assumption~\ref{assump:bounded}\ref{assump:bounded_sets}), we see that $\dom J\cap\calY$ is compact.  Furthermore, since $J$ is closed and convex, it is continuous on $\dom J\cap\calY$. Therefore, we can invoke Berge's maximum theorem (see Theorem~\ref{thm:Berge} in Appendix~\ref{app:Berge}) to conclude that $x\mapsto{\max}_{y\in\calY\bigcap\dom J} \Phi(x,y) - J(y)$ is continuous on $\bbX$. Since $g$ is closed and convex, it is continuous on $\dom g\cap \calX$, which is compact since $\dom g$ is closed and $\calX$ is compact (cf.\ Assumption~\ref{assump:bounded}\ref{assump:bounded_sets}). Therefore, we conclude that $\barS$ is continuous on $\dom g\cap\calX$. Similarly, we can show that $\uS$ is continuous on $\dom J\cap\calY$. Hence, there exists a positive constant $\Gamma<+\infty$ such that 
\begin{equation}
{\sup}_{x\in\dom g\cap \calX}\;{\sup}_{y\in\dom J\cap \calY}G(x,y) = {\sup}_{x\in\dom g\cap \calX} \;\barS(x) - {\inf}_{y\in\dom J\cap \calY} \;\uS(y)\le \Gamma.\label{eq:bound_sup_gap} %<+\infty. 
\end{equation}
Based on this observation, we can easily derive the following result. 

\begin{theorem}\label{thm:sc_expectation}
Assume that $\calB(0,1)\subseteq \dom h_\calX$. Also, let Assumptions~\ref{assump:bounded}\ref{assump:bounded_sets},~\ref{assump:noise}\ref{assump:unbiased},~\ref{assump:noise}\ref{assump:sub-Gaussian}, and~\ref{assump:regularity_g_J} hold. 
In Algorithm~\ref{algo:restart_stoc}, for any $x_0\in\calX^o$ and desired accuracy $\varepsilon\in(0,\mu U^2/2]$, if we choose $K =\big\lceil\log_2\big(\mu U^2/(2\varepsilon)\big)\big\rceil+1$ and $\nu=\min\{\varepsilon/(2\Gamma),1\}$, then $\bbE[G(x_{K+1},y_{K+1})] \le \varepsilon$. Moreover, the oracle complexity of Algorithm~\ref{algo:restart_stoc} is
\begin{equation}
O\left(\bigg(\sqrt{\frac{L}{\mu}}+\frac{L_{xx}}{\mu}\bigg)\log\left(\frac{1}{\varepsilon}\right)+ \frac{L_{yx}}{\sqrt{\mu\varepsilon}} + \frac{L_{yy}}{\varepsilon}+ \left(\frac{(\sigma_{x,f}+\sigma_{x,\Phi})^2}{\mu\varepsilon}+\frac{\sigma_{y,\Phi}^2}{\varepsilon^2}\right)\log\left(\frac{1}{\varepsilon}\right)\right).\label{eq:complexity_sc_exp}
\end{equation}
\end{theorem}

\proof{Proof.}
Define the event $\calG_{K,\varepsilon}\defeq \{G(x_{K+1},y_{K+1})\le \varepsilon/2\}$ and denote its complement as $\calG_{K,\varepsilon}^\rmc$. From Theorem~\ref{thm:restart_stoc}, we know that by the choice of $K$, $\Pr\{\calG_{K,\varepsilon}\}\ge 1-\nu$, for any $\nu\in(0,1]$. Therefore, 
\begin{align}
\bbE[G(x_{K+1},y_{K+1})] &=  \bbE[G(x_{K+1},y_{K+1})\bbI_{\calG_{K,\varepsilon}}] + \bbE[G(x_{K+1},y_{K+1})\bbI_{\calG^\rmc_{K,\varepsilon}}]\nn\\
&\lea (\varepsilon/2)\Pr\{\calG_{K,\varepsilon}\} + \Gamma\Pr\{\calG^\rmc_{K,\varepsilon}\}\nn\\
&\leb \varepsilon/2 + \Gamma\varepsilon/(2\Gamma)\le \varepsilon, 
\end{align}
where (a) follows from~\eqref{eq:bound_sup_gap} since $x_{K+1}\in\dom g\cap\calX$ and $y_{K+1}\in\dom J\cap\calY$, and (b) follows from the choice of $\nu$ (which implies that $\nu\le \varepsilon/(2\Gamma)$). To derive the oracle complexity in~\eqref{eq:complexity_sc_exp}, we simply substitute $\epsilon=\varepsilon/2$ and $\nu=\varepsilon/(2\Gamma)$ into~\eqref{eq:complexity_sc_hp}, and note that $O(\log(\log(1/\varepsilon)/\varepsilon)) = O(\log(1/\varepsilon))$. 
\Halmos
\endproof

%\begin{remark}
%Similar to~\eqref{eq:complexity_cvx_f}, 
%\end{remark}

%\newpage

%\subsection{Restart schemes with unknown $\mu$.}

\subsection{Discussions.}
We conclude this section by discussing some technical issues. 

%\begin{enumerate}[labelindent=0pt]
\subsubsection{Restart schemes.} \label{sec:restart}
Note that our restart scheme in Algorithm~\ref{algo:restart_stoc}, which is developed for stochastic SPPs, greatly differs from the usual approaches for convex minimization problems (see e.g.,~\citet{Nest_13},~\citet{Ghad_13b} and in particular,~\citet{Renegar_18} for a comprehensive review).
To be specific, consider $\min_{x\in\calX} f(x)$, where $f$ is $L$-smooth and $\mu$-strongly convex on $\calX$ (cf.\ Section~\ref{sec:intro}) and denote the minimizer by $x^*$. Also, assume that we have access to stochastic gradients of $f$ that satisfy Assumption~\ref{assump:noise}. Indeed, most of the existing restart schemes developed for this problem use the (primal) sub-optimality $f(x)-f(x^*)$, either in expectation or with high probability, as the reduction criterion. To illustrate, consider the expected sub-optimality. These restart schemes require the subroutine to satisfy the following: for any starting point $\barx^1\in\calX$  and $\epsilon,\delta>0$, there exists $T\in\bbN$ such that %the iterates produced by this subroutine, denoted by $\{\}$ %to satisfy that  such that 
\begin{align}
\bbE[\normt{\barx^1-x^*}^2]\le \delta \quad \Longrightarrow\quad \bbE[f(\barx^{T})-f(x^*)]\le \epsilon, \label{eq:recur_sc_min}
\end{align}
where $\barx^{T}$ denotes $T$-th iterate produced by the subroutine. 
By the strong convexity of $f$, we have 
\begin{equation}
\bbE[\normt{\barx^1-x^*}^2]\le ({2}/{\mu})\bbE[f(\barx^1)-f(x^*)],\label{eq:bound_dist}
\end{equation}
and therefore, we can establish an recursion between $\bbE[f(\barx^{T})-f(x^*)]$ and $\bbE[f(\barx^1)-f(x^*)]$. (Indeed, by properly choosing $T$, this recursion becomes a contraction.) However, somewhat interestingly, this approach cannot be straightforwardly extended to solve stochastic SPPs,  %we cannot use the seemingly analogous 
which would involve using the expected duality gap as the reduction criterion. %This is due to two reasons. 
There are at least two reasons behind this. First, by the definition of the duality gap in~\eqref{eq:duality_gap}, we do not have the analog of~\eqref{eq:recur_sc_min}, i.e., for any starting point $(\barx^1,\bary^1)\in\calX^o\times\calY^o$ and $\epsilon,\delta>0$, there exists $T'\in\bbN$ such that
\begin{align}
\max\big\{\bbE[\normt{\barx^1-x^*}^2],\bbE[\normt{\bary^1-y^*}^2]\big\}\le \delta \quad \Longrightarrow\quad \bbE[G(\barx^{T'},\bary^{T'})]\le \epsilon, \label{eq:recur_sc_spp}
\end{align}
where recall that $(x^*,y^*)$ denotes a saddle-point of~\eqref{eq:main}. 
Second, even if~\eqref{eq:recur_sc_spp} holds, since we do not assume strong convexity on the dual function $\uS$, there is no analog of~\eqref{eq:bound_dist} on the dual variable $\bary^1$. Therefore, we still cannot establish a recursion between $\bbE[G(\barx^{T'},\bary^{T'})]$ and $\bbE[G(\barx^{1},\bary^{1})]$. As such, more sophisticated restart schemes need to be developed. 
%However, observe that the subroutine used in our restart scheme, i.e., Algorithm~\ref{algo:convex_f_scaled}, does not satisfy this requirement. From the definition of the duality gap in~\eqref{eq:duality_gap}, we see that there is no analog of the ``reference point'' $x^*$, and the (expected) duality gap is bounded by $\sup_{x\in\calX}D_{h_\calX}(x,x_c)$ and $\sup_{y\in\calY}D_{h_\calY}(y,y_c)$ (as opposed to $\normt{x^*-x_c}^2$ and $\normt{y^*-y_c}^2$). 

Inspired by~\citet{Juditsky_12b}, in the deterministic setting, we use the quantity $R$ as the reduction criterion (cf.\ Section~\ref{sec:det_restart}). The quantity $R$ not only serves as an upper bound of the distance between the primal iterate $\barx^T$ to $x^*$, but also of the duality gap $G(\barx^T,\bary^T)$ (cf.\ Proposition~\ref{prop:conv_scaled}). Thus by sufficiently reducing $R$ (possibly over multiple stages), we can reduce the duality gap below any desired accuracy $\epsilon>0$. 
To achieve this, we modify Algorithm~\ref{algo:convex_f} into Algorithm~\ref{algo:convex_f_scaled}, whose convergence rate depends on $R$ and hence can be used as a subroutine in Algorithm~\ref{algo:restart_det}. % (i.e., Algorithm~\ref{algo:convex_f_scaled}). %, so as to incorporating $R$ in the . 
Although our deterministic restart scheme is in the same spirit as the one in~\citep{Juditsky_12b}, they are different in two important aspects.  First, since our subroutine (i.e., Algorithm~\ref{algo:convex_f_scaled}) is based on PDHG, but the subroutine used in~\citep{Juditsky_12b} is based on Mirror-Prox, our approach for developing the subroutine (cf.\ Section~\ref{sec:algo_modified}) is very different from that in~\citep{Juditsky_12b}.  As developing subroutines plays a crucial role in developing restart schemes,  this  difference clearly distinguishes our scheme from that in~\citep{Juditsky_12b}. 
 Second, in contrast to the approach in~\citep{Juditsky_12b}, by using different DGFs  (cf.~\eqref{eq:new_DGFs}) and properly choosing stepsizes in our subroutine (cf.\ Propositions~\ref{prop:conv_scaled} and~\ref{prop:conv_scaled_stoc}), our approach do not require defining new norms on $\bbX$ or $\bbY$. Therefore, it is simpler and easier to understand. 

In the stochastic setting, the problem becomes even more challenging. As discussed above (see also Section~\ref{sec:intuition}), the restart scheme based on the expected duality gap fails to work. Instead, we develop the restart scheme (i.e., Algorithm~\ref{algo:restart_stoc}) by leveraging the framework of finite-state Markov chain, which enables us to analyze the oracle complexity to obtain an $\epsilon$-duality gap with high probability. In addition, if the duality gap function $G(\cdot,\cdot)$ is bounded on $\calX\times\calY$, we can also derive the oracle complexity to achieve an  $\epsilon$-expected duality gap (cf. Section~\ref{sec:conv_expectation}). 

Shortly after the initial appearance of our paper on arXiv~\citep{Zhao_19a},~\citet{Yan_19} also developed a stochastic restart scheme for a different class of SPP from the one studied in this work. Specifically, they assume $f\equiv 0$, $\bbX$ and $\bbY$ are Hilbertian, and $\Phi(x,\cdot)$ is  a linear function for any $x\in\calX$. Also, instead of assuming $g$ and $J$ have easily computable BPPs on $\calX$ and $\calY$, they assume a primal-dual stochastic sub-gradient oracle exists for $S(\cdot,\cdot)$. Moreover, such a stochastic sub-gradient is {\em uniformly} bounded on $\calX\times\calY$ {\em almost surely}. As another important difference, they are interested in obtaining an $\epsilon$-primal sub-optimality gap, in expectation or with high probability, rather than an $\epsilon$-duality gap. Due to these differences, their restart scheme, as well as the subroutines therein, is different from ours. Interestingly, they consider the case where the primal function $\barS$ (cf.~\eqref{eq:primal_dual_function}) satisfies the local error bound condition~\citep[Definition~1]{Yan_19}, which might be worth exploring in our context. However, we defer this to future work.  
%This condition 
%Compared to the class of SPP studied in this work, i.e.,~\ref{eq:main}, there are 
%structure hidden in the restart scheme, and 

\subsubsection{Sub-Gaussian gradient noises.}\label{sec:light_tailed}
In our stochastic restart scheme, we assume that the gradient noises $\{\delta_{y,\Phi}^{t}\}_{t\in\bbN}$, $\{\delta_{x,\Phi}^{t}\}_{t\in\bbN}$ and $\{\delta_{x,f}^{t}\}_{t\in\bbN}$ not only have uniformly bounded (conditional) second moments, but also follow sub-Gaussian distributions. The reason is that sub-Gaussianity leads to concentration, so that we have large-deviation-type results at each stage (cf.\ Proposition~\ref{prop:conv_scaled_stoc}). %Moreover, we still have such type of results over multiple stages (cf.\ Theorem~\ref{thm:restart_stoc}). 
However, these results fail to hold if the noises are ``heavy-tailed''.  %we do not have % (e.g., not sub-Gaussian), 
Although we can still apply Markov's inequality at each stage of Algorithm~\ref{algo:restart_stoc} (in lieu of concentration inequalities), 
%following the same approach, 
the oracle complexity in~\eqref{eq:complexity_sc_hp} will have a significantly worse dependence on $\nu$, i.e., $O(1/\nu)$. Since in the proof of Theorem~\ref{thm:sc_expectation}, we need to choose $\nu=\Theta(\epsilon)$, this worse dependence on $\nu$ %$O(1/\nu)$ complexity 
leads to a worse dependence on $\epsilon$ in the complexity in~\eqref{eq:complexity_sc_exp}. For example,  to obtain an $\epsilon$-expected duality gap, the complexities corresponding to $(\sigma_{x,f}+\sigma_{x,\Phi})^2$ and $\sigma_{y,\Phi}^2$  degrade to $O((\sigma_{x,f}+\sigma_{x,\Phi})^2/(\mu\epsilon^2))$ and $O({\sigma_{y,\Phi}^2}/{\epsilon^3})$, respectively. Indeed, these complexities are worse than those in~\eqref{eq:complexity_cvx_f}, which corresponds to the case where $\mu=0$. 
%, since in the proof of Theorem~\ref{thm:sc_expectation} we need to choose $\nu=\epsilon/(2\Gamma)$. 
As future work, it would be interesting to develop new methodology that can effectively deal with the ``heavy-tailed'' noises. 
%Consequently, the oracle complexity 

\subsubsection{Boundedness of $\calX$.}% and $\Omega_{h_\calY}$.}
\label{sec:bounded_discussion}
 We notice that the boundedness assumption on $\calX$ %and $\Omega_{h_\calY}$, i.e., Assumption~\ref{assump:bounded}\ref{assump:bounded_sets}, 
 is required by all the complexity results regrading both of the deterministic and stochastic restart schemes (i.e., Algorithms~\ref{algo:restart_det} and~\ref{algo:restart_stoc}). 
 Since $f$ is $\mu$-strongly convex ($\mu>0$), in the following, we provide two ways to relax the assumption $D_\calX<+\infty$ in Algorithm~\ref{algo:restart_det}, by assuming that %$\calY$ is bounded. %, i.e., 
 $D_\calY<+\infty$. Note that this setting (i.e., $\mu>0$, $\calX$ is unbounded, and $\calY$ is bounded) occurs in many applications. For example, consider the three-composite optimization problem (see e.g.,~\citet{Zhao_19}):
 \begin{align}
 \min_{x\in\calX} \; f(x) + g(x) + h^*(\rvA x), \label{eq:three_composite}
 \end{align}
 where recall from Section~\ref{sec:related_work} that $\rvA:\bbX\to\bbY^*$ is a linear operator, and $h^*:\bbY^*\to\bbR$ is $M_h$-Lipschitz on $\bbY^*$ ($M_h<+\infty$). Instances of this problem include graph-guided fused Lasso (with full-column-rank data matrix)~\citep{Kim_09} and elastic-net regularization~\citep{Zou_05}. Note that the problem in~\eqref{eq:three_composite} can be equivalently formulated into the SPP in~\eqref{eq:main}, with $\Phi(x,y)=\lranglet{\rvA x}{y}$, and $\calY = \dom h$. In addition, since $M_h<+\infty$, we see that $\calY$ is indeed bounded. 
 % (which implies the boundedness of $\calY$). 
%  Nevertheless, we discuss below how it can be relaxed in some specific scenarios. 

% \subsubsection{Relaxing the boundedness assumption of $\calX$.}\label{sec:Boundedness_of_X}
%The first approach is  specific to Algorithm~\ref{algo:restart_det}. %the deterministic restart scheme (i.e.,). 
To start with, let us observe that the assumption $D_\calX<+\infty$ is only used in the input condition $U\ge D_\calX$, % (in both Algorithms~\ref{algo:restart_det} and~\ref{algo:restart_stoc}), 
as in general we can only use $D_\calX$ to bound the initial distance $\normt{x^1-x^*}$, i.e., $D_\calX\ge \normt{x^1-x^*}$. However, since $\barS$ is $\mu$-strongly convex on $\calX$, if we can find $\ell\in\bbR$ such that $\ell\le \min_{x\in\calX}\barS(x) = \barS(x^*)$, then $\normt{x^1-x^*}^2 \le ({2}/{\mu}) (\barS(x^1) - \barS(x^*)) \le ({2}/{\mu}) (\barS(x^1) - \ell)$ and consequently, 
%Note that the boundedness assumption of $\calX$  However, we can remove this assumption if we have additional information, e.g.,  To see this, note that , and therefore
\begin{equation}
%\;\Longrightarrow\; 
\normt{x^1-x^*}\le \sqrt{{2}(\barS(x^1) - \ell)/{\mu} }. \label{eq:bound_init_dist}
\end{equation}
Thus, the assumption $D_\calX<+\infty$ is no longer needed. 
%In this case, it suffices to let $U\ge \sqrt{{2}/{\mu}} (\barS(x^1) - \ell)^{1/2}$. 
We now provide a general approach to derive the lower bound $\ell$. Let us focus on the deterministic restart scheme (i.e., Algorithm~\ref{algo:restart_det}) first. Fix any $y\in\calY$. We then find an $\eta$-approximate optimal solution (where $\eta>0$) of the problem $\min_{x\in\calX} S(x,y)$, denoted by $x^*_\eta(y)$, such that 
\begin{equation}
S(x^*_\eta(y),y) -{\min}_{x\in\calX} S(x,y)\le \eta. \label{eq:bdd_eta}
\end{equation} 
Since $f$ is $\mu$-strongly convex and  $L$-smooth on $\calX$ and $\Phi(\cdot,y)$ is $L_{xx}$-smooth on $\calX$, from Lemma~\ref{lem:grad_mapping},  to find $x^*_\eta(y)$, it suffices to find $\barx\in\calX$ such that its gradient mappings, i.e., $G_\lambda(\barx) \defeq (\barx-\barx^+)/\lambda$ and $\barG_\lambda(\barx)\defeq \big(\nabla h_\calX(\barx)-\nabla h_\calX(\barx^+)\big)/\lambda$, satisfy that 
\begin{align}
\normt{G_\lambda(\barx)}^2 + \normt{\barG_\lambda(\barx)}_*^2\le \mu\eta,\label{eq:grad_map_norm_uB}
\end{align}
where $0<\lambda\le (L+L_{xx})^{-1}$ and %(recall that ) and 
\begin{align}
\barx^+ \defeq {\argmin}_{x\in\calX} \;\lranglet{\nabla f(\barx) + \nabla_x \Phi(\barx,y)}{x} + g(x) + \lambda^{-1}D_{h_\calX}(x,\barx), \label{eq:def_x^+}
\end{align}
and take $x^*_\eta(y) = \barx^+$. 
%For details, we refer readers to Appendix~\ref{app:grad_map}. 
%check~\eqref{eq:bdd_eta}, 
Clearly, we can set the lower bound of $\barS(x^*)$, i.e., $\ell=S(x^*_\eta(y),y)-\eta$, since by the definition of $\barS$ in~\eqref{eq:primal_dual_function},  
\begin{equation}
{\min}_{x\in\calX} \barS(x)\ge {\min}_{x\in\calX} S(x,y)\ge S(x^*_\eta(y),y) - \eta. 
\end{equation}
A potential issue of the approach above is that $\barS(x^1)$ in~\eqref{eq:bound_init_dist} may not be evaluated exactly. In this case,  we need to find an upper bound $\upsilon\ge \barS(x^1)$. The way to find $\upsilon$ is similar to that for finding $\ell$. Specifically, we can find an $\eta$-approximate optimal solution of the problem $\max_{y\in\calY}\,\Phi(x^1,y) - J(y)-(\eta/\Upsilon_{h_\calY})h_\calY(y)$, denoted by $y^*_\eta(x^1)$, where $\Upsilon_{h_\calY} = \max_{y\in\calY}\abst{h_\calY(y)}<+\infty$. 
 As a result,  % to $\eta$%
\begin{align}
\barS(x^1)\le \upsilon\defeq f(x^1)+g(x^1)+\Phi(x^1,y^*_\eta(x^1)) - J(y^*_\eta(x^1))+2\eta. 
\end{align} 
To find $y^*_\eta(x^1)$, following Lemma~\ref{lem:grad_mapping}, it suffices to find $\bary\in\calY$ such that
\begin{align}
\norm{(\bary-\bary^+)/\lambda}^2 + \norm{\big(\nabla h_\calY(\bary)-\nabla h_\calY(\bary^+)\big)/\lambda}_*^2\le \eta^2/\Upsilon_{h_\calY},\label{eq:grad_map_norm_uB_Y}
\end{align}
where $0<\lambda\le L_{yy}^{-1}$ (recall that $\Phi(x,\cdot)$ is $L_{yy}$-smooth on $\calY$) and 
\begin{align}
\bary^+ \defeq {\argmin}_{y\in\calY} \;-\lranglet{\nabla_y \Phi(x^1,\bary)}{y} + J(y) +(\eta/\Upsilon_{h_\calY})h_\calY(y) + \lambda^{-1}D_{h_\calY}(y,\bary), 
\end{align}
and set $y^*_\eta(x^1) = \bary^+$. 
% For both schemes, as noted , the assumption $D_\calX<+\infty$ can be removed as long as
% The first method is illustrated 
%
% in Section~\ref{sec:Boundedness_of_X}, for the deterministic restart scheme. %Namely, we find a lower bound $\ell\le \min_{x\in\calX}\barS(x)$, and choose  and as illustrated in Remark~\ref{sec:Boundedness_of_X}, there exists a generic approach to construct the lower bound $\ell$.  %If $\barS$ is the loss function in many statistical learning problems (e.g., sparse group LASSO~\citep{Simon_13}), we can simply set $\ell=0$. 
 %In the generic way 
 
%Now, to apply the approach above to the stochastic scheme (i.e., Algorithm~\ref{algo:restart_stoc}),   the difficulty lies in the fact that the function values and gradients of $f$ and $\Phi(\cdot,\cdot)$ %in~\eqref{eq:bound_init_dist} 
%may not be evaluated exactly, due to high-dimensional integrals over $\Xi$ and $\calZ$ (cf.~\eqref{eq:expectation}). 
%However,  we can find an upper bound of $\barS(x^1)$, denoted by $\upsilon$, with high probability (say, at least $1-\delta/4$), by using sample average approximation or Monte Carlo methods. Then~\eqref{eq:bound_init_dist} becomes 
%\begin{equation}
%\bbP\left(\normt{x^1-x^*}^2 \le {2}(\upsilon - \ell)/{\mu} \right)\ge 1-\delta/4. 
%\end{equation}
%As a result, if we choose $U\ge \sqrt{{2}(\upsilon - \ell)/{\mu}} $, then $\bbP\left(\normt{x^1-x^*}\le U\right)\ge 1-\delta/4$. Then we can condition on the event $\{\normt{x^1-x^*}\le U\}$ and perform the analysis in Section~\ref{sec:stoc_restart}. 
 
The second approach aims to write the SPP in~\eqref{eq:main} into an equivalent SPP that has a bounded primal constraint set, and so %is more general than the first one, 
it is independent of Algorithm~\ref{algo:restart_det}. 
For any $y\in\calY$, let us define 
\begin{equation}
x^*(y) \defeq {\argmin}_{x\in\calX} \big[\psi^\rmP(x,y)\defeq f(x) + g(x) + \Phi(x,y)\big]. \label{eq:def_x^*}
\end{equation}
From Lemma~\ref{lem:lips_x^*} (see Appendix~\ref{app:Lips}), we know that $x^*(\cdot)$ is $(L_{yx}/\mu)$-Lipschitz on $\calY$. Thus, we have %for any $y,y'\in\calY$, we have 
\begin{align}
\normt{x^*(y')} \le \normt{x^*(y)} + (L_{yx}/\mu)\normt{y-y'}\le \normt{x^*(y)} + L_{yx}D_\calY/\mu,\quad \forall\;y,y'\in\calY. \label{eq:bound_x^*(y')}
\end{align}
To bound $\normt{x^*(y)}$, let us find $x^*_\eta(y)$ in the same way as in the first approach, and by definition, $x^*_\eta(y)$ satisfies~\eqref{eq:bdd_eta}. %(Note that )
 As a result, 
\begin{equation}
\normt{x^*(y) - x_\eta^*(y)}\le \sqrt{2\eta/\mu}\quad\Longrightarrow\quad \normt{x^*(y)}\le \normt{x_\eta^*(y)} +  \sqrt{2\eta/\mu}. \label{eq:bound_x^*(y)}
\end{equation} 
If we combine~\eqref{eq:bound_x^*(y')} and~\eqref{eq:bound_x^*(y)}, then  we have %for any $y'\in\calY$, 
\begin{equation}
\normt{x^*(y')}\le \varkappa\defeq \normt{x_\eta^*(y)} +  \sqrt{2\eta/\mu} + L_{yx}D_\calY/\mu,\quad \forall\,y'\in\calY. \label{eq:bound_x^*(y')_final}
\end{equation}
If we define $\tilcalX\defeq \{x^*(y):y\in\calY\}\ne\emptyset$, then from~\eqref{eq:bound_x^*(y')_final}, we see that $\tilcalX\subseteq\calX\cap\calB(0,\varkappa)$. As a result, the SPP in~\eqref{eq:main} can be equivalently written as 
\begin{equation}
{\min}_{x\in\calX\cap\calB(0,\varkappa)}\;{\max}_{y\in\calY}\; S(x,y).\label{eq:main_alt}
\end{equation}
Note that in this new SPP, the primal constraint set $\calX\cap\calB(0,\varkappa)$ is nonempty, convex and compact. 

Although the two approaches above work well in Algorithm~\ref{algo:restart_det}, there is one major difficulty that prevents us from extending either of them to the stochastic case. Recall from Section~\ref{sec:stoc_oracle} that we only have access to the stochastic gradients of $f$ and $\Phi(\cdot,\cdot)$. Therefore, in the first approach, the proximal point $\barx^+$ defined in~\eqref{eq:def_x^+} becomes stochastic, hence the gradient mappings $G_\lambda(\barx)$ and $\barG_\lambda(\barx)$,  as well as $x^*_\eta(x)$ (since $x^*_\eta(x)=\barx^+$). One of the most natural ways to handle such stochasticity is to let the event in~\eqref{eq:grad_map_norm_uB} hold with high probability (say, no less than $1-\delta/4$), and perform the rest of the analysis by conditioning on this event. However, note that we only have information about the distributions of $\hat{\nabla} f(x)$ and $\hat{\nabla} \Phi(x,y)$. In addition, the mappings from $\hat{\nabla} f(x)+\hat{\nabla} \Phi(x,y)$ to $\barx^+$, as well as from $\barx^+$ to  $\normt{G_\lambda(\barx)}^2 + \normt{\barG_\lambda(\barx)}_*^2$, are rather complicated. Therefore, it is unclear how to lower bound the probability of the event in~\eqref{eq:grad_map_norm_uB}. Certainly, there may exist more sophisticated approaches to relax the assumption $D_\calX<+\infty$ in the stochastic case, and we leave the investigation of these approaches to  future work.   
%However, we remark that %the aforementioned difficulties do not exist 
As a final remark, if both of the sample spaces ${\Xi}$ and $\calZ$ have finite sizes (cf.~\eqref{eq:finite-sum}), %since in that case 
we will still have access to the exact gradients $f$ and $\Phi(\cdot,\cdot)$. 
As a result, the aforementioned two approaches continue to work in this case. 
%  %As a final remark, note that 
%in checking the gradient mappings fo

\subsubsection{Unknown problem parameters.} %$\Omega_{h_\calX}$ and $\Omega_{h_\calY}$.}
\label{sec:estimation_diameters}

Note that our restart schemes (i.e., Algorithms~\ref{algo:restart_det} and~\ref{algo:restart_stoc}) require the knowledge of several problem parameters. These include the strong convexity parameter $\mu$, the smoothness parameters  $L$, $L_{xx}$, $L_{xy}$ and $L_{yy}$,  the diameters $\calD_\calX$ and $\Omega_{h_\calY}$, and the variance parameters $\sigma^2_{x,f}$, $\sigma^2_{x,\Phi}$ and $\sigma^2_{y,\Phi}$. In the following, we will discuss possible ways to estimate some of these parameters, by assuming the knowledge of the rest.

To begin with, let us focus on unknown $\calD_\calX$, or more precisely, unknown $\normt{x^1-x^*}$. Note that in Section~\ref{sec:bounded_discussion}, we already provide two ways to estimate $\normt{x^1-x^*}$, by assuming that $D_\calY<+\infty$. These methods apply to both the deterministic case and the stochastic case with finite sample spaces $\Xi$ and $\calZ$ (i.e., the finite-sum case). In the future, it would  be interesting to develop a procedure that starts with an initial guess of $\normt{x^1-x^*}$, and based on the duality gap (and other information),  to search for an upper bound on $\normt{x^1-x^*}$. %We leave this to future work. 

Next, we discuss how to handle unknown $\mu$. Note that even in the deterministic case, it is non-trivial to develop a modified version of our restart scheme (i.e., Algorithm~\ref{algo:restart_det}), %a restart scheme for~\eqref{eq:main}, 
which automatically adapts to $\mu$ and has the same oracle complexities up as Algorithm~\ref{algo:restart_det} (up to a log factor). 
%to logarithmic factor is a task, .  
%restart scheme f without knowing $\mu$ 
In the context of strongly convex (smooth and non-smooth) minimization, many adaptive  restart schemes have been proposed to avoid the knowledge of $\mu$ (see e.g.,~\citet{Nest_13},~\citet{Juditsky_14},~\citet{Lin_15b} and~\citet{Renegar_18}). Among them, the restart scheme in~\citep{Juditsky_14} is most relevant to ours, since it also uses the quantity $R$ as the reduction criterion (as opposed to the primal sub-optimality gap; cf.\ Section~\ref{sec:restart}). %based on the domain-shrinkage technique. 
The basic idea is  to notice that $\mu$ resides in a bounded interval $(0,B]$, and partition this interval as $(0,\mu_1]\cup(\mu_1,\mu_2]\cup\ldots\cup(\mu_{K-1},B]$, where $K$ denotes the total number of stages (cf.~Algorithm~\ref{algo:restart_det}), and 
\begin{align*}
\mu_k = \mu_{k+1}/2, \quad \forall\,k\in\calK\defeq \{1,\ldots,K-1\}.
\end{align*}
Furthermore, given the total number of iterations $T$,  number of iterations in each stage is the same and equal to $\floor{T/K}$ (i.e., no longer exponentially increasing as in Algorithm~\ref{algo:restart_det}). By choosing $K$ properly (and independent of $\mu$), they obtain a convergence rate $O(\log T/T)$ of the {\em primal sub-optimality gap}. However, a close inspection reveals that their approach cannot be easily extended to Algorithm~\ref{algo:restart_det}. Let $k^*\in\calK$ satisfy $\mu\in(\mu_{k^*},\mu_{k^*+1}]$. One of the major difficulties is that in their scheme, for $k\ge k^*+1$, we no longer have the condition in~\eqref{eq:suff_close}, which is crucial to establish $\hatG_{R_k}(\barx_{k+1},y_{k+1})=G(\barx_{k+1},y_{k+1})$. Therefore, in this case, we cannot derive the convergence rate of the {\em duality gap}. Due to this, we leave developing $\mu$-free versions of Algorithms~\ref{algo:restart_det} and~\ref{algo:restart_stoc} as open questions. However, since these difficulties are mainly related to theoretical (complexity) analysis, in practice, we can still use the aforementioned restart schemes as heuristics. 
%the optimal solution $x^*$ is no longer in the 

Finally, we discuss how to handle the unknown smoothness parameters, i.e., $L$, $L_{xx}$, $L_{xy}$ and $L_{yy}$. For simplicity, let us focus on unknown $L$ (as the other three parameters can be handled similarly). Indeed, in the deterministic case,~\citet[Section~3]{Nest_13} has proposed a backtracking method for estimating $L$. %Specifically, 
This method starts from a lower bound on $L$, and increase it by a constant factor until certain termination criterion is satisfied. This method can be straightforwardly integrated into the subroutine of Algorithm~\ref{algo:restart_det}, i.e., Algorithm~\ref{algo:convex_f_scaled}. However, in the stochastic case, this backtracking method fails to work properly.  Specifically, since the termination criterion now involves random variables, we need it to hold in expectation or with high probability. However, both of which are hard to check, since we have little knowledge about the distribution of random variables involved. For this reason, we leave estimating smoothness parameters in the stochastic case to future work.  % (e.g., $\mu$). 

\section{Application and numerical experiments.}\label{sec:experiment}
We consider the problem of two-player zero-sum game with non-linear payoff in~\citet{Hien_17} (see also~\citet{Chen_17}). The problem has the following form:
\begin{equation}
\label{eq:game2player}
\min_{x\in\calX} \max_{y\in\calY} \;\;\frac{\varpi}{2}\lranglet{x}{Qx} + \sum_{i=1}^{n}\log\left(1+\frac{y_{i}}{c_{i}+x_{i}}\right),
\end{equation}
where %Define the following sets:
\begin{equation}
\calX\defeq \{x\in\bbR^n:x\geq 0,\textstyle\sum_{i=1}^n x_i=N\} \quad \mbox{and}\quad  \calY\defeq \{y\in\bbR^n: y\geq 0, \textstyle\sum_{i=1}^n y_i=P\}.  \label{eq:game_constraint_sets}
\end{equation}
There are several problem constants in~\eqref{eq:game2player} and~\eqref{eq:game_constraint_sets}, including $\varpi> 0$, $c_i>0$, for all $i\in [n]$, $Q\in\bbS_{++}^n$ (i.e., $Q\in\bbR^{n\times n}$ is a symmetric and positive definite matrix), $N>0$, and $P>0$. %Let us consider the following problem:   

Indeed, the problem in~\eqref{eq:game2player} corresponds to the water-filling problem in information theory (see e.g.,~\citet{Cover_06}). Specifically, let Alice and Bob be two players, both of whom face $n$  Gaussian communication channels, each with noise power $c_i$, for all $i\in[n]$. Alice has noise power $N$, and wishes to minimize the total channel capacity, by allocating noise power $x_i$ to channel $i$, for each $i\in[n]$. Bob has transmission power $P$, and wishes to maximize the total channel capacity, by allocating transmission power $y_i$ to channel $i$, for each $i\in[n]$. % and transmission power $P$, respectively. 
%Specifically,  and each channel $i$ has   Alice  whereas Bob aims to allocate his transmission power  to these $n$ channels to maximize the total channel capacity. 
%Let $x = (x_1,\ldots,x_n)$ and $y = (y_1,\ldots,y_n)$ denote the allocation plans of Alice and Bob, respectively. 
The total channel capacity is given by $\sum_{i=1}^{n}\log\left(1+{y_{i}}/({c_{i}+x_{i}})\right)$. %where $\{\beta_i\}_{i\in[n]}$ are normalization constants. 
In addition, for Alice, injecting noise power to the channels incurs cost $(1/{2})\lranglet{x}{Qx}$ (where $Q$ denotes the cost matrix), and $\varpi$ reflects the trade-off between minimizing her cost and minimizing the total channel capacity. In fact, by finding a saddle point of~\eqref{eq:game2player} (which always exists by Sion's minimax theorem~\citep{Sion_58}), we  find a Nash equilibrium of this two-player game. 
%the allocation plans $x$ and $y$ incur costs $({\varpi}/{2})\|x\|_2^2$ for Alice and $({\gamma}/{2})\langle Qy,y\rangle$ for Bob, respectively. 

%As such, Problem~\eqref{eq:game2player} reflects the zero-sum game where Alice (resp.\ Bob) wishes to minimize (resp.\ maximize) the total channel capacity, while accounting for her (resp.\ his) allocation cost. The saddle points of Problem~\eqref{eq:game2player} are precisely the Nash equilibria of this game. 

%\subsection{Experimental setup.}
Let us complete setting up the problem in~\eqref{eq:game2player}, by specifying the normed spaces $\bbX$ and $\bbY$. Indeed, if we set $\bbX = (\bbR^n,\norm{\cdot}_2)$ and $\bbY = (\bbR^n,\norm{\cdot}_1)$, then in~\eqref{eq:game2player}, by letting 
\begin{equation}
f(x) \defeq ({\varpi}/{2})\lranglet{x}{Qx} \quad \mbox{and}\quad \Phi(x,y) \defeq \textstyle\sum_{i=1}^{n}\log\left(1+{y_{i}}/({c_{i}+x_{i}})\right), \label{eq:f_Phi}
\end{equation}
the objective function  $S(x,y) \defeq f(x) + \Phi(x,y)$ satisfies all the assumptions stated in Section~\ref{sec:intro}, where the constants $L$, $\mu$, $L_{xx}$, $L_{xy}$ and $L_{yy}$ are given by 
\begin{align}
L = \varpi \lambda_{\max}(Q), \quad \mu = \varpi \lambda_{\min}(Q), \quad L_{xx} = c_{\min}^{-2} - (c_{\min} + 1)^{-2}, \quad L_{xy} = L_{yy} = c_{\min}^{-2}. \label{eq:constants_game}
\end{align}
In~\eqref{eq:constants_game}, $\lambda_{\max}(Q)>0$ and $\lambda_{\min}(Q)>0$ denote the largest and smallest eigenvalues of $Q$, respectively, and $c_{\min}\defeq \min_{i\in[n]} c_i>0$. (For the derivation of the constants $L_{xx}$, $L_{xy}$ and $L_{yy}$, see Appendix~\ref{app:deriv_constants}.)  In addition, under this setting, it is natural to choose the DGFs $h_\calX$ and $h_\calY$ to be
\begin{equation}
h_\calX(x) = (1/2)\norm{x}_2^2 \quad \mbox{and}\quad  h_\calY(y) = \textstyle\sum_{i=1}^n y_i\log y_i, \quad \forall\,(x,y)\in\calX\times\calY. 
\end{equation}

In what follows, we will the stochastic Mirror-Prox algorithm, i.e., Algorithm~1 in~\citet{Juditsky_11}, as the benchmark.  For ease of reference, we will refer to our stochastic restart scheme (i.e., Algorithm~\ref{algo:restart_stoc}) and Algorithm~1 in~\citep{Juditsky_11} as {\sf SRS} and {\sf SMP}, respectively. 

\subsection{Implementation details.}
We conducted four sets of experiments, corresponding to different values of $n\in\{1000,2000,3000,4000\}$. For simplicity, in all our experiments, we set $\varpi=0.1$ and $c_i=1$, for all $i\in[n]$, in~\eqref{eq:game2player}. In addition, in~\eqref{eq:game_constraint_sets}, we set $P=N=1$, so $\calX= \calY = \Delta_n$, i.e., the $(n-1)$-dimensional probability simplex. 
For each value of $n$, we generated the cost matrix  $Q$ as follows. First we generated a matrix $\barQ\in\bbR^{n\times n}$ such that its entries $\{\barQ_{ij}\}_{i,j=1}^n$ are \iid standard normal. (Therefore, with probability one, $\barQ$ has full rank.) Then we let $Q = \barQ^T \barQ$, and consequently, 
\begin{equation}
f(x) = (1/2)\lranglet{x}{Qx}  = (1/2)\normt{\barQ x}_2^2 = \textstyle (1/2)\sum_{i=1}^n (\barQ_i^T x)^2, \label{eq:f_finite_sum}
\end{equation}
where $\barQ_i$ denotes the $i$-th column of $\barQ^T$. From~\eqref{eq:f_finite_sum} and~\eqref{eq:f_Phi}, we observe that both $f$ and $\Phi$ have finite-sum forms. Therefore, to generate stochastic estimators of the gradients $\nabla f$, $\nabla_x \Phi(\cdot,\cdot)$ and $\nabla_y\Phi(\cdot,\cdot)$, we uniformly randomly drew a mini-batch $\calB\subseteq [n]$ without replacement, and set 
\begin{align}
\hat{\nabla} f(x)  \defeq \frac{n}{\abs{\calB}} \sum_{i\in\calB} (\barQ_i^Tx)R_i, \quad \hat{\nabla}_x \Phi(x,y) &\defeq \frac{n}{\abs{\calB}} \sum_{i\in\calB} -\frac{y_i}{(1+x_i)(1+x_i+y_i)} e_i,\\
\hat{\nabla}_y \Phi(x,y) &\defeq \frac{n}{\abs{\calB}} \sum_{i\in\calB} \frac{1}{1+x_i+y_i} e_i,  
\end{align}
where $e_i$ denotes the $i$-th standard coordinate vector. Based on these estimators, we can set the variance parameters in Assumption~\ref{assump:noise}\ref{assump:variance_bound} as 
\begin{align}
\sigma_{x,f}^2 &= \frac{n^2 (n-\abs{\calB})}{\abs{\calB}(n-1)}\left(n^{-1}\textstyle\sum_{i=1}^n \norm{\barQ_i}_2^2\norm{\barQ_i}_\infty^2 - n^{-3}\sigma_{\min}^2(Q)\right),\quad \sigma_{x,\Phi}^2=\sigma_{y,\Phi}^2 = \frac{n^2 (n-\abs{\calB})}{\abs{\calB}(n-1)}. 
\label{eq:var_game}
\end{align}
(For the derivation of these parameters, see Appendix~\ref{app:deriv_constants}.) For simplicity, we set $\abs{\calB} = \ceil{n/2}$ in all our experiments. 

For each value of $n$, after generating $Q\in\bbS_{++}^n$ and setting the problem parameters, we were able to run both algorithms {\sf SRS} and {\sf SMP}. Specifically, for both algorithms, we used the same initial primal-dual pair $(x_0,y_0)$, which was randomly chosen from $\Delta_n\times \Delta_n$. As stated in Remark~\ref{rmk:abs_cosntant}, we replaced all the absolute constants in the definition of $T_k$ in Algorithm~\ref{algo:restart_stoc} by a constant $c=0.7$. For either {\sf SRS} or {\sf SMP}, we stopped the algorithm once the duality gap fell below $\varepsilon = 0.001$, and recorded the running time accordingly. For each value of $n$, we repeated this process five times, and calculated the mean and standard deviation (std) of the running times for each algorithm. All the codes were implemented in Python 3.7.6 on a machine with 1.90 GHz CPU and 16 GB RAM, and the evaluation of the duality gap was done by the disciplined convex programming solvers in the CVXPY package~\citep{Diamond_16}. 
%we first chose a proportion $\gamma\in (0,1)$ and le

\subsection{Numerical results.}

The results are shown in Table~\ref{table:RegularizeProb}. From this table, we can observe that our algorithm, i.e., {\sf SRS}, outperforms {\sf SMP} for different values of $n$. Specifically, to reach a duality gap that is no more than $0.001$, {\sf SRS} uses less time than {\sf SMP} on average.   This is consistent with the superiority of the oracle complexity of {\sf SRS} over {\sf SMP} (cf.\ Tables~\ref{table:sc} and~\ref{table:cvx}). 

\setlength{\tabcolsep}{3pt}
\begin{table}[t!]\small%\scriptsize
\centering 
\caption{Comparison of the running times in seconds (mead $\pm$ std) obtained by {\sf SRS} and {\sf SMP} to reach $\varepsilon = 0.001$ over five runs.}
\label{table:RegularizeProb}
\begin{tabular}{|c|c|c|c|} %
 \hline
 $n$ &  {\sf SRS} & {\sf SMP} \\ 
 \hline 
1000 &   ${29.05 \pm 1.32}$ & $33.45\pm 0.95$ \\
2000 &    ${66.19 \pm 1.52}$ &$ 73.64  \pm 1.98$ \\
3000 &    $79.69\pm 2.33$ & $ 88.75 \pm 2.67$ \\
4000 &    ${100.63 \pm 2.85}$ & $ 112.80 \pm 2.69$ \\
\hline 
\end{tabular}
\end{table} 

\section{Proof of Proposition~\ref{prop:recursion_cvx}.}\label{sec:analysis}

Before proving Proposition~\ref{prop:recursion_cvx}, we first present the {\em Bregman proximal inequality} associated with the Bregman proximal projection in~\eqref{eq:Bregman_projection} and its corollary. The proof of this inequality can be found in many works, e.g.,~\citet[Lemma~2]{Ghad_12}.  

\begin{lemma}\label{lem:Bregman_proximal_inequality}
In~\eqref{eq:Bregman_projection}, for any $u\in\calU$, we have
\begin{align}
\varphi(u^+)-\varphi(u)%&\le \lranglet{u^*}{u-u^+} + \lambda^{-1}(D_{h_\calU}(u,u')-D_{h_\calU}(u,u^+)-D_{h_\calU}(u^+,u'))\\
&\le \lranglet{u^*}{u-u^+} + \lambda^{-1}(D_{h_\calU}(u,u')-D_{h_\calU}(u,u^+)) - (2\lambda)^{-1} \normt{u^+-u'}^2. \label{eq:Bregman_proximal_ineq}
\end{align}
%Fix $x\in\calX$. For any CCP function $\varphi:\bbU\to\barbbR$ and $u^*\in\bbU^*$, define
\end{lemma}

%Based on Lemma~\ref{lem:Bregman_proximal_inequality}, we can further derive another useful lemma. 

\begin{corollary}\label{cor:Bregman_proximal_inequality2}
In~\eqref{eq:Bregman_projection}, for any $u\in\calU$, we have
\begin{equation}
\varphi(u^+)-\varphi(u)\le \lranglet{u^*}{u-u'} + \lambda^{-1}(D_{h_\calU}(u,u')-D_{h_\calU}(u,u^+)) + (\lambda/2) \normt{u^*}_*^2. 
\end{equation}
%Fix $x\in\calX$. For any CCP function $\varphi:\bbU\to\barbbR$ and $u^*\in\bbU^*$, define
\end{corollary}
\proof{Proof.}
First, by Young's inequality, for any $u\in\bbU$ and $u^*\in\bbU^*$, % (see Section~\ref{sec:prelim}), 
we have
\begin{equation}
\abst{\lranglet{u^*}{u}} \le \normt{u^*}_*\normt{u}\le (\eta/2)\normt{u^*}_*^2 + (2\eta)^{-1} \normt{u}^2,\quad\forall\,\eta>0. \label{eq:Young}
\end{equation}
Therefore,
%\begin{align}
$\lranglet{u^*}{u'-u^+}\le (2\lambda)^{-1}\normt{u^+-u'}^2+(\lambda/2)\normt{u^*}_*^2$. \label{eq:Bregman_proximal_ineq2}
%\end{align}
Add up this inequality with~\eqref{eq:Bregman_proximal_ineq}, we then complete the proof. \Halmos
\endproof

%Based on Lemma~\ref{eq:Bregman_proximal_ineq2}, we can derive a lemma that will be used in the sequel. 
%\begin{lemma}
%For any $t\in\bbZ^+$ and $u_t^*\in\calU^*$, define $u^{t+1}\defeq \argmin_{u\in\calU} \lranglet{u_t^*}{u}+\lambda_t^{-1}D_{h_\calU}(u,u^t)$. Then 
%\end{lemma}

\proof{Proof of Proposition~\ref{prop:recursion_cvx}.}
For convenience, for any $t\in\bbZ^+$ and $(x,y)\in\calX\times\calY$, let us define 
\begin{align}
\tilG(\barx^{t},\bary^{t};x,y)\defeq &S(\barx^{t},y) - S(x,\bary^{t})\nn\\
=& (f(\barx^{t})-f(x)) + (g(\barx^{t}) - g(x)) + (\Phi(\barx^{t},y) - \Phi(x,\bary^{t})) + (J(\bary^{t})-J(y)). \label{eq:pseudo_gap}
\end{align}
The most crucial step in our proof is to establish the recursion between $\tilG(\barx^{t+1},\bary^{t+1};x,y)$ and $\tilG(\barx^{t},\bary^{t};x,y)$, which requires us to establish the recursion for each of the four terms in~\eqref{eq:pseudo_gap}. By the convexity of $g$ and $J$ and convexity-concavity of $\Phi$, %Jensen's inequality, 
we easily see that
\begin{align}
g(\barx^{t+1}) - g(x) &\le (1-\beta_t)(g(\barx^t)- g(x)) + \beta_t (g(x^{t+1})- g(x)),\label{eq:bound_g}\\
J(\bary^{t+1}) - J(y) &\le (1-\beta_t)(J(\bary^t)- J(y)) + \beta_t (J(y^{t+1})- J(y)), \label{eq:bound_J}\\
\Phi(\barx^{t+1},y) - \Phi(x,\bary^{t+1}) &\le [(1-\beta_t)\Phi(\barx^{t},y) + \beta_t\Phi(x^{t+1},y)] - [(1-\beta_t)\Phi(x,\bary^t) + \beta_t\Phi(x,y^{t+1})]\nn\\
&\le (1-\beta_t)(\Phi(\barx^{t},y) - \Phi(x,\bary^t)) + \beta_t (\Phi(x^{t+1},y) - \beta_t\Phi(x,y^{t+1})). \label{eq:bound_Phi}
\end{align}
%By using the $L$-smoothness and convexity of $f$, we also have
%Since $f$ is $L$-smooth, by the descent lemma (e.g.,~\citet{Bert_99}), we have
To connect $f(\barx^{t+1})$ and $f(\barx^{t})$, we have 
\begin{align}
f(\barx^{t+1}) - f(x)  \lea \;\;&f(\tilx^{t+1})  + \lranglet{\nabla f(\tilx^{t+1})}{\barx^{t+1}- \tilx^{t+1}} + (L/2)\normt{\barx^{t+1} - \tilx^{t+1}}^2 - f(x)\\
\leb \;\;& (1-\beta_t)(f(\tilx^{t+1})   + \lranglet{\nabla f(\tilx^{t+1})}{\barx^{t}- \tilx^{t+1}}- f(x))\nn\\ 
&+ \beta_t(f(\tilx^{t+1})   + \lranglet{\nabla f(\tilx^{t+1})}{x^{t+1}- \tilx^{t+1}}- f(x)) +(L/2)\normt{\barx^{t+1} - \tilx^{t+1}}^2\\
\lec \;\;& (1-\beta_t)(f(\barx^{t}) - f(x) )\nn\\
& + \beta_t(f(\tilx^{t+1}) + \lranglet{\nabla f(\tilx^{t+1})}{x^{t+1}- \tilx^{t+1}}- f(x)  ) +(L\beta_t^2/2)\normt{x^{t+1}-x^t}^2\\
\led \;\;&(1-\beta_t)(f(\barx^{t}) -f(x)) +\beta_t\lranglet{\nabla f(\tilx^{t+1})}{x^{t+1}- x} +(L\beta_t^2/2)\normt{x^{t+1}-x^t}^2, \label{eq:bound_f}
\end{align}
where in (a) we use the descent lemma (e.g.,~\citet{Bert_99}), resulted from the $L$-smoothness of $f$; in~(b) we use the step~\eqref{eq:ave_primal} in Algorithm~\ref{algo:convex_f}; in~(c) we use the convexity of $f$ and that $\barx^{t+1}-\tilx^{t+1} = \beta_t(x^{t+1}-x^t)$ (resulted from the steps~\eqref{eq:interp_primal} and~\eqref{eq:ave_primal}); in~(d) we again use  the convexity of $f$. 

Combining~\eqref{eq:bound_g},~\eqref{eq:bound_J},~\eqref{eq:bound_Phi} and~\eqref{eq:bound_f}, we have the following recursion
\begin{align}
&\tilG(\barx^{t+1},\bary^{t+1};x,y) \le \;(1-\beta_t)\tilG(\barx^{t},\bary^{t};x,y) + \beta_t(g(x^{t+1}) - g(x)) + \beta_t(J(y^{t+1}) - J(y))\nn\\
& + \beta_t(\Phi(x^{t+1},y)-\Phi(x,y^{t+1})) + \beta_t\lranglet{\nabla f(\tilx^{t+1})}{x^{t+1}- x} +(L\beta_t^2/2)\normt{x^{t+1}-x^t}^2.\label{eq:bound_pseudo_gap}
\end{align} 
We now can apply the Bregman proximal inequality in~\eqref{eq:Bregman_proximal_ineq} to the steps~\eqref{eq:interp_primal} and~\eqref{eq:ave_primal}, and obtain bounds on $g(x^{t+1}) - g(x)$ and $J(y^{t+1}) - J(y)$, i.e., 
\begin{align}
\hspace{-.0cm} J(y^{t+1}) - J(y) 
\le & -\lranglet{s^t}{y-y^{t+1}}+ \alpha_t^{-1}(D_{h_\calY}(y,y^t) - D_{h_\calY}(y,y^{t+1})) - ({2\alpha_t})^{-1}\normt{y^{t+1}-y^t}^2\nn\\
= &-\lranglet{(1+\theta_t)\hat{\nabla}_y\Phi(x^t,y^t,\zeta_y^t)-\theta_t\hat{\nabla}_y\Phi(x^{t-1},y^{t-1},\zeta_y^{t-1})}{y-y^{t+1}}\nn\\
& + \alpha_t^{-1}(D_{h_\calY}(y,y^t) - D_{h_\calY}(y,y^{t+1})) - ({2\alpha_t})^{-1}\normt{y^{t+1}-y^t}^2\label{eq:prox_J0}\\%\,\forall\,y\in\calY,
= & -\lranglet{(1+\theta_t)\delta^t_{y,\Phi}-\theta_t\delta_{y,\Phi}^{t-1}}{y-y^{t+1}} - \lranglet{{\nabla}_y\Phi(x^t,y^t)}{y-y^{t+1}}\nn\\
&-\theta_t\lranglet{{\nabla}_y\Phi(x^t,y^t)-{\nabla}_y\Phi(x^{t-1},y^{t-1})}{y-y^{t+1}}\nn\\
& + \alpha_t^{-1}(D_{h_\calY}(y,y^t) - D_{h_\calY}(y,y^{t+1})) - ({2\alpha_t})^{-1}\normt{y^{t+1}-y^t}^2,\label{eq:prox_J}\\
g(x^{t+1}) - g(x)
\le& \lranglet{\hat{\nabla}_x\Phi(x^t,y^{t+1},\zeta_x^t)+\hat{\nabla}f(\tilx^{t+1},\xi^t)}{x-x^{t+1}}\nn\\
 & + \tau_t^{-1}(D_{h_\calX}(x,x^t) - D_{h_\calX}(x,x^{t+1}))-({2\tau_t})^{-1}\normt{x^{t+1}-x^t}^2\\
 \le& \lranglet{\delta_{x,\Phi}^{t}+\delta_{x,f}^t}{x-x^{t+1}}+ \lranglet{{\nabla}_x\Phi(x^t,y^{t+1})+{\nabla}f(\tilx^{t+1})}{x-x^{t+1}}\nn\\
 & + \tau_t^{-1}(D_{h_\calX}(x,x^t) - D_{h_\calX}(x,x^{t+1}))-({2\tau_t})^{-1}\normt{x^{t+1}-x^t}^2,\label{eq:prox_g}
\end{align}
%where in~\eqref{eq:prox_J0} we use the fact that $s^0 = \hat{\nabla}_y\Phi(x^0,y^0,\zeta_y^0) = (1+\theta_0)\hat{\nabla}_y\Phi(x^0,y^0,\zeta_y^0)-\theta_0\hat{\nabla}_y\Phi(x^{-1},y^{-1},\zeta_y^{-1})$ 
Note that when $t=1$,~\eqref{eq:prox_J0} holds for any $\hat{\nabla}_y\Phi(x^{0},y^{0},\zeta_y^{0})\in\bbY^*$ since $\theta_1=0$. 

To bound $\Phi(x^{t+1},y)-\Phi(x,y^{t+1})$, we have
\begin{align}
\Phi(x^{t+1},y)-\Phi(x,y^{t+1})=& [\Phi(x^{t+1},y)-\Phi(x^{t+1},y^{t+1})]\nn\\
&+[\Phi(x^t,y^{t+1})-\Phi(x,y^{t+1})] + [\Phi(x^{t+1},y^{t+1}) - \Phi(x^t,y^{t+1})]\\
\le & \lranglet{\nabla_y \Phi(x^{t+1},y^{t+1})}{y-y^{t+1}} + \lranglet{\nabla_x\Phi(x^t,y^{t+1})}{x^t-x}\nn\\
& + \lranglet{\nabla_x\Phi(x^t,y^{t+1})}{x^{t+1}-x^t} + (L_{xx}/2)\normt{x^{t+1}-x^t}^2\label{eq:recur_Phi_interim}\\
\le & \lranglet{\nabla_y \Phi(x^{t+1},y^{t+1})}{y-y^{t+1}} + \lranglet{\nabla_x\Phi(x^t,y^{t+1})}{x^{t+1}-x}\nn\\
&+ (L_{xx}/2)\normt{x^{t+1}-x^t}^2,\label{eq:recur_Phi}
\end{align}
where in~\eqref{eq:recur_Phi_interim}, we use the concavity of $\Phi(x^{t+1},\cdot)$, the convexity of $\Phi(\cdot,y^{t+1})$ and the $L_{xx}$-smoothness of $\Phi(\cdot,y^{t+1})$, respectively. 
If we multiply both sides of~\eqref{eq:bound_pseudo_gap} by $\gamma_t\beta_t^{-1}$, and then substitute~\eqref{eq:prox_J},~\eqref{eq:prox_g} and~\eqref{eq:recur_Phi} into the resulting inequality, we have % to obtain 
\begin{align}
&\;\gamma_{t}\beta_{t}^{-1}\tilG(\barx^{t+1},\bary^{t+1};x,y)\le \gamma_t(\beta_t^{-1}-1)\tilG(\barx^{t},\bary^{t};x,y) + (\gamma_t/2)(L\beta_t+L_{xx}-1/\tau_t)\normt{x^{t+1}-x^t}^2\nn\\
&+ \gamma_t\tau_t^{-1}(D_{h_\calX}(x,x^t) - D_{h_\calX}(x,x^{t+1})) + \gamma_t\alpha_t^{-1}(D_{h_\calY}(y,y^t) - D_{h_\calY}(y,y^{t+1}))-{\gamma_t}({2\alpha_t})^{-1}\normt{y^{t+1}-y^t}^2\nn\\
& +\gamma_t\lranglet{\nabla_y \Phi(x^{t+1},y^{t+1})-\nabla_y \Phi(x^{t},y^{t})}{y-y^{t+1}}-\gamma_t\theta_t\lranglet{{\nabla}_y\Phi(x^t,y^t)-{\nabla}_y\Phi(x^{t-1},y^{t-1})}{y-y^{t}}\nn\\
&+\underbrace{\gamma_t\theta_t(-\lranglet{{\nabla}_y\Phi(x^t,y^t)-{\nabla}_y\Phi(x^{t-1},y^{t-1})}{y^t-y^{t+1}})}_{\rm (I)}+\underbrace{\gamma_t\lranglet{\delta_{x,\Phi}^{t}+\delta_{x,f}^t}{x-x^{t+1}}}_{\rm (II)}\nn\\
&+\underbrace{\gamma_t(1+\theta_t)\lranglet{\delta^t_{y,\Phi}}{y^{t+1}-y}}_{\rm (III)}+\underbrace{\gamma_t\theta_t\lranglet{\delta_{y,\Phi}^{t-1}}{y-y^{t+1}}}_{\rm (IV)}.\label{eq:recur_before_Exp}
\end{align}
 %Also, denote $\bbE[\cdot\,|\,\calF_{t}]$ by $\bbE_t[\cdot]$. 
%In addition, based on $\{\theta_t\}_{t\in\bbN}$, 

Before we proceed, recall that $\gamma_t\theta_t=\gamma_{t-1}$ and $\gamma_t\beta_t^{-1}=\gamma_{t+1}(\beta_{t+1}^{-1}-1)$, for any $t\in\bbN$. This enables us to observe certain recursion patterns (e.g., on $\gamma_t(\beta_t^{-1}-1)\tilG(\barx^{t},\bary^{t};x,y)$) in~\eqref{eq:recur_before_Exp}. We now bound the terms (I), (II) and (III) in~\eqref{eq:recur_before_Exp}. 

To bound (I), we make use of Young's inequality (cf.~\eqref{eq:Young}), i.e.,
\begin{align}
{\rm (I)}=&-\gamma_t\theta_t\lranglet{{\nabla}_y\Phi(x^t,y^t)-{\nabla}_y\Phi(x^{t},y^{t-1})}{y^t-y^{t+1}}-\gamma_t\theta_t\lranglet{{\nabla}_y\Phi(x^t,y^{t-1})-{\nabla}_y\Phi(x^{t-1},y^{t-1})}{y^t-y^{t+1}}\nn\\
\le &\gamma_t\theta_t\left\{(2L_{yy})^{-1}\normt{{\nabla}_y\Phi(x^t,y^t)-{\nabla}_y\Phi(x^{t},y^{t-1})}_*^2+(L_{yy}/2)\normt{y^{t+1}-y^t}^2\right\}\nn\\
& + \gamma_t\theta_t\left\{2\alpha_t\theta_t \normt{{\nabla}_y\Phi(x^t,y^{t-1})-{\nabla}_y\Phi(x^{t-1},y^{t-1})}_*^2 + (8\alpha_t\theta_t)^{-1}\normt{y^{t+1}-y^t}^2\right\}\nn\\
\lea &\gamma_t\theta_t\left\{(L_{yy}/2)\normt{y^t-y^{t-1}}^2+(L_{yy}/2)\normt{y^{t+1}-y^t}^2\right\}\nn\\
& + \gamma_t\theta_t\left\{2\alpha_t\theta_tL_{yx}^2 \normt{x^t-x^{t-1}}_*^2 + (8\alpha_t\theta_t)^{-1}\normt{y^{t+1}-y^t}^2\right\}\nn\\
 \leb &(\gamma_t/2)(\theta_tL_{yy}+(4\alpha_t)^{-1})\normt{y^{t+1}-y^t}^2 \nn\\
 &+ (\gamma_{t-1}L_{yy}/2)\normt{y^t-y^{t-1}}^2 +  2\gamma_{t-1}\alpha_{t-1}L_{yx}^2 \normt{x^t-x^{t-1}}^2,\label{eq:bound_I}
\end{align}
where in (a) we use the Lipschitz continuity of ${\nabla}_y\Phi(x^{t},\cdot)$ and ${\nabla}_y\Phi(\cdot,y^{t-1})$ respectively and in (b) we use the conditions that $\gamma_t\theta_t=\gamma_{t-1}$ and $\alpha_t\theta_t\le \alpha_{t-1}$ for any $t\in\bbN$. %(Again, when $t=0$, since $\theta_0=0$, $\alpha_{-1}$ can be any nonnegative number.)

%To bound~(II), we again employ Young's inequality, i.e.,
%\begin{align}
%{\rm (II)}\le & \gamma_t\tau_t 1\normt{\delta_{x,\Phi}^{t}+\delta_{x,f}^{t}}_*^2 + \gamma_t 1 (4\tau_t)^{-1}\normt{x^t-x^{t+1}}^2\nn\\
%& + \gamma_t\alpha_t 1\normt{(1+\theta)\delta_{y,\Phi}^{t}-\theta_t\delta_{y,\Phi}^{t-1}}_*^2 + \gamma_t 1 (4\alpha_t)^{-1}\normt{y^t-y^{t+1}}^2.
%\end{align}

To bound~(II), we need to use a technique introduced in~\citet{Nemi_09}, which involves the sequence $\{\hatx^t\}_{t\in\bbN}$ defined in~\eqref{eq:def_hat_x}. 
Based on $\{\hatx^t\}_{t\in\bbN}$, we decompose (II) into three parts, i.e.,
\begin{align}
{\rm (II)} =& \underbrace{\gamma_t\lranglet{\delta_{x,\Phi}^{t}+\delta_{x,f}^t}{x-\hatx^{t}}}_{\rm (II.A)} + \underbrace{\gamma_t\lranglet{\delta_{x,\Phi}^{t}+\delta_{x,f}^t}{\hatx^t-x^{t}}}_{\rm (II.B)} + \underbrace{\gamma_t\lranglet{\delta_{x,\Phi}^{t}+\delta_{x,f}^t}{x^t-x^{t+1}}}_{\rm (II.C)}.
\end{align}
To see the benefit of doing this, note that in (II.B), both $\hatx^t,x^t\in\calF_{t-1}$, i.e., $\hatx^t$ and $x^t$ are measurable with respect to $\calF_{t-1}$. Therefore by Assumption~\ref{assump:noise}\ref{assump:unbiased}, $\{\lranglet{\delta_{x,\Phi}^{t}+\delta_{x,f}^t}{\hatx^t-x^{t}}\}_{t\in\bbZ_+}$ is an MDS adapted to $\{\calF_t\}_{t\in\bbZ_+}$. Moreover, (II.A) and (II.C) can also be bounded using Corollary~\ref{cor:Bregman_proximal_inequality2} and Young's inequality respectively, i.e., 
\begin{align}
{\rm (II.A)} \le &\gamma_t\tau_t^{-1}(D_{h_\calX}(x,\hatx^t)-D_{h_\calX}(x,\hatx^{t+1})) + (\gamma_t\tau_t/2)\normt{\delta_{x,\Phi}^{t}+\delta_{x,f}^{t}}_*^2,\label{eq:bound_II.A}\\
{\rm (II.C)} \le  & \gamma_t\tau_t \normt{\delta_{x,\Phi}^{t}+\delta_{x,f}^{t}}_*^2 + \gamma_t  (4\tau_t)^{-1}\normt{x^t-x^{t+1}}^2. \label{eq:bound_II.C}
\end{align}
In summary, we have
\begin{align}
{\rm (II)} \le& \gamma_t\tau_t^{-1}(D_{h_\calX}(x,\hatx^t)-D_{h_\calX}(x,\hatx^{t+1})) + 2\gamma_t\tau_t \normt{\delta_{x,\Phi}^{t}+\delta_{x,f}^{t}}_*^2\nn\\
&+ \gamma_t  (4\tau_t)^{-1}\normt{x^t-x^{t+1}}^2+\gamma_t\lranglet{\delta_{x,\Phi}^{t}+\delta_{x,f}^t}{\hatx^t-x^{t}}.\label{eq:bound_II}
\end{align}

We can bound (III) and (IV) in a similar fashion. Indeed, define $\haty^0\defeq y^0$ and for any $t\in\bbN$, 
\begin{equation}
\haty^{t+1}\defeq \argmin_{y\in\calY}-\lranglet{\delta_{y,\Phi}^t}{y}+\alpha_t^{-1}D_{h_\calY}(y,\haty^t). 
\end{equation}
We then have
\begin{align}
{\rm (III)}=&\gamma_t(1+\theta_t)\lranglet{\delta^t_{y,\Phi}}{y^{t+1}-y^t} + \gamma_t(1+\theta_t)\lranglet{\delta^t_{y,\Phi}}{y^{t}-\haty^t}+\gamma_t(1+\theta_t)\lranglet{\delta^t_{y,\Phi}}{\haty^{t}-y}&\nn\\
\le& (1+\theta_t)((1+\theta_t)\gamma_t\alpha_t \normt{\delta^t_{y,\Phi}}_*^2+(1+\theta_t)^{-1}\gamma_t(4\alpha_t)^{-1}\normt{y^{t+1}-y^t}^2)\nn\\
& + (1+\theta_t)(\gamma_t\alpha_t^{-1}(D_{h_\calY}(y,\haty^t)-D_{h_\calY}(y,\haty^{t+1}))+(\gamma_t\alpha_t/2)\normt{\delta_{y,\Phi}^t}_*^2) + \gamma_t(1+\theta_t)\lranglet{\delta^t_{y,\Phi}}{y^{t}-\haty^t}\nn\\
\le& 5\gamma_t\alpha_t \normt{\delta^t_{y,\Phi}}_*^2+\gamma_t(4\alpha_t)^{-1}\normt{y^{t+1}-y^t}^2\nn\\
& + (1+\theta_t)\gamma_t\alpha_t^{-1}(D_{h_\calY}(y,\haty^t)-D_{h_\calY}(y,\haty^{t+1})) + \gamma_t(1+\theta_t)\lranglet{\delta^t_{y,\Phi}}{y^{t}-\haty^t},\label{eq:bound_III}
\end{align}
where %$\haty^0\defeq y^0$ and $\haty^{t+1}\defeq \argmin_{y\in\calY}-\lranglet{\delta_{y,\Phi}^t}{y}+\alpha_t^{-1}D_{h_\calY}(y,\haty^t)$, for any $t\in\bbZ_+$ and 
in the last inequality we use the fact that $\theta_t\in[0,1]$. In addition, %we bound (IV) as follows:
\begin{align}
{\rm (IV)}=& \gamma_{t-1}\lranglet{\delta_{y,\Phi}^{t-1}}{y-\hat{y}^{t-1}}+\gamma_{t-1}\lranglet{\delta_{y,\Phi}^{t-1}}{\hat{y}^{t-1}-y^{t-1}} + \gamma_{t-1}\lranglet{\delta_{y,\Phi}^{t-1}}{y^{t-1}-y^{t+1}} \nn\\
\le & \gamma_{t-1}(\alpha_{t-1}^{-1}(D_{h_\calY}(y,\haty^{t-1})-D_{h_\calY}(y,\haty^{t})) + (\alpha_{t-1}/2)\normt{\delta_{y,\Phi}^{t-1}}_*^2) \nn\\
& + \gamma_{t-1}(16\alpha_{t-1}\normt{\delta^{t-1}_{y,\Phi}}_*^2 + (64\alpha_{t-1})^{-1}\normt{y^{t-1}-y^{t+1}}^2) + \gamma_{t-1}\lranglet{\delta_{y,\Phi}^{t-1}}{\hat{y}^{t-1}-y^{t-1}}\nn\\
\le & \gamma_{t-1}\alpha_{t-1}^{-1}(D_{h_\calY}(y,\haty^{t-1})-D_{h_\calY}(y,\haty^{t})) + 17\gamma_{t-1}\alpha_{t-1}1\normt{\delta_{y,\Phi}^{t-1}}_*^2 + \gamma_{t-1}(32\alpha_{t-1})^{-1}\normt{y^{t-1}-y^{t}}^2\nn\\
& + \gamma_{t}(32\alpha_{t})^{-1}\normt{y^{t+1}-y^{t}}^2 + \gamma_{t-1}\lranglet{\delta_{y,\Phi}^{t-1}}{\hat{y}^{t-1}-y^{t-1}},\label{eq:bound_IV}
\end{align}
where in the last inequality we use $\normt{a+b}^2\le 2(\normt{a}^2+\normt{b}^2)$, for any $a,b\in\bbY$ and $\gamma_{t-1}/\alpha_{t-1}\le \gamma_t/\alpha_t$ (since $\alpha_t\theta_t\le \alpha_{t-1}$ and $\gamma_t\theta_t=\gamma_{t-1}$). 

Now we can substitute~\eqref{eq:bound_I},~\eqref{eq:bound_II},~\eqref{eq:bound_III} and~\eqref{eq:bound_IV} into the recursion~\eqref{eq:recur_before_Exp} to obtain 
\begin{align}
&\;\gamma_{t+1}(\beta_{t+1}^{-1}-1)\tilG(\barx^{t+1},\bary^{t+1};x,y)\le \gamma_t(\beta_t^{-1}-1)\tilG(\barx^{t},\bary^{t};x,y)\nn\\
& + (\gamma_t/2)(L\beta_t+L_{xx}-1/(2\tau_t))\normt{x^{t+1}-x^t}^2+  2\gamma_{t-1}\alpha_{t-1}L_{yx}^2 \normt{x^t-x^{t-1}}^2\nn\\
& + (\gamma_t/2)(\theta_tL_{yy}-3(16\alpha_t)^{-1})\normt{y^{t+1}-y^t}^2 + (\gamma_{t-1}/2)(L_{yy}+1/(16\alpha_{t-1}))\normt{y^t-y^{t-1}}^2 \nn\\
&+ \gamma_t\tau_t^{-1}\{(D_{h_\calX}(x,x^t)+D_{h_\calX}(x,\hatx^t)) - (D_{h_\calX}(x,x^{t+1})+D_{h_\calX}(x,\hatx^{t+1}))\}\nn\\
& + \gamma_t\alpha_t^{-1}(D_{h_\calY}(y,y^t) - D_{h_\calY}(y,y^{t+1}))+ (1+\theta_t)\gamma_t\alpha_t^{-1}(D_{h_\calY}(y,\haty^t)-D_{h_\calY}(y,\haty^{t+1}))\nn\\
& + \gamma_{t-1}\alpha_{t-1}^{-1}(D_{h_\calY}(y,\haty^{t-1})-D_{h_\calY}(y,\haty^{t}))\nn\\
&+ 2\gamma_t\tau_t \normt{\delta_{x,\Phi}^{t}+\delta_{x,f}^{t}}_*^2 + 5\gamma_t\alpha_t \normt{\delta^t_{y,\Phi}}_*^2 + 17\gamma_{t-1}\alpha_{t-1}\normt{\delta_{y,\Phi}^{t-1}}_*^2\nn\\
& +\gamma_t\lranglet{\nabla_y \Phi(x^{t+1},y^{t+1})-\nabla_y \Phi(x^{t},y^{t})}{y-y^{t+1}}-\gamma_{t-1}\lranglet{{\nabla}_y\Phi(x^t,y^t)-{\nabla}_y\Phi(x^{t-1},y^{t-1})}{y-y^{t}}\nn\\
&+\gamma_t\lranglet{\delta_{x,\Phi}^{t}+\delta_{x,f}^t}{\hatx^t-x^{t}}+ \gamma_t(1+\theta_t)\lranglet{\delta^t_{y,\Phi}}{y^{t}-\haty^t}+\gamma_{t-1}\lranglet{\delta_{y,\Phi}^{t-1}}{\hat{y}^{t-1}-y^{t-1}}.\label{eq:recur_before_Tel}
\end{align}

We  then sum up the inequality~\eqref{eq:recur_before_Tel} over $t=1,\ldots,T-1$ to obtain
\begin{align}
&\gamma_T(\beta_T^{-1}-1)\tilG(\barx^T,\bary^T;x,y) %- \gamma_0(\beta_0^{-1}-1)\tilG(\barx^T,\bary^T;x,y)
\nn\\
\le &\sum_{t=1}^{T-1}\frac{\gamma_t}{2}\left(L\beta_t+L_{xx}-({2\tau_t})^{-1}+{4\alpha_tL_{yx}^2}\right)\normt{x^{t+1}-x^t}^2 - {2\gamma_{T-1}\alpha_{T-1}L_{yx}^2}\normt{x^T-x^{T-1}}^2\nn\\
&+\sum_{t=1}^{T-1}\frac{\gamma_t}{2}\left((1+\theta_t)L_{yy}-\frac{1}{8\alpha_t}\right)\normt{y^{t+1}-y^t}^2-\frac{\gamma_{T-1}}{2}\left({L_{yy}+\frac{1}{16\alpha_{T-1}}}\right)\normt{y^{T}-y^{T-1}}^2\nn\\
&+ \sum_{t=1}^{T-1}\left(\frac{\gamma_t}{\tau_t}-\frac{\gamma_{t-1}}{\tau_{t-1}}\right)(D_{h_\calX}(x,x^t)+D_{h_\calX}(x,\hatx^t))% - \frac{\gamma_{T-1}}{\tau_{T-1}}(D_{h_\calX}(x,x^T)+D_{h_\calX}(x,\hatx^T))\nn\\
+ \sum_{t=1}^{T-1}\left(\frac{\gamma_t}{\alpha_t}-\frac{\gamma_{t-1}}{\alpha_{t-1}}\right)D_{h_\calY}(y,y^t)- \frac{\gamma_{T-1}}{\alpha_{T-1}}D_{h_\calY}(y,y^T)\nn\\
&+ \sum_{t=1}^{T-1}\left((2+\theta_t)\frac{\gamma_t}{\alpha_t}-(2+\theta_{t-1})\frac{\gamma_{t-1}}{\alpha_{t-1}}\right)D_{h_\calY}(y,\haty^t)%- (1+\theta_{T-1})\frac{\gamma_{T-1}}{\alpha_{T-1}}D_{h_\calY}(y,\haty^T)\nn\\
+ {2}\sum_{t=1}^{T-1}\gamma_t\tau_t \normt{\delta_{x,\Phi}^{t}+\delta_{x,f}^{t}}_*^2 + {22}\sum_{t=1}^{T-1}\gamma_t\alpha_t \normt{\delta^t_{y,\Phi}}_*^2 \nn\\
%+ \frac{17}{1}\sum_{t=1}^{T-2}\gamma_{t}\alpha_{t}\normt{\delta_{y,\Phi}^{t}}_*^2
&+\underbrace{\gamma_{T-1}\lranglet{\nabla_y \Phi(x^{T},y^{T})-\nabla_y \Phi(x^{T-1},y^{T-1})}{y-y^{T}}}_{\rm (V)}\nn\\
&+\sum_{t=1}^{T-1}\gamma_t\lranglet{\delta_{x,\Phi}^{t}+\delta_{x,f}^t}{\hatx^t-x^{t}}+ \sum_{t=1}^{T-2}\gamma_t\theta_t\lranglet{\delta^t_{y,\Phi}}{y^{t}-\haty^t}+\gamma_{T-1}(1+\theta_{T-1})\lranglet{\delta_{y,\Phi}^{T-1}}{{y}^{T-1}-\haty^{T-1}},
\label{eq:telescoped}
\end{align}
where we have used the facts that %$\gamma_{-1}=\gamma_0=0$ and 
$\gamma_{t-1}/\alpha_{t-1}\le \gamma_t/\alpha_t$. 

In addition, we can bound (V) in a similar fashion to bounding (I) (cf.~\eqref{eq:bound_I}), i.e.,
\begin{align}
{\rm (V)}=&\gamma_{T-1}\{\lranglet{\nabla_y \Phi(x^{T},y^{T})-\nabla_y \Phi(x^{T},y^{T-1})}{y-y^{T}} + \lranglet{\nabla_y \Phi(x^{T},y^{T-1})-\nabla_y \Phi(x^{T-1},y^{T-1})}{y-y^{T}}\}\nn\\
\le & \gamma_{T-1}\{\alpha_{T-1}L_{yy}^{2}\normt{y^{T}-y^{T-1}}^2+ (4\alpha_{T-1})^{-1}\normt{y-y^T}^2\nn\\
&\hspace{3cm} + \alpha_{T-1}L_{yx}^2\normt{x^T-x^{T-1}}^2 + (4\alpha_{T-1})^{-1}\normt{y-y^T}^2\}\nn\\
\le & (\gamma_{T-1}L_{yy}/2)\normt{y^{T}-y^{T-1}}^2 + \gamma_{T-1}\alpha_{T-1}L_{yx}^2\normt{x^T-x^{T-1}}^2 + \gamma_{T-1}\alpha_{T-1}^{-1}D_{h_\calY}(y,y^T),\label{eq:bound_inner_product_T}
\end{align}
where in the last inequality we have use $\alpha_{T-1}\le (2L_{yy})^{-1}$ and $(1/2)\normt{y-y^T}^2\le D_{h_\calY}(y,y^T)$. 

%Although~\eqref{eq:telescoped} looks complicated, it can be greatly simplified 
We now substitute~\eqref{eq:bound_inner_product_T} into~\eqref{eq:telescoped} and simplify the resulting inequality by noting that
\begin{enumerate}[label={\arabic*)}]
\item $L\beta_t+L_{xx}-({2\tau_t})^{-1}+{4\alpha_tL_{yx}^2}\le 0$, $(1+\theta_t)L_{yy}-({8\alpha_t})^{-1}\le 0$.
\item $D_{h_\calX}(x,x^T),D_{h_\calX}(x,\hatx^T)\le \Omega_{h_\calX}$, $\gamma_{t-1}/\alpha_{t-1}\le \gamma_t/\alpha_t$, \\
$(2+\theta_{t-1})\gamma_{t-1}/\alpha_{t-1}\le (2+\theta_{t})\gamma_t/\alpha_t$ (since $\theta_{t-1}\le \theta_t$), $D_{h_\calY}(y,y^T),D_{h_\calY}(y,\haty^T)\le \Omega_{h_\calY}$. 
\end{enumerate}
As a result, we have
\begin{align}
&\gamma_T(\beta_T^{-1}-1)\tilG(\barx^T,\bary^T;x,y) \le \sum_{t=1}^{T-1}\left(\frac{\gamma_t}{\tau_t}-\frac{\gamma_{t-1}}{\tau_{t-1}}\right)(D_{h_\calX}(x,x^t)+D_{h_\calX}(x,\hatx^t))+\frac{4\gamma_{T-1}}{\alpha_{T-1}}\Omega_{h_\calY}\nn\\
& + {2}\sum_{t=1}^{T-1}\gamma_t\tau_t \normt{\delta_{x,\Phi}^{t}+\delta_{x,f}^{t}}_*^2 + {22}\sum_{t=1}^{T-1}\gamma_t\alpha_t \normt{\delta^t_{y,\Phi}}_*^2\nn\\
&+\sum_{t=1}^{T-1}\gamma_t\lranglet{\delta_{x,\Phi}^{t}+\delta_{x,f}^t}{\hatx^t-x^{t}}+ \sum_{t=1}^{T-2}\gamma_{t-1}\lranglet{\delta^t_{y,\Phi}}{y^{t}-\haty^t}+\gamma_{T-1}(1+\theta_{T-1})\lranglet{\delta_{y,\Phi}^{T-1}}{{y}^{T-1}-\haty^{T-1}}.\label{eq:final_ineq_before_sup}
\end{align}
In addition, since $\gamma_{t-1}/\tau_{t-1}\le \gamma_t/\tau_t$, we have 
\begin{align}
\sum_{t=1}^{T-1}\left(\frac{\gamma_t}{\tau_t}-\frac{\gamma_{t-1}}{\tau_{t-1}}\right)(D_{h_\calX}(x,x^t)+D_{h_\calX}(x,\hatx^t))\le \frac{2\gamma_{T-1}}{\tau_{T-1}}\Omega_{h_\calX}. \label{eq:bound_sum_x}
\end{align}
%If we substitute , note that~\eqref{eq:final_ineq_before_sup} holds for any $(x,y)\in\calX\times\calY$, so we can replace the LHS of~\eqref{eq:final_ineq_before_sup} by  $\gamma_T(\beta_T^{-1}-1)G(\barx^T,\bary^T)$. As a result,~\eqref{eq:final_ineq_before_sup} becomes an almost sure bound on $G(\barx^T,\bary^T)$. %the duality gap at time $T$. 

We obtain the proof of parts~\ref{item:conv_exp} and~\ref{item:conv_whp} by treating the inequality~\eqref{eq:final_ineq_before_sup} in different ways. For part~\ref{item:conv_exp}, we simply take expectation on both sides of \eqref{eq:final_ineq_before_sup}; whereas for part~\ref{item:conv_whp}, we need to use concentration inequalities and the Chernoff bound. The details are shown below.

\emph{{Proof of Part~\ref{item:conv_exp}.}} For any $t\in\bbN$, %denote $\bbE[\cdot\,|\,\calF_{t}]$ by $\bbE_t[\cdot]$. 
since $x^t,\hatx^t,y^t,\haty^t\in\calF_{t-1}$, both $\{\lranglet{\delta_{x,\Phi}^{t}+\delta_{x,f}^t}{\hatx^t-x^{t}}\}_{t\in\bbN}$ and $\{\lranglet{\delta^t_{y,\Phi}}{y^{t}-\haty^t}\}_{t\in\bbN}$ are MDS adapted to $\{\calF_t\}_{t\in\bbZ_+}$. Therefore,
\begin{align}
\bbE[\lranglet{\delta_{x,\Phi}^{t}+\delta_{x,f}^t}{\hatx^t-x^{t}}] = \bbE[\lranglet{\delta^t_{y,\Phi}}{y^{t}-\haty^t}]=0.\label{eq:unbiased}
\end{align}
In addition, by Assumption~\ref{assump:noise}\ref{assump:variance_bound}, we have
\begin{align}
&\bbE[\normt{\delta_{y,\Phi}^{t}}_*^2]=\bbE[\bbE_{t-1}[\normt{\delta_{y,\Phi}^{t}}_*^2]]\le \sigma_{y,\Phi}^2, \label{eq:var_bounded1}\\ 
&\bbE[\normt{\delta_{x,f}^{t+1}+\delta_{x,\Phi}^{t+1}}_*^2]=\bbE[\bbE_{t-1}[\normt{\delta_{x,f}^{t+1}+\delta_{x,\Phi}^{t+1}}_*^2]]\le 2(\sigma_{x,f}^2+\sigma_{x,\Phi}^2),\label{eq:var_bounded2}
\end{align}
We then take expectation on both sides of~\eqref{eq:final_ineq_before_sup} and substitute~\eqref{eq:bound_sum_x},~\eqref{eq:unbiased},~\eqref{eq:var_bounded1} and~\eqref{eq:var_bounded2} to the resulting inequality to obtain~\eqref{eq:recur_res_exp}.
In addition, by Markov's inequality, for any $\varsigma\in(0,1]$,
\begin{align*}
\Pr\big\{G(\barx^T,\bary^T)>\varsigma^{-1}({\beta_T^{-1}-1})^{-1}(B_1(T) + B_2(T))\big\}\le \frac{\bbE[G(\barx^T,\bary^T)]}{\varsigma^{-1}({\beta_T^{-1}-1})^{-1}(B_1(T) + B_2(T))} \le \varsigma,
\end{align*}
where the last step follows from~\eqref{eq:recur_res_exp}. This proves~\eqref{eq:conv_markov}. 

{\em Proof of Part~\ref{item:conv_whp}.} %For completeness, 
We first present the Asuma-Hoeffding lemma for sub-Gaussian MDS. 
%The version that we will present  %which can be regraded as Asuma-Hoeffding with relaxed conditions. Its proof 
%can be found in~\citet{Nemi_09}.

\begin{lemma}[{\citet{Nemi_09}}]\label{lem:Asuma}
Let $\{\epsilon_t\}_{t\in\bbN}$ be a real-valued MDS adapted to a filtration $\{\calF_t\}_{t\in\bbZ_+}$, such that for any $t\in\bbN$, $\bbE_{t-1}[\epsilon_t]=0$ and there exists a constant $d_t>0$ such that $\bbE_{t-1}[\epsilon_t^2/d_t^2]\le \exp(1)$. Then for any $p>0$ and $T\in\bbN$, 
\begin{equation}
\Pr\left\{\textstyle{\sum}_{t=1}^{T-1}\epsilon_t>p\left({\sum_{t=1}^{T-1} d_t^2}\right)^{1/2}\right\}\le \exp(-p^2/4). 
\end{equation}
\end{lemma}

For convenience, define $C_T\defeq \sum_{t=1}^{T-1}\gamma_t\alpha_t$ and $C'_T\defeq \sum_{t=1}^{T-1}\gamma_t\tau_t$. Then 
\begin{align}
\Pr\left\{\sum_{t=1}^{T-1}\gamma_t\alpha_t\normt{\delta_{y,\Phi}^{t}}_*^2>(1+p)C_T\sigma_{y,\Phi}^2\right\} &= \Pr\left\{\exp\left(\frac{1}{C_T}\sum_{t=1}^{T-1}\gamma_t\tau_t\frac{\normt{\delta_{y,*}^{t}}_*^2}{\sigma_{y,\Phi}^2}\right)>\exp(1+p)\right\}\nn\\
&\lea \exp(-1-p) \bbE\left[\exp\left(\frac{1}{C_T}\sum_{t=1}^{T-1}\gamma_t\tau_t\frac{\normt{\delta_{y,*}^{t}}_*^2}{\sigma_{y,\Phi}^2}\right)\right]\nn\\
&\leb \exp(-1-p) \frac{1}{C_T}\sum_{t=1}^{T-1}\gamma_t\tau_t\bbE\left[\exp\left({\normt{\delta_{y,\Phi}^{t}}_*^2}/{\sigma_{y,\Phi}^2}\right)\right]\nn\\
&\lec \exp(-p), %C_T^{-1}\sum_{t=0}^{T-1}c_t \exp(\theta)\le \exp(-\theta^2), 
\label{eq:bound_var_y}
\end{align}
where in (a) we use Markov's inequality, in~(b) we use the convexity of $\exp(\cdot)$ and in~(c) we use Assumption~\ref{assump:noise}\ref{assump:sub-Gaussian}. 
Similarly, we can also show that for any $\tilp>0$, 
\begin{align}
&\Pr\left\{\sum_{t=1}^{T-1}\gamma_t\alpha_t\normt{\delta_{x,\Phi}^{t}}_*^2>(1+\tilp)C'_T\sigma_{x,\Phi}^2\right\} \le \exp(-\tilp), \label{eq:bound_var_x_1}\\
&\Pr\left\{\sum_{t=1}^{T-1}\gamma_t\alpha_t\normt{\delta_{x,f}^{t}}_*^2>(1+\tilp)C'_T\sigma_{x,f}^2\right\} \le \exp(-\tilp).  \label{eq:bound_var_x_2}
\end{align}

Next, %from~\eqref{eq:unbiased}, we observe that and 
since $\{\gamma_t\lranglet{\delta_{x,f}^{t}}{x^{t}- \hatx^t}\}_{t\in\bbN}$, $\{\gamma_t\lranglet{\delta_{x,\Phi}^{t}}{x^{t}- \hatx^t}\}_{t\in\bbN}$ and $\{\gamma_{t-1}\lranglet{\delta_{y,\Phi}^{t}}{\haty^t-y^{t}}\}_{t\in\bbN}$ are MDS, we can apply Lemma~\ref{lem:Asuma} to the last three terms in~\eqref{eq:final_ineq_before_sup}.
Specifically, let us define %$\tilde{gamma}_t=\gamma_t$, for $t=1,\ldots,T-3$ and $\tilde{\gamma}_{T-2} = \gamma_{T-1}(1+\theta_{T-1})$. According, we define 
\begin{align}
d_{t,y}^2\defeq 2\gamma_{t-1}^2\sigma_{y,\Phi}^2 D^2_{\calY}, \quad d_{t,x,\Phi}^2\defeq  \gamma_t^2\sigma_{x,\Phi}^2 D^2_{\calX}, \quad  d_{t,x,f}^2\defeq  \gamma_t^2\sigma_{x,f}^2 D^2_{\calX}.
\end{align}
Then by Assumption~\ref{assump:noise}\ref{assump:sub-Gaussian}, for any $t=1,\ldots,T-2$, we have  
% In addition, since $\abst{\gamma_t\lranglet{\delta_{y,*}^{t+1}}{\haty^t-y^{t}}}^2\le \gamma_t^2\normt{\delta_{y,*}^{t+1}}_*^2 \normt{\haty^t-y^{t}}^2 \le 2\gamma_t^2\normt{\delta_{y,*}^{t+1}}_*^2 1\Omega_{h_\calY}^2$ and $\abst{\gamma_t\lranglet{\delta^{t+1}+\delta_{x,*}^{t+1}}{x^{t}- \hatx^t}}^2\le \gamma_t^2\normt{\delta^{t+1}+\delta_{x,*}^{t+1}}_{*}^2\normt{x^{t}- \hatx^t}^2\le 2\gamma_t^2\normt{\delta^{t+1}+\delta_{x,*}^{t+1}}_*^21\Omega_{h_\calX}^2 $, we have that 
\begin{align}
\bbE_{t-1}\left[\exp\left(\abst{\gamma_{t-1}\lranglet{\delta_{y,\Phi}^{t}}{\haty^t-y^{t}}}^2/d_{t,y}^2\right)\right]\le& \bbE_{t-1}\left[\exp\left(\gamma_{t-1}^2\normt{\delta_{y,\Phi}^{t}}_*^2\normt{\haty^t-y^{t}}^2/d_{t,y}^2\right)\right]\nn\\
\le &\bbE_{t-1}\left[\exp\left(\gamma_{t-1}^2\normt{\delta_{y,\Phi}^{t}}_*^2D^2_\calY/d_{t,y}^2\right)\right] \le \exp(1),
\end{align}
and 
\begin{align}
\bbE_{t-1}\left[\exp\left(\abs{(\gamma_{T-2}+\gamma_{T-1})\lranglet{\delta_{y,\Phi}^{T-1}}{\haty^{T-1}-y^{T-1}}}^2\big/\big(2(\gamma_{T-2}^2+\gamma_{T-1}^2)\sigma_{y,\Phi}^2 D^2_\calY\big)\right)\right]\le \exp(1). 
\end{align}
Similarly, for any $t=1,\ldots,T-1$, we have
\begin{align}
& \bbE_{t-1}\left[\exp\left(\abst{\gamma_t\lranglet{\delta_{x,\Phi}^{t}}{\hatx^t-x^{t}}}^2/d_{t,x,\Phi}^2\right)\right] \le \exp(1),\\
 &  \bbE_{t-1}\left[\exp\left(\abst{\gamma_t\lranglet{\delta_{x,f}^{t}}{\hatx^t-x^{t}}}^2/d_{t,x,f}^2\right)\right] \le \exp(1).  
\end{align}
%from Assumption~\ref{assump:noise}\ref{assump:variance_LDbound}. 
Thus by Lemma~\ref{lem:Asuma}, for any $q,\tilq>0$, we have % to obtain that for any $\bar{\theta},\bar{\vartheta}>0$, 
\begin{align}
&\Pr\Bigg\{\sum_{t=1}^{T-2}\gamma_{t-1}\lranglet{\delta^t_{y,\Phi}}{y^{t}-\haty^t}+(\gamma_{T-2}+\gamma_{T-1})\lranglet{\delta_{y,\Phi}^{T-1}}{{y}^{T-1}-\haty^{T-1}}> {2q\sigma_{y,\Phi}D_\calY}\left(\sum_{t=1}^{T-1}\gamma_{t}^2\right)^{1/2}\Bigg\}\nn\\
&\hspace{12cm}\le \exp(-q^2/4),\label{eq:bound_inner_prod_y}\\
&\Pr\Bigg\{ \sum_{t=1}^{T-1}\gamma_t\lranglet{\delta_{x,\Phi}^{t}+\delta_{x,f}^t}{\hatx^t-x^{t}}> {{\tilq(\sigma_{x,\Phi}+\sigma_{x,f})D_\calX}}\left(\sum_{t=1}^{T-1}\gamma_{t}^2\right)^{1/2}\Bigg\}\le 2\exp(-\tilq^2/4).\label{eq:bound_inner_prod_x}
\end{align}
We then combine~\eqref{eq:final_ineq_before_sup},~\eqref{eq:bound_var_y},~\eqref{eq:bound_var_x_1},~\eqref{eq:bound_var_x_2},~\eqref{eq:bound_inner_prod_y} and~\eqref{eq:bound_inner_prod_x}, 
%substitute the choices of $\beta_t$, $\alpha_t$, $\tau_t$ and $\theta_t$ in Proposition~\ref{thm:main_cvx} 
and take $p=\tilp={\log(1/\varsigma)}$ and $q=\tilq=2\sqrt{\log(1/\varsigma)}$ to obtain~\eqref{eq:conv_whp_pseudo}. To obtain~\eqref{eq:recur_res_hp}, we simply substitute~\eqref{eq:bound_sum_x} into~\eqref{eq:conv_whp_pseudo}.
\Halmos\endproof

\endproof

%\newpage
%\section{Appendix.}
\begin{APPENDICES}
\section{Unique output of the Bregman proximal projection (BPP).}\label{app:unique}
\begin{lemma}\label{lem:unique_sln}
%In Problem~
For the minimization problem in~\eqref{eq:Bregman_projection}, if $\varphi^*\defeq\inf_{u\in\calU}\varphi(u)>-\infty$ and $\calU^o\cap\dom \varphi\ne \emptyset$, then a unique solution exists in $\calU^o\cap\dom \varphi$. 
\end{lemma}
\proof{Proof.}
%Denote the objective function in Problem~\eqref{eq:Bregman_projection} as $\varphi_\lambda:\bbU\to\barbbR$. 
%We first prove the existence of solutions.  
By condition~\eqref{eq:quad_lower_bound} and $\varphi^*>-\infty$, we see that $\varphi_\lambda$ is coercive on $\bbU$, i.e., $\lim_{\norm{u}\to+\infty}\varphi_\lambda(u)=+\infty$. In addition, $\varphi_\lambda$ is CCP since both $\varphi$ and $h_\calU$ are CCP. Choose any point $u\in\calU\cap\dom \varphi_\lambda\ne\emptyset$ %(note that $\dom \varphi_\lambda=\dom \varphi\cap\dom h$) 
and any $\alpha\ge \varphi_\lambda(u)$. The closedness and coercivity of $\varphi_\lambda$ together imply that the sub-level set $\calS_\alpha(\varphi_\lambda)\defeq \{u\in\bbU:\varphi_\lambda(u)\le \alpha\}$ is compact. Since $\calU$ is closed, $\calU\cap\calS_\alpha(\varphi_\lambda)$ is compact and nonempty. This, together with the lower semi-continuity of $\varphi_\lambda$, implies that the solution set of Problem~\eqref{eq:Bregman_projection}, denoted by $\calU_{\rm opt}$, is nonempty and contained in $\calU\cap\dom \varphi_\lambda$. %We next prove uniqueness. 
%In addition, 
Since $\inter\dom h_\calU\cap\dom \varphi\ne\emptyset$, by Moreau-Rockafellar theorem (see e.g.,~\citet[Theorem~3.30]{Peyp_15}), $\partial \varphi_\lambda(u) = \partial\varphi(u)+\partial h_\calU(u) + u^* - \nabla h(u')$, for any $u\in\calU$. Since $h_\calU$ is essentially smooth, for any $u\in\dom h_\calU\setminus\inter\dom h_\calU$, $\partial h_\calU(u)=\emptyset$, and hence $\partial \varphi_\lambda(u)=\emptyset$. As a result, $\calU_{\rm opt}\subseteq\calU^o\cap\dom\varphi$. Finally, since $h_\calU$ is 1-strongly convex on $\calU^o$, so is $\varphi_\lambda$, and therefore, $\calU_{\rm opt}$ must be a singleton. \Halmos
\endproof
\section{Berge's maximum theorem.}\label{app:Berge}
The following theorem is adapted from~\citet[Section~E.3]{Ok_07}, which was originally proved by~\citet{Berge_63}. 
\begin{theorem}\label{thm:Berge}
Let $\bbU$ and $\bbV$ be two metric spaces and $\Psi:\bbU\times\bbV\to\bbR$ be a continuous function on $\bbU\times\bbV$. If  $\calV\subseteq\bbV$ is nonempty and compact, then $\psi(u)\defeq \sup_{v\in\calV} \Psi(u,v)$ is continuous on $\bbU$. 
\end{theorem}
\section{Stopping criterion via gradient mappings.}\label{app:grad_map} % via Bregman Projection.}
Consider the BPP in~\eqref{eq:Bregman_projection} with $u^*=\nabla f(u')$, where the function $f$ is $L$-smooth and $\mu$-strongly convex on $\calU$. Accordingly, let us define two gradient mappings, i.e.,
\begin{equation}
G_\lambda (u')\defeq \frac{u'-u^+}{\lambda},\qquad \barG_\lambda (u')\defeq \frac{\nabla h_\calU(u')-\nabla h_\calU(u^+)}{\lambda}. 
\end{equation}
As we will show in the lemma below (i.e., Lemma~\ref{lem:grad_mapping}), the norm of the two gradient mappings, i.e., $G_\lambda(u')$ and $\barG_\lambda(u')$, can be used as the stopping criterion for the following optimization problem:
\begin{equation}
P^*\defeq \min_{u\in\bbU} \big\{P(u)\defeq f(u) + \varphi(u) + \iota_\calU(u)\big\}. \label{eq:grad_map_problem}
\end{equation}
Clearly, if $\dom\varphi\cap\calU\ne \emptyset$, then~\eqref{eq:grad_map_problem} has a unique optimal solution, which we denote by $u_{\rm opt}\in\calU$. Also, note that if $u_{\rm opt}\in\calU^o$, then $G_\lambda (u_{\rm opt})=\barG_\lambda (u_{\rm opt})=0$, since 
\begin{align*}
u_{\rm opt} = {\min}_{u\in\bbU} P(u) &\quad \Longleftrightarrow \quad 0\in \nabla f(u_{\rm opt}) + \partial (\varphi+\iota_\calU)(u_{\rm opt}) \\
&\quad \Longleftrightarrow \quad 0\in \nabla f(u_{\rm opt}) + \partial (\varphi+\iota_\calU)(u_{\rm opt}) + \frac{\nabla h_\calU(u_{\rm opt})-\nabla h_\calU(u_{\rm opt})}{\lambda}\\
&\quad \Longleftrightarrow \quad u_{\rm opt} = {\argmin}_{u\in\calU}\; \varphi(u) + \lrangle{\nabla f(u_{\rm opt})}{u} + \lambda^{-1}D_{h_\calU}(u,u_{\rm opt}),
\end{align*}
where in the last step, we use the uniqueness of the BPP (cf.\ Lemma~\ref{lem:unique_sln}). 

\begin{lemma}\label{lem:grad_mapping}
Consider the setting in Appendix~\ref{app:grad_map}. For any $\epsilon>0$ and $0<\lambda\le 1/L$, 
we have
\begin{equation}
\normt{\barG_\lambda (u')}_*^2 + \normt{G_\lambda (u')}^2\le \mu\epsilon \quad \Longrightarrow \quad P(u^+) - P^*\le \epsilon. 
\end{equation}
%if $\normt{\barG_\lambda (u')}_*^2 + \normt{G_\lambda (u')}^2\le \mu\epsilon$, then $P(u^+) - P^*\le \epsilon$. 
\end{lemma}

\proof{Proof.}
By the first-order optimality condition of~\eqref{eq:Bregman_projection}, we have 
\begin{equation}
0\in \partial( \varphi+\iota_\calU)(u^+) + \nabla f(u') - \barG_\lambda (u') \quad \Longleftrightarrow \quad \barG_\lambda (u') - \nabla f(u')\in \partial( \varphi+\iota_\calU)(u^+), 
\end{equation}
where $\calN_\calU(u^+)$ denotes the normal cone of $\calU$ at $u^+\in\calU$. Consequently, we have
\begin{equation}
\xi\defeq \barG_\lambda (u') + \nabla f(u^+)- \nabla f(u') \in \nabla f(u^+) + \partial( \varphi+\iota_\calU)(u^+)= \partial P(u^+). 
\end{equation}
The $\mu$-strong convexity of $P$ on $\calU$ indicates that %\supseteq\dom P$ indicates that for any $u\in\bbU$,
\begin{align}
P^* = P(u_{\rm opt}) & \ge P(u^+) + \lranglet{\xi}{u_{\rm opt}-u^+} + (\mu/2)\normt{u_{\rm opt}-u^+}^2\nn\\
&\ge P(u^+) -\normt{\xi}_*\norm{u_{\rm opt}-u^+} + (\mu/2)\normt{u_{\rm opt}-u^+}^2\nn\\
&\ge P(u^+) - \normt{\xi}_*^2/(2\mu). \label{eq:grad_map_bound1}
\end{align}
On the other hand, we have
\begin{align}
\normt{\xi}_*^2 \le 2\big(\normt{\barG_\lambda (u')}_*^2 + \normt{\nabla f(u^+)- \nabla f(u')}_*^2\big) &\le 2\big(\normt{\barG_\lambda (u')}_*^2 + L^2\normt{u^+ - u'}_*^2\big)\nn \\
&\le 2\big(\normt{\barG_\lambda (u')}_*^2 + \normt{G_\lambda (u')}^2\big), \label{eq:grad_map_bound2}
\end{align}
where in the last step we use that $0<\lambda\le 1/L$. Now, combining~\eqref{eq:grad_map_bound1} and~\eqref{eq:grad_map_bound2}, we have
\begin{equation}
P(u^+) - P^*\le \big(\normt{\barG_\lambda (u')}_*^2 + \normt{G_\lambda (u')}^2\big)/\mu. 
\end{equation}
We hence complete the proof. 
\Halmos
\endproof
%to bound the sub-optimality of 

\begin{remark}
Two remarks are in order. 
First, note that both $G_\lambda (u')$ and $\barG_\lambda (u')$ can be easily computed, since 1) we assume that the BPP in~\eqref{eq:Bregman_projection} has an easily computable solution (cf.\ Section~\ref{sec:DGF_BPP}), and 2) the DGF $h_\calU$ typically has simple structure, so that $\nabla h_\calU$ can be obtained easily. Second, if $\nabla h_\calU$ is $\eta$-Lipschitz on $\calU$ ($\eta\ge 1$), as in the Hilbertian case where $h_\calU = (1/2)\normt{\cdot}^2$, then $\normt{\barG_\lambda (u')}_*\le \eta \normt{G_\lambda (u')}$ and we need not compute $\nabla h_\calU$. 
\end{remark}

\section{Lipschitz continuity of $x^*(\cdot)$.}\label{app:Lips}
The following lemma is similar to~\citet[Lemma~3.4]{Zhao_20}. For completeness, we provide it proof here. 

\begin{lemma}\label{lem:lips_x^*}
Let $x^*(\cdot):\calY\to\bbX$ be defined in~\eqref{eq:def_x^*}. Then it is $(L_{yx}/\mu)$-Lipschitz on $\calY$. 
\end{lemma}

\proof{Proof.}
By the $\mu$-strong convexity of $f$ (hence $\psi^\rmP(\cdot,y)$, for any $y\in\calY$), we have
\begin{align}
\psi^\rmP(x^*(y'),y) - \psi^\rmP(x^*(y),y)\ge (\mu/2)\norm{x^*(y')-x^*(y)}^2,\label{eq:lb_1}\\
\psi^\rmP(x^*(y),y') - \psi^\rmP(x^*(y'),y')\ge (\mu/2)\norm{x^*(y')-x^*(y)}^2.\label{eq:lb_2}
\end{align}
On the other hand,
\begin{align}
\psi^\rmP(x^*(y'),y) - \psi^\rmP(x^*(y'),y') &= \int_0^1 \lranglet{\nabla_y\psi^\rmP(x^*(y'),y'+t(y-y'))}{y-y'}\; \rmd t,\label{eq:ub_1}\\
\psi^\rmP(x^*(y),y) - \psi^\rmP(x^*(y),y') &= \int_0^1 \lranglet{\nabla_y\psi^\rmP(x^*(y),y'+t(y-y'))}{y-y'}\; \rmd t. \label{eq:ub_2}
\end{align}
Combining~\eqref{eq:lb_1} to~\eqref{eq:ub_2}, we have
\begin{align*}
\mu \norm{x^*(y')-x^*(y)}^2 &\le \int_0^1 \lranglet{\nabla_y\psi^\rmP(x^*(y'),y'+t(y-y')) - \nabla_y\psi^\rmP(x^*(y),y'+t(y-y'))}{y-y'}\; \rmd t\\
&\le \int_0^1 \normt{\nabla_y\psi^\rmP(x^*(y'),y'+t(y-y')) - \nabla_y\psi^\rmP(x^*(y),y'+t(y-y'))}_*\normt{y-y'}\; \rmd t\\
&\le L_{yx} \normt{x^*(y') - x^*(y)} \normt{y-y'}. \nt\label{eq:final_bd}
\end{align*}
Note that~\eqref{eq:final_bd} implies that $\norm{x^*(y')-x^*(y)}\le (L_{yx}/\mu)\normt{y-y'}$, for both $x^*(y') = x^*(y)$ and $x^*(y') \ne x^*(y)$. We hence complete the proof. \Halmos %Specifically, it is clear if $x^*(y') \ne x^*(y)$, 
\endproof

\section{Technical proofs in Section~\ref{sec:experiment}.}\label{app:deriv_constants}
We first derive the constants $L_{xx}$, $L_{xy}$ and $L_{yy}$ in~\eqref{eq:constants_game}.  
Let $i,j\in[n]$ and $\delta_{ij}$ be the Kronecker delta, i.e., $\delta_{ij}=1$ if $i = j$ and $0$ otherwise. It is easy to see that 
\begin{align}
\frac{\partial}{\partial x_i} \Phi(x,y) &= -\frac{y_i}{(c_i+x_i)(c_i+x_i+y_i)},\\
 \frac{\partial}{\partial y_i} \Phi(x,y) &= (c_i+x_i+y_i)^{-1},\\
\frac{\partial^2}{\partial x_i\partial x_j} \Phi(x,y) &= \big((c_i + x_i)^{-2} - (c_i + x_i + y_i)^{-2}\big)\delta_{ij},\\
\frac{\partial^2}{\partial x_i\partial y_j} \Phi(x,y) = \frac{\partial^2}{\partial y_i\partial y_j} \Phi(x,y) &= -(c_i + x_i + y_i)^{-2} \delta_{ij}.
\end{align}
In addition, we have that for any $i\in[n]$,
\begin{align}
&(c_i + x_i)^{-2} - (c_i + x_i + y_i)^{-2} \le (c_i + x_i)^{-2} - (c_i + x_i + 1)^{-2}\lea  c_{\min}^{-2} - (c_{\min} + 1)^{-2},\\
&(c_i + x_i + y_i)^{-2} \le c_{\min}^{-2}. 
\end{align}
where in (a) we use $c_i + x_i\ge c_{\min}$ and that the function $a \mapsto a^{-2} - (a+1)^{-2}$ is strictly decreasing on $(0,+\infty)$. Therefore, we can let $L_{xx} = c_{\min}^{-2} - (c_{\min} + 1)^{-2}$ and $L_{xy} = L_{yy} = c_{\min}^{-2}$. 

Next, we derive the variance parameters in~\eqref{eq:var_game}. To do that, let us first present a well-known lemma on the variance of sampling without replacement from a finite population (see e.g.,~\citet{Banerjee_12}). 

\begin{lemma}\label{lem:FPC}
Given a finite set $\calZ\defeq \{z_i\}_{i=1}^n\subseteq \bbR^d$, uniformly randomly draw a subset $\calB\subseteq [n]$ without replacement. %, and denote them by $\{z_i\}_{i=1}^n$. 
Define the population mean $\barz \defeq n^{-1}\sum_{i=1}^{n} z_i$. Then 
\begin{equation}
\bbE\left[\norm{\abs{\calB}^{-1}\textstyle\sum_{i\in\calB} z_i - \barz}_2^2\right] \le \frac{n-\abs{\calB}}{\abs{\calB}(n-1)}\left(n^{-1}\textstyle\sum_{i=1}^n \norm{z_i}_2^2 - \norm{\barz}_2^2\right). 
\end{equation}
\end{lemma}
Based on Lemma~\ref{lem:FPC}, to derive $\sigma_{x,f}^2$, it suffices to upper bound $\norm{z_i}_2^2$ and lower bound $\norm{\barz}_2^2$, where $z_i = \barQ_i\barQ_i^T x$, and $x\in\Delta_n$.  By H\"older's inequality, it is clear that $\norm{z_i}_2^2 = \norm{\barQ_i}_2^2(\barQ_i^T x)^2\le \norm{\barQ_i}_2^2\norm{\barQ_i}_\infty^2$. In addition, since $\barz = n^{-1}\sum_{i=1}^n \barQ_i\barQ_i^T x = n^{-1} Qx$, we have $\norm{\barz}_2^2 = n^{-2}\norm{Qx}_2^2\ge n^{-2}\norm{x}_2^2\sigma^2_{\min}(Q)\ge n^{-3}\sigma_{\min}^2(Q)$. As a result, for any $(x,y)\in\Delta_n\times\Delta_n$, 
\begin{align}
\bbE\left[\normt{\hat{\nabla} f(x) - \nabla f(x)}_2^2\right] \le \frac{n^2 (n-\abs{\calB})}{\abs{\calB}(n-1)}\left(n^{-1}\textstyle\sum_{i=1}^n \norm{\barQ_i}_2^2\norm{\barQ_i}_\infty^2 - n^{-3}\sigma_{\min}^2(Q)\right). 
\end{align}
Using the same argument, we have
\begin{align}
\bbE\left[\normt{\hat{\nabla}_x \Phi(x,y) - {\nabla}_x \Phi(x,y)}_2^2\right] &\le \frac{n^2 (n-\abs{\calB})}{\abs{\calB}(n-1)},\\
\bbE\left[\normt{\hat{\nabla}_y \Phi(x,y) - {\nabla}_y \Phi(x,y)}_2^2\right] &\le \frac{n^2 (n-\abs{\calB})}{\abs{\calB}(n-1)}. 
\end{align}

\end{APPENDICES}

\section*{Acknowledgments.}
The author would like to thank Robert Freund for his constructive feedback and helpful discussions during the preparation of this manuscript. The author would also like to thank two anonymous referees for their helpful and insightful suggestions, which have greatly improved the exposition of the paper. The author's research is supported by AFOSR Grant No. FA9550-19-1-0240.

% Acknowledgments here

% References here (outcomment the appropriate case) 

% CASE 1: BiBTeX used to constantly update the references 
%   (while the paper is being written).
\bibliographystyle{informs2014} % outcomment this and next line in Case 1
\bibliography{math_opt,dataset,mach_learn,stat_ref,stoc_ref} % if more than one, comma separated

% CASE 2: BiBTeX used to generate mypaper.bbl (to be further fine tuned)
%\input{mypaper.bbl} % outcomment this line in Case 2

\end{document}